\documentclass{amsart}

\usepackage[a4paper, top=2cm, bottom=2cm, left=3cm, right=3cm, heightrounded, centering]{geometry}

\usepackage[backend=biber, style=alphabetic,maxbibnames=99,maxalphanames=99]{biblatex}
\usepackage{graphicx}
\usepackage{float}
\usepackage[T1]{fontenc}
\usepackage{amsmath,amsthm,amssymb,amscd,mathtools}
\usepackage{thmtools}
\usepackage{amsfonts}
\usepackage{lastpage}
\usepackage{fancyhdr}
\usepackage{mathrsfs}
\usepackage{tikz-cd}
\usepackage{microtype} 
\usepackage{stmaryrd}
\usepackage{xr} 
\usepackage[dvipsnames]{xcolor}
\usepackage[bookmarksopen,bookmarksdepth=2]{hyperref}
\hypersetup{colorlinks=true,citecolor=NavyBlue,linkcolor=NavyBlue,urlcolor=Green,hypertexnames=false} 
\newtheorem{thm}{Theorem}[section]
\newtheorem{prop}[thm]{Proposition}

\newtheorem{lem}[thm]{Lemma}
\newtheorem{cor}[thm]{Corollary}

\theoremstyle{remark}
\newtheorem{rmk}[thm]{Remark}
\newtheorem{eg}[thm]{Example}
\theoremstyle{definition}
\newtheorem{dfn}[thm]{Definition}
\newtheorem{notation}[thm]{Notation}
\newtheorem{construction}[thm]{Construction}

\usepackage{cleveref}
\crefname{thm}{Theorem}{Theorems}
\crefname{prop}{Proposition}{Propositions}
\crefname{fact}{Fact}{Facts}
\crefname{lem}{Lemma}{Lemmas}
\crefname{cor}{Corollary}{Corollaries}
\crefname{conj}{Conjecture}{Conjectures}
\crefname{question}{Question}{Questions}
\crefname{rmk}{Remark}{Remarks}
\crefname{eg}{Example}{Examples}
\crefname{dfn}{Definition}{Definitions}
\crefname{notation}{Notation}{Notations}
\crefname{construction}{Construction}{Constructions}

\Crefname{thm}{Theorem}{Theorems}
\Crefname{prop}{Proposition}{Propositions}
\Crefname{fact}{Fact}{Facts}
\Crefname{lem}{Lemma}{Lemmas}
\Crefname{cor}{Corollary}{Corollaries}
\Crefname{conj}{Conjecture}{Conjectures}
\Crefname{question}{Question}{Questions}
\Crefname{rmk}{Remark}{Remarks}
\Crefname{eg}{Example}{Examples}
\Crefname{dfn}{Definition}{Definitions}
\Crefname{notation}{Notation}{Notations}
\Crefname{construction}{Construction}{Constructions}

\newcommand{\mbb}[1]{\mathbb{#1}}
\newcommand{\mcal}[1]{\mathcal{#1}}

\newcommand{\NN}{\mathbb{N}}
\newcommand{\RR}{\mathbb{R}}
\newcommand{\ZZ}{\mathbb{Z}}
\newcommand{\QQ}{\mathbb{Q}}
\newcommand{\CC}{\mathbb{C}}

\newcommand{\FF}{\mathbb{F}}
\newcommand{\GG}{\mathbb{G}}

\newcommand{\Gal}{\mathrm{Gal}}

\DeclareMathOperator{\im}{Im}

\DeclareMathOperator{\Hom}{Hom}
\DeclareMathOperator{\Ext}{Ext}

\DeclareMathOperator{\Nm}{Nm}

\DeclareMathOperator{\Spf}{Spf}
\DeclareMathOperator{\Spec}{Spec}

\DeclareMathOperator{\Lie}{Lie}
\DeclareMathOperator{\Mod}{Mod}
\newcommand{\Res}{\mathrm{Res}}

\newcommand{\GL}{\mathrm{GL}}

\newcommand{\mi}{\mathrm{i}}

\DeclareMathOperator{\Fil}{Fil}

\newcommand{\p}{\mathfrak{p}}

\DeclareMathOperator{\pr}{pr}
\newcommand{\mO}{\mathcal{O}}
\newcommand{\mP}{\mathbb{P}}
\newcommand{\mA}{\mathbb{A}}

\DeclareMathOperator{\SL}{SL}

\DeclareMathOperator{\End}{End}

\newcommand{\id}{\mathrm{id}}

\newcommand{\mF}{\mathcal{F}}

\DeclareMathOperator{\Ker}{Ker}

\newcommand{\mX}{\mathcal{X}}

\newcommand{\et}{\mathrm{\acute{e}t}}

\newcommand{\proket}{\mathrm{prok\acute{e}t}}

\newcommand{\dR}{\mathrm{dR}}

\newcommand{\OBdR}{\mO\mbb{B}_\dR}

\newcommand{\OBdRKp}{\mO\mbb{B}_{\dR,K^{p}}}

\newcommand{\BdR}{\mbb{B}_\dR}

\newcommand{\BdRKp}{\mbb{B}_{\dR,K^{p}}}
\newcommand{\Bdif}{\mbb{B}_{\dif}}
\newcommand{\Bdifcf}{\mbb{B}_{\dif,L}}

\newcommand{\BdifKpcf}{\mbb{B}_{\mrm{dif},K^{p},L}}

\DeclareMathOperator{\gr}{gr}
\newcommand{\HT}{\mathrm{HT}}
\newcommand{\GM}{\mathrm{GM}}

\newcommand{\Fl}{{\mathscr{F}\!l}}

\newcommand{\TT}{\mbb{T}}

\newcommand{\fg}{\mathfrak{g}}
\newcommand{\sm}{{\mathrm{sm}}}
\newcommand{\la}{{\mathrm{la}}}
\newcommand{\an}{{\mathrm{an}}}

\newcommand{\Kpp}{K^pK_p}

\newcommand{\mrm}[1]{\mathrm{#1}}
\newcommand{\mfk}[1]{\mathfrak{#1}}
\newcommand{\mbf}[1]{\mathbf{#1}}

\newcommand{\RHom}{R\!\Hom}

\newcommand{\mscr}[1]{\mathscr{#1}}

\newcommand{\dif}{\mrm{dif}}

\newcommand{\Rlim}{R\!\varprojlim}
\newcommand{\Fib}{\mrm{Fib}}

\newcommand{\oInd}{\otimes\!-\!\mrm{Ind}}

\newcommand{\kw}{(\Bbbk,w)}

\DeclareMathOperator{\Fr}{Fr}
\newcommand{\QCoh}{\mrm{QCoh}}

\newcommand{\tor}{\mrm{tor}}
\DeclareMathOperator{\Rep}{Rep}

\DeclareMathOperator{\AnSpec}{AnSpec}

\DeclareMathOperator{\oIndFQ}{\otimes-\mrm{Ind}^{\QQ}_{F}}

\usepackage{mathtools}

\newcommand{\mS}{\mcal{S}}
\newcommand{\ti}[1]{\tilde{#1}}

\newcommand{\SGTKp}[1]{\mS_{K^{p}}(G_{#1})}

\newcommand{\SGTKpTo}[1]{\mS_{K^{p}}^{\tor}(G_{#1})}

\newcommand{\Darith}{D_{\mrm{arith}}}

\newcommand{\pdR}{\mrm{alg}.\dR}

\newcommand{\hatotimes}{\mathbin{\hat{\otimes}}}
\newcommand{\unl}[1]{\underline{#1}}
\newcommand{\std}{\mrm{std}}
\newcommand{\Per}{\mscr{P}\!er}

\newcommand{\cycl}{\mrm{cycl}}

\newcommand{\Vin}{\mrm{Vin}}

\newcommand{\cf}{L} 
\newcommand{\cfi}{\cf_{\infty}}
\newcommand{\Darithcf}{D_{\mrm{arith},\cf}}

\newcommand{\rk}{\mrm{rk}}
\newcommand{\Eu}{\mrm{Eu}}
\newcommand{\loc}{\mrm{loc}}

\begin{document}
\title{Fontaine operators and Vinberg semigroups}
\author{Yuanyang Jiang}
\address{Max Planck Institute for Mathematics, Bonn}
\email{jiang@mpim-bonn.mpg.de}

\begin{abstract}
For any Shimura datum, we construct a \( B_{\dR}^{+} \)-thickening of the locally analytic Shimura variety constructed in \cite{Jiang2026Shla}, for which we construct a local model based on Vinberg semigroups. This reduces the study of geometric Fontaine operators considered in \cite{Pan2209.06II} to a calculation on Vinberg semigroups. 

As an application, let \(F\) be a totally real number field,
and consider the associated Hilbert Shimura datum. We construct a certain plectic structure on the period sheaf using the local model. This allows us to extend our previous results in \cite{Jiang2025HMF} on classicality and on the Langlands-Clozel-Fontaine-Mazur conjecture to the cases of possibly non-parallel weights.
\end{abstract}

\maketitle
\tableofcontents

\section{Introduction}
In this note, we will show that  certain subvarieties of Vinberg semigroups (\cite{Vinberg1995reductive}, \cite{Vinberg1995asymptotic}) are local models of the \( B_{\dR}^{+} \)-thickening of the locally analytic Shimura varieties as formulated in \cite{Jiang2026Shla}. 
This gives a new perspective on the calculation of Fontaine operators in \cite{Pan2209.06II}. 

As concrete applications, we will 
\begin{enumerate}
\item 
give an alternative proof of the identity \( N=d\circ \bar{d} \) of \cite{Pan2209.06II} (\S \ref{subsec:Fontaine-op-mod-curve}); 
\item extend the classicality theorem of \cite{Jiang2025HMF} to the cases of non-parallel weights. 
\end{enumerate}

More precisely, let \(F\) be a totally real number field of dimension \(d\). Let \(L\subset\bar{\QQ}_{p}\) be a finite Galois extension of \(\QQ_{p}\), with ring of integer \( \mO_{L} \) and residue field \( \FF \). 

Let \(G:=\mrm{Res}_{F/\QQ}(\GL_{2})\). For any neat open compact subgroup \(K\subset G(\mA_{f})\), we have the \(d\)-dimensional Hilbert modular variety \(X_{K}\) of level \(K\) defined over \(\QQ\).
Following \cite{Emerton06}, we fix a neat open compact subgroup \(K^{p}\subset G(\mA_{f}^{p})\) and define the \(i\)-th \emph{completed cohomology} of the Hilbert modular varieties as
\[
\ti{H}^{i}(K^{p},L):=\left(\varprojlim_{n}\varinjlim_{K_{p}}H^{i}_{\et}(X_{\Kpp,\bar{\QQ}},\ZZ/p^{n})\right)\hatotimes_{\ZZ_{p}}L,
\]
where the colimit is taken over all open compact subgroups \(K_{p}\subset G(\QQ_{p})\). 

Let \(S\) be a fixed finite set of places of \(\QQ\), including \(p\) and \(\infty\), as well as all places that ramify in \( F \) or \(K^{p}\). Denote by \(\TT^{S}\) the spherical Hecke algebra over \(\ZZ_{p}\) generated by the spherical Hecke operators of \(G\) at finite places away from \(S\). Then \(\ti{H}^{i}(K^{p},L)\) is equipped with an action of \(\TT^{S}\times G(\QQ_{p})\times \Gal_{\QQ}\). We will prove the following \emph{classicality} theorem, removing the parallel weight condition in \cite[Theorem \ref{2-IntrothmClassicalityCompleteCoho}]{Jiang2025HMF}.
\begin{thm}[\cref{thm:ClassicalityCompleteCoho}]
	\label{Intro:thm:ClassicalityCompleteCoho}
	Let \(\rho:\Gal_{F,S}\to \GL_{2}(L)\) be a continuous representation.
	Let \(\chi:\TT^{S} \to L\) be the associated character determined by the Eichler-Shimura relation (\ref{alignEicherShimuraDetermin}).


	We assume that 
	
	(1) \(\ti{H}^{d}(K^{p},L)[\chi]:=\Hom_{\TT^{S}}(\chi,\ti{H}^{d}(K^{p},L))\ne 0\);

	(2) \(\rho\) is absolutely irreducible;
	
	(3) \(\rho\) is de Rham of regular Hodge-Tate weights.

	Then the character \(\chi\) is classical, i.e. there exists a Hilbert modular eigenform \(f\),
	such that \(\TT^{S}\) acts on \(f\) via \(\chi\). In other words, \(\rho\cong \rho_{f}\).
\end{thm}
As an application, we can prove the following new cases of Langlands-Clozel-Fontaine-Mazur conjecture, removing the parallel weight condition in \cite[\cref{2-thmConjecFontaineMazur}]{Jiang2025HMF}.
\begin{thm}[\cref{thmConjecFontaineMazur}]\label{Intro:thmConjecFontaineMazur}
Assume that \(p\) is inert in \(F\).
Let 
\(\rho:\Gal_{F,S}\to \GL_{2}(L)\) be a continuous representation.
Let \(\bar{\rho}:\Gal_{F}\to \GL_{2}(\FF)\) be the reduction of \(\rho\). 
Let \(\chi:\TT^{S}\to \mO_{L}\) be the character determined by \(\rho\) via (\ref{alignEicherShimuraDetermin}).

Assume that 

(1) \(\bar{\rho}\) is residual modular, that is, there exists a Hilbert modular eigenform \(f_{0}\)
such that \(\bar{\rho}\cong \bar{\rho}_{f_{0}}\).

(2) \(\bar{\rho}|_{\Gal_{F(\sqrt[p]{1})}}\)
is absolutely irreducible.

(3) \(\bar{\rho}|_{\Gal_{F_{p}}}\) is semisimple generic in the sense of \cite[p.100 (i)]{BreuilHerzigHuMorraSchraen2023GKdim}.

(4) For any \(w\in S\backslash \{p\}\),
the universal framed deformation ring of \(\bar{\rho}|_{\Gal_{F_{w}}}\) over \(W(\FF)\)
is formally smooth over \(W(\FF)\).

(5) Let \(\psi\) denote the Teichmüller lift of \(\chi_{\mrm{cycl}}\det(\bar{\rho})\), and we 
assume that \(\det(\rho)=\psi \chi_{\mrm{cycl}}^{-1}\).



(6) \(\rho\) is de Rham of regular weights.

Then there exists a Hilbert modular eigenform \(f\)
such that \(\rho\cong \rho_{f}\).
\end{thm}

The proof of \cref{Intro:thm:ClassicalityCompleteCoho} follows the strategy of \cite{Pan2209.06II}. We refer readers to the introduction of \cite{Jiang2025HMF} for a summary of the strategy in the Hilbert case. Roughly, by a classical theorem of \cite{Fontaine2004arithmetique}, we can understand the de Rhamness of a  Galois representation of finite dimension through a certain \emph{Fontaine operator}. Concretely, for a Hodge representation \( \rho \) of \( \Gal_{K} \) (where \( K\subset \cf \)), one can consider \[
D_{\mrm{pdR},\cf}(\rho):=E_{0}^{\theta}(\Darithcf(\rho\otimes B_{\dR})), 
\] where \( \cfi:=\cf(\zeta_{p^{\infty}}) \), \( H_{\cf}:=\Gal_{\cfi} \), \( \Gamma_{\cf}:=\Gal(\cfi/\cf)\subset \ZZ_{p}^{\times} \), \( \Darithcf(M):=H^{0}(H_{\cf},M)^{\Gamma_{\cf}-\la} \), \( \theta \) denotes a fixed generator of \( \Lie(\Gamma_{\cf}) \),
and \( E_{0}^{\theta} \) is the functor taking the subspace on which \( \theta \) acts nilpotently (taken with respect to the \( t \)-adic filtration). 
Then \( D_{\mrm{pdR},\cf}(\rho) \) has a finite filtration \( \Fil_{t}^{*} \) (inherited from the \( t \)-adic filtration on the RHS), the action of \( \theta \) preserves the filtration, and is zero on the graded pieces. We refer to this nilpotent operator as the \emph{Fontaine operator} \( \theta \) of \( \rho \). Then \( \rho \) is de Rham if and only if \( \theta=0 \)

In the locally analytic approach initiated by \cite{Pan22}, \cite{Pan2209.06II}, one identifies the Fontaine operators of locally analytic completed cohomology with differential operators on Shimura varieties of infinite level, as considered in \cite{Scholze15}. One can then calculate the kernels of these  differential operators, and relate them to automorphic forms.

This strategy works well in the case of Shimura curves as considered in \cite{Pan22} and \cite{QiuSu2025locallyanalyticvectorscompleted}. 
However, when trying to generalize to the case of Hilbert modular varieties, one runs into problems related to the \emph{plectic conjecture} (\cite{NekovarScholl2017Plectic}).
Let \( \Gal_{F}^{plec}:=\mrm{Aut}_{F}(F\otimes_{\QQ}\bar{\QQ}) \), which contains \( \Gal_{\QQ} \) as a subgroup. Then \cite{NekovarScholl2017Plectic} predicts that the action of \( \Gal_{\QQ} \)
on \( R\Gamma(X_{K,\bar{\QQ}},\ZZ/p^{m}) \)
extends naturally to an action of \( \Gal_{F}^{plec} \). In particular, conjecturally there is also a natural action of \( \Gal_{F}^{plec} \)
on \( \ti{H}^{i}(K^{p},\QQ_{p}) \). 

Note that non-canonically, we have \( \Gal_{F}^{plec}\cong \Gal_{F}^{\Sigma}\rtimes \mrm{Aut}(\Sigma)\supset \Gal_{\cf}^{\Sigma} \), with \( \Sigma:=\{\tau:F\hookrightarrow \cf \}=\{\tau_{1},\ldots,\tau_{d}\} \).
If \( \rho \) is a representation of \( \Gal_{F}^{plec} \), then we can consider \(
D_{\mrm{pdR},\cf}(\rho|_{\Gal_{\QQ_{p}}})
\). However, in this case, we have \[
D_{\mrm{pdR},\cf}(\rho|_{\Gal_{\QQ_{p}}})\cong \left(D_{\pdR,\cf}^{\tau_{1}}\circ D_{\pdR,\cf}^{\tau_{2}}\cdots \circ D_{\pdR,\cf}^{\tau_{d}}(\rho)\right)\otimes_{\cfi^{\otimes \Sigma}}\cfi,
\] where \( D_{\pdR,\cf}^{\tau} \) denotes the functor \( D_{\mrm{pdR},\cf} \) with respect to the \( \tau \)-th factor \( \Gal_{\cf} \), and \( \cfi^{\otimes\Sigma}\to \cfi \) is given by multiplication. The operations \( D_{\pdR,\cf}^{\tau} \)
commute with each other, and each one of them produces a \( t_{\tau} \)-adic filtration, and a \( \tau \)-partial Fontaine operator \( \theta_{\tau} \), such that the original Fontaine operator  \( \theta \) on the LHS can be written as \( \theta=\sum_{\tau\in\Sigma}\theta_{\tau} \), and the \( t \)-adic filtration on the LHS can be obtained as the total filtration of the multi-filtration on the RHS. 

If we believe in the plectic conjecture, we would wish to see such ``plectic structure'' (i.e. multi-filtration and partial Fontaine operators) from geometry. Let us introduce the geometric set-up.  
Let
\( X^{\tor}_{K^{p}K_{p}} \) be the toroidal compactification of \( X_{K^{p}K_{p}} \) for a fixed cone decomposition,
\( \mX^{\tor}_{K^{p}K_{p},\cf} \) be the analytification of \( X^{\tor}_{K^{p}K_{p}}\otimes_{\QQ} \cf \), and \( \mX^{\tor}_{K^{p}K_{p}}:=\mX^{\tor}_{K^{p}K_{p},\cf} \hatotimes_{\cf} \CC_{p} \).
We consider the Hilbert modular varieties of infinite level \[
\mX^{\tor}_{K^{p},\proket}:=\varprojlim_{K_{p}}\mX^{\tor}_{\Kpp},
\] where the inverse limit is taken in the pro-Kummer \'etale site of \( \mX_{K^{p}K_{p,0}}^{\tor} \) (for a fixed \( K_{p,0} \)) over all open compact subgroup \( K_{p}\subset G(\QQ_{p}) \). 
One can consider the pro-Kummer \'etale period sheaf \( \BdR^{+} \) defined in \cite{Scholze13} and \cite{DLLZ2022logarithmicJAMS}.
Denote by \( \BdRKp^{+} \) the push-forward of this period sheaf along the natural map \[ 
\nu:(\mX^{\tor}_{K^{p}K_{p,0}})_{\proket/\mX^{\tor}_{K^{p},\proket}}
\to 
|\mX^{\tor}_{K^{p},\cf}|_{\an}:=\varprojlim_{K_{p}}|\mX^{\tor}_{K^{p}K_{p},\cf}|_{\an}, \] which is filtered (up to a Tate twist) by the complete structure sheaf \( \hat{\mO}_{K^{p}} \). 
We denote by \( \BdRKp^{+,\la} \)
and \( \mO^{\la}_{K^{p}} \) their subsheaves consisting of \( G(\QQ_{p}) \)-locally analytic vectors. Then \cite{Pan22} and \cite{Juan2022.09locallyShi} show that the functor \( \Darithcf \) is well-behaved for these sheaves. Put \[
\mO^{\la}_{K^{p},\cfi}:=\Darithcf(\mO^{\la}_{K^{p}}),\;
\BdifKpcf^{+,\la}:=\Darithcf(\BdRKp^{+,\la}),\;\BdifKpcf^{\la}:=\BdifKpcf^{+,\la}[1/t].
\] Then  \( \BdifKpcf^{\la} \) is filtered by \( \mO^{\la}_{K^{p},\cfi}(i) \) for \( i\in\ZZ \). 

For a fixed infinitesimal character \( \lambda^{\kw}:Z(\fg)\to \QQ_{p} \) associated to a multiweight \( \kw \) (\cref{dfn:weight-HMF}),
we have by primitive comparison (\cite{Scholze12}, \cite[Theorem 6.2.6]{Juan2022.09locallyShi}) the following Hodge-Tate decomposition
\begin{align}\label{align:plec-HT-decomp-intro}
\ti{H}^{i}(K^{p},\QQ_{p})^{\la,\lambda^{\kw}}\hatotimes_{\QQ_{p}} \CC_{p}\cong \bigoplus_{T\subset \Sigma}H^{i}(|\mX^{\tor}_{K^{p},\cf}|_{\an},\mO^{\la,(s_{T}\cdot \Bbbk,w)}_{K^{p},\cfi})\left(\sum_{\tau\in\Sigma}\frac{w-(s_{T}\cdot \Bbbk)_{\tau}}{2}\right),
\end{align}
 where \( \mO^{\la,(s_{T}\cdot \Bbbk,w)}_{K^{p},\cfi} \) is (up to a Tate twist) a certain subsheaf of \( \mO^{\la}_{K^{p},\cfi} \) defined using horizontal actions (\cite{Pan22}, \cite{Juan2022.09locallyShi}), and  \( (s_{T}\cdot \Bbbk)_{\tau}:=k_{\tau} \) if \( \tau\notin T \), and \( :=2-k_{\tau} \) if \( \tau\in T \). This decomposition can be thought of as a \emph{plectic Hodge-Tate decomposition}, where the Sen operator \( \theta \) decomposes into a sum of partial Sen operators \( \theta_{\tau} \) which acts on the \( T \)-th component as \( \frac{w-(s_{T}\cdot \Bbbk)_{\tau}}{2} \). Note that this idea is used in \cite{JiangPan2025Plectic} to construct a certain plectic Lie algebra action on \( \ti{H}^{i}(K^{p},\QQ_{p})^{\la} \). 

We can continue to apply the functor \( E_{0}^{\theta} \). It follows easily from (\ref{align:plec-HT-decomp-intro})
that \begin{align}\label{align:DpdR-class-intro}
D_{\mrm{pdR},\cf}(\ti{H}^{d}(K^{p},\QQ_{p})^{\la,\lambda^{\kw}}) 
\cong H^{d}(|\mX^{\tor}_{K^{p},\cf}|_{\an},E_{0}^{\theta}(\BdifKpcf^{\la,\lambda^{\kw}})),
\end{align}
and \( E_{0}^{\theta}(\BdifKpcf^{\la,\lambda^{\kw}})) \) is equipped with an induced filtration \( \Fil^{*}_{t} \) with non-trivial graded pieces being \( \mO^{\la,(s_{T}\cdot \Bbbk,w)}_{K^{p},\cfi} \)
for \( T\subset \Sigma \). 
The following theorem constructs the desired ``plectic structure'' on \( E_{0}^{\theta}(\BdifKpcf^{\la,\lambda^{\kw}})) \):

\begin{thm}[\cref{thm:HMV-Fontaine=diff}]\label{intro:thm:HMV-Fontaine=diff}
Let \( \kw \) be a multiweight with \( k_{\tau}\ge 2 \) for all \( \tau\in\Sigma \). 

(1) \( E_{0}^{\theta}(\BdifKpcf^{\la,\lambda^{\kw}}) \) admits a descending \( \ZZ^{\Sigma} \)-filtration \( \Fil^{*}_{(t_{\tau})_{\Sigma}} \) whose total filtration is \( \Fil^{*}_{t} \), and for which the only \emph{non-trivial} graded pieces are 
\begin{equation}\label{eq:HT-grade-HM-intro}
\gr^{\unl{j}_{T}}_{(t_{\tau})_{\Sigma}}E_{0}^{\theta}(\BdifKpcf^{\la,\lambda^{\kw}})\cong \mO^{\la,(s_{T}\cdot \Bbbk,w)}_{K^{p},\cfi},
\end{equation}
 for \( T\subset \Sigma \)
and \( \unl{j}_{T}:=\left(-\frac{w-(s_{T}\cdot\Bbbk)_{\tau}}{2}\right)_{\tau} \). 

(2) There exist commuting operators \( \theta_{\tau}\in \End(E_{0}^{\theta}(\BdifKpcf^{\la,\lambda^{\kw}})) \) for \( \tau\in\Sigma \)
such that \(
\theta=\sum_{\tau\in\Sigma}\theta_{\tau}.
\) 

(3)
For \( j\in\ZZ \), denote by \( \Fil_{t_{\tau_{0}}}^{j}E_{0}^{\theta}(\BdifKpcf^{\la,\lambda^{\kw}}) \) the filtration obtained by forgetting all the filtrations except the \( \tau_{0} \)-th component,
and denote by \( \gr_{t_{\tau_{0}}}^{j} \) the associated graded piceces.

For \( \tau_{0}\in\Sigma \), the action of \( \theta_{\tau_{0}} \) preserves the multi-filtration \( \Fil^{*}_{(t_{\tau})_{\Sigma}} \), is zero on \( \gr^{*}_{t_{\tau_{0}}}(E_{0}^{\theta}(\BdifKpcf^{\la,\lambda^{\kw}}) ) \), and its action on \( E_{0}^{\theta}(\BdifKpcf^{\la,\lambda^{\kw}}) \)
induces monodromy operators 
\[
\ti{N}_{\tau_{0},k_{\tau_{0}}}:
\gr^{\unl{j}_{T\cup\{\tau_{0}\}}}_{(t_{\tau})_{\Sigma}}(E_{0}^{\theta}(\BdifKpcf^{\la,\lambda^{\kw}}))
\to 
\gr^{\unl{j}_{T}}_{(t_{\tau})_{\Sigma}}(E_{0}^{\theta}(\BdifKpcf^{\la,\lambda^{\kw}})),
\] for any \( \tau_{0}\notin T\subset \Sigma \).
Moreover, up to the isomorphism (\ref{eq:HT-grade-HM}), we have
\begin{equation}\label{eq:HMF-Nk-theta-intro}
\ti{N}_{\tau_{0},k_{\tau_{0}}}=c_{k}\cdot \theta^{k_{\tau}-1}_{\Fl^{\std}_{G,\mu}}\boxtimes \theta^{k_{\tau}-1}_{\Fl_{G,\mu}},
\end{equation}
where \( c_{k_{\tau}}\in\QQ^{\times} \) is as in \cref{thm:Vin-GL2-Fontaine=Diff}, and \( \theta^{k_{\tau}-1}_{\Fl^{\std}_{G,\mu}}\boxtimes \theta^{k_{\tau}-1}_{\Fl_{G,\mu}} \) is as in \cref{notation:Ola-kw}.
\end{thm}
\begin{rmk}
Unwinding the definition, the formula (\ref{eq:HMF-Nk-theta-intro}) is exactly the identity \( N=d\circ \bar{d} \) of \cite{Pan2209.06II}. 
\end{rmk}

\cref{intro:thm:HMV-Fontaine=diff} plus the results in \cite{Jiang2025HMF} implies easily \cref{Intro:thm:ClassicalityCompleteCoho}. Although the statement of \cref{intro:thm:HMV-Fontaine=diff} only involves concrete sheaves on a topological space, our proof involves the deep theory of analytic stacks developed by Clausen-Scholze. Roughly, we define analytic stacks \( \BdR^{+}(\mX^{\tor}_{K^{p}})^{\la} \)
and \( \Bdifcf^{+}(\mX^{\tor}_{K^{p}})^{\la} \) as the (analytic) relative spectrum of \( \BdRKp^{+,\la} \)
and \( \BdifKpcf^{+,\la} \) over the topological space \( |\mX^{\tor}_{K^{p},\cf}|_{\an} \). See \cref{construction:BdR-Xla} for details. In particular, \( \BdR^{+}(\mX^{\tor}_{K^{p}})^{\la} \)
is a \( B_{\dR}^{+} \)-thickening of the locally analytic Shimura variety \( \mX^{\tor,\la}_{K^{p}} \)
defined in \cite{Jiang2026Shla}. 

We want to formulate a Grothendieck-Messing theory for \( \BdRKp^{+,\la} \)
and \( \BdifKpcf^{+,\la} \), in the spirit of \cite[\cref{1-thmLAstrucGMHTperIntro}]{Jiang2026Shla}. For this, we define a certain ``algebraic de Rham stack'' \( \bar{\mX}^{\tor,\la}_{K^{p},\cfi} \) of the infinite-level Hilbert modular varieties, such that \[
\bar{\mX}^{\tor,\la}_{K^{p},\cfi}\hatotimes_{\cfi}\CC_{p}\cong (\mX^{\tor,\la,0}_{K^{p}})^{\pdR/\CC_{p},++}_{\log}\times_{[\mA^{1}/\GG_{m}]^{2}}[\GG_{m}/\GG_{m}]^{2},
\] where \( (\mX^{\tor,\la,0}_{K^{p}})^{\pdR/\CC_{p},++}_{\log} \) is defined in \cite[\cref{1-dfnBiFilAgLa0}]{Jiang2026Shla}. 
We have natural maps \[
\mX^{\tor,\la}_{K^{p}}\to \BdR^{+}(\mX^{\tor}_{K^{p}})^{\la}
\to \Bdifcf^{+}(\mX^{\tor}_{K^{p}})^{\la}
\xrightarrow{h^{\Bdif}} \bar{\mX}^{\tor,\la}_{K^{p},\cfi}\to |\mX^{\tor}_{K^{p},\cf}|_{\an},
\] and \( h^{\Bdif}_{*}(\mO_{\Bdifcf^{+}(\mX^{\tor}_{K^{p}})^{\la}}) \) offers a natural lift of \( \BdifKpcf^{+,\la} \)
from \( |\mX^{\tor}_{K^{p},\cf}| \)
to \( \bar{\mX}^{\tor,\la}_{K^{p},\cfi} \). 
\begin{rmk}
\( \Bdifcf^{+}(\mX^{\tor}_{K^{p}})^{\la}/K_{p}^{\la} \) can be reinterpreted as (a logarithmic variant of) the infinitesimal neighborhood of the Hodge-Tate locus of the \( p \)-adic twistor stack of \( \mX^{\tor}_{\Kpp} \) of \cite{ALBRCS2026rational}.
\end{rmk}

The proof of \cref{intro:thm:HMV-Fontaine=diff} requires a good understanding of \( h^{\Bdif} \). For this end, we prove the following theorem for \emph{any} Shimura datum, which provides a local model of \( \Bdifcf^{+}(\mX^{\tor}_{K^{p}})^{\la} \):
\begin{thm}[\cref{cor:BdR-Base-Change}]
	\label{intro:thm:GM-theory-BdR}
There is a Cartesian diagram over \( \cfi \) \begin{equation}\label{eq:Alg-DR-Bdif2}
\begin{tikzcd}
\Bdifcf^{+}(\mX^{\tor}_{K^{p}})^{\la}\arrow[d,"h^{\Bdif}"]\arrow[r,"\pi_{\GM\HT}^{\Bdif}"] & \Vin^{*}_{G^{c},\mu,\cfi}|_{\hat{\mA}^{1}}/G^{c}\arrow[d,"\hat{h}_{0}"]\arrow[r,"\hat{i}"]
& \Vin^{*}_{G^{c},\mu,\cfi}/G^{c}\arrow[d,"h"]\\
\bar{\mX}^{\tor,\la}_{K^{p},\cfi}\arrow[r,"\bar{\pi}_{\GM\HT}"] & (\Vin^{*}_{G^{c},\mu,\cfi}|_{0})^{\pdR/\cfi}/G^{c}
\arrow[r,"i^{\dR}"] & (\Vin^{*}_{G^{c},\mu,\cfi})^{\pdR/\cfi}/G^{c},
\end{tikzcd}
\end{equation}
where \( \Vin^{*}_{G^{c},\cfi} \) denotes the non-degenerate locus of the Vinberg semigroup of \( G^{c} \) (base changed to \( \cfi \)), \( \Vin^{*}_{G^{c},\mu,\cfi}:=\Vin^{*}_{G^{c},\cfi}\times_{\bar{T}^{c}_{ad}}\mA^{1} \to \mA^{1}\) where \( \mA^{1}\to \bar{T}^{c}_{ad} \) is determined by \( \mu \), and \( -|_{\hat{\mA}^{1}} \)
and \( -|_{0} \) denotes its fiber over \( \hat{\mA}^{1} \)
and \( \{0\} \) respectively. 

Moreover, 
the action of \( \theta\in\Lie(\Gamma_{\cf}) \)
on \( \Bdifcf^{+}(\mX^{\tor}_{K^{p}})^{\la} \) coincides with the action of \( \theta_{\mu} \) (as in \cref{dfn:Vin-G-mu}) on the RHS.
\end{thm}
See \S \ref{sec:Vinberg} for a review of the theory of Vinberg semigroups. 
The construction of the commutative diagram uses the logarithmic Riemann-Hilbert correspondence (\cite{LiuZhu2017rigidity}, \cite{DLLZ2022logarithmicJAMS}). To prove that the diagram is Cartesian, since \( \{0\}\to \hat{\mA}^{1} \) is an epimorphism, we can localize to \( \Vin^{*}_{G^{c},\mu,\cfi}|_{0} \), which is precisely the Grothendieck-Messing-Hodge-Tate period domain \( \Per_{G,\mu} \)
considered in \cite{Jiang2026Shla}. Thus we can conclude using the Grothendieck-Messing theory for locally analytic infinite-level Shimura varieties established in \cite[\cref{1-thmLAstrucGMHTperIntro}]{Jiang2026Shla}. 
\begin{eg}
In the case of modular curves,
we have \[ \Vin_{\GL_{2}}^{*}\cong \Vin_{\GL_{2},\mu}^{*}\cong M_{2\times 2}^{*}\times \GG_{m}\xrightarrow{(X,a)\mapsto \det(X)a^{-1}} \mA^{1}, \] where \( M_{2\times 2}^{*} \) denotes the space of non-zero \( 2\times 2 \)-matrices, and along this isomorphism, the operator \( \theta_{\mu} \) is the derivative of the diagonal scalar action of \( \GG_{m} \) on \( M_{2\times 2}^{*} \)
and \( \GG_{m} \). 

Then
by Riemann-Hilbert correspondence, we have a natural isomorphism \[
T_{p}E_{K}\otimes_{\ZZ_{p}}\BdRKp\cong R\Gamma(\nabla, \mcal{H}_{1}^{\dR}(E_{K}/\mX_{K}^{\tor})\otimes_{\mO_{\mX_{K}^{\tor}}} \OBdR),
\] where \( E_{K}\to \mX_{K}^{\tor} \) denotes the (generalized) universal elliptic curve, \( \mcal{H}_{1}^{\dR} \) denotes the relative (log) de Rham homology, and
\( \OBdRKp \) is the period sheaf in \cite{Scholze12} and \cite{DLLZ2022logarithmicJAMS}. Multiplying this morphism by \( t \) induces a morphism \[
T_{p}E_{K}\otimes_{\ZZ_{p}}\BdRKp^{+}\hookrightarrow  R\Gamma(\nabla, \mcal{H}_{1}^{\dR}(E_{K}/\mX_{K}^{\tor})\otimes_{\mO_{\mX_{K}^{\tor}}} \OBdR^{+}).
\]
Let \( \mcal{G}_{\dR} \) denote the de Rham \( G \)-torsor on \( \bar{\mX}^{\tor,\la}_{K^{p}} \). 
Passing to infinite level, we have a trivialization of \( T_{p}E_{K} \), and pulling back to \( \mcal{G}_{\dR} \), we have a trivialization of the RHS. 
We thus obtain a map \[
\mrm{Std}\otimes h^{\Bdif,*}\mcal{G}_{\dR}\to \mrm{Std}\otimes h^{\Bdif,*}\mcal{G}_{\dR},
\] which induces a \( G \)-equivariant morphism \( h^{\Bdif,*}\mcal{G}_{\dR}\to M_{2\times 2}^{*} \). The map to the other component \( \GG_{m} \) is obtained by considering \( \det \) instead of \( \mrm{Std} \). 
\end{eg}
Finally, 
\cref{intro:thm:HMV-Fontaine=diff} follows from \cref{intro:thm:GM-theory-BdR} by smooth base change and a concrete calculation of \( G^{c} \)-equivariant \( D \)-modules on \( \Vin^{*}_{G^{c},\mu,\cfi} \). Concretely, we have a Cartesian diagram 
\[
\begin{tikzcd}
\Vin^{*}_{G^{c},\mu}\arrow[r]\arrow[d,"\mfk{t}"] & \Vin^{*}_{G^{c}}\arrow[d,"\mfk{t}"]\\
\mA^{1}\arrow[r,"\Delta"]
& (\mA^{1})^{\Sigma}.
\end{tikzcd}
\]
Thus we have a natural enlargement of \( \Vin_{G^{c},\mu}^{*}|_{\hat{\mA}^{1}} \)
\[
\begin{tikzcd}
\Vin_{G^{c},\mu}^{*}|_{\hat{\mA}^{1}}\arrow[r,"\Delta",hook]\arrow[d,"\hat{h}_{0}"] & 
\Vin_{G^{c}}^{*}|_{(\hat{\mA}^{1})^{\Sigma}}\arrow[dl,"\hat{h}"]\\
(\Vin_{G^{c},\mu}^{*}|_{0})^{\pdR/\cfi},
\end{tikzcd}
\]  \( \hat{h}_{*}(\mO_{\Vin^{*}_{G^{c}}|_{(\hat{\mA}^{1})^{\Sigma}}}) \) is endowed with a \( (t_{\tau}^{\bullet}) \)-adic \( \ZZ^{\Sigma} \)-filtration, where \( t_{\tau} \) is the coordinate of the \( \tau \)-th component \( \mA^{1} \), and we have a natural morphism \[p:\hat{h}_{*}(\mO_{\Vin^{*}_{G^{c}}|_{(\hat{\mA}^{1})^{\Sigma}}})\twoheadrightarrow 
\hat{h}_{0,*}(\mO_{\Vin^{*}_{G^{c},\mu}|_{\hat{\mA}^{1}}})
\cong \hat{h}_{*}(\mO_{\Vin^{*}_{G^{c}}|_{(\hat{\mA}^{1})^{\Sigma}}})\otimes_{\QQ[t_{\tau}:\tau\in\Sigma]}\QQ[t],
\] where \( \QQ[t_{\tau}:\tau\in\Sigma]\to \QQ[t],t_{\tau}\mapsto t \). 
Now the action of \( \theta_{\mu} \) on \( \hat{h}_{*}(\mO_{\Vin^{*}_{G^{c}}|_{(\hat{\mA}^{1})^{\Sigma}}}) \) naturally decomposes into a sum of \( \theta_{\tau} \) for \( \tau\in\Sigma \), and  \( p \) induces an isomorphism \[
E_{0}^{\theta_{\tau_{1}}}\circ \cdots E_{0}^{\theta_{\tau_{d}}}(\hat{h}_{*}(\mO_{\Vin^{*}_{G^{c}}|_{(\hat{\mA}^{1})^{\Sigma}}})^{\lambda^{\kw}})\cong E_{0}^{\theta_{\mu}}(\hat{h}_{0,*}(\mO_{\Vin^{*}_{G^{c},\mu}|_{\hat{\mA}^{1}}})^{\lambda^{\kw}}).
\] This allows us to translate the \( \ZZ^{\Sigma} \)-filtration
and the action of \( \theta_{\tau} \) on the LHS to those on the RHS. Note that up to center, \( \Vin_{G^{c}}^{*} \) is essentially a product over \( \Sigma \) of \( \Vin_{\GL_{2}}^{*} \). Thus the proof of the  formula (\ref{eq:HMF-Nk-theta-intro}) also reduces to the case of \( \GL_{2} \), where it essentially follows from the classical Capelli identity (\ref{align:capelli}) (\cite{Capelli1887zuruckfuhrung}). 

In reality, we will use \( \ti{h}_{0,\natural} \) instead of \( \ti{h}_{0,*} \) so that we could apply smooth base change to \cref{intro:thm:GM-theory-BdR}.
See \cref{lem:pullback-compare-diagonal} for a general statement. 

\subsection{Organization of the paper}
The paper is organized as follows:
\S \ref{subsec:technical-lem} collects some technical results about equivariant \( D \)-modules. \S \ref{sec:Vinberg} gives a review of Vinberg semigroups. 
In
\S \ref{sec:BdR-thick}, we define the stacks \( \bar{\mX}^{\tor,\la}_{K^{p},\cfi} \) and \( \Bdifcf^{+}(\mX^{\tor}_{K^{p}})^{\la} \) appearing above, and establish \cref{intro:thm:GM-theory-BdR}. 
In \S \ref{sec:Font-op-on-Vin}, we do the concrete calculation on Vinberg semigroups.
In \S \ref{sec:Fontaine-Ope-HMV}, we deduce \cref{intro:thm:HMV-Fontaine=diff} from \cref{intro:thm:GM-theory-BdR} and the results in \S \ref{sec:Font-op-on-Vin}. 
In \S \ref{sec:Classicality-Non-Parallel}, we apply the previous theorem to deduce \cref{Intro:thm:ClassicalityCompleteCoho} and \cref{Intro:thmConjecFontaineMazur}. 
\subsection{Acknowledgement}
I am extremely grateful to Arthur-C\'esar Le Bras for suggesting the idea of studying Fontaine operators using Vinberg semigroups, an idea motivated by his joint work in progress with Ansch\"utz, Rodr{\'\i}guez Camargo, and Scholze.
I also thank him for helpful comments on an earlier draft of this paper. 
I would like to thank Johannes Anschütz, George Boxer, Lue Pan, Vincent Pilloni, Juan Esteban Rodr\'iguez Camargo, Peter Scholze, and Qixiang Wang for helpful discussions. 
I am grateful to Max Planck Institute for Mathematics in Bonn for its hospitality and financial support.

\subsection{AI disclosure}
Gemini 3.1 Pro (March 2026) was used for reference searching. 
All the writings and mathematics were done by the author.

\subsection{Convention}
We will work with solid modules in the sense of \cite{CS19}. To ease the notation, all the tensor products will be taken in the solid sense, which we will simply denote as \(-\otimes-\).

We fix \(\bar{\QQ}\subset\CC\),
and we fix an isomorphism \(\iota:\CC\cong\bar{\QQ}_{p}\). 
For a group \( G \) and a space \( X \), a \( G \)-action on \( X \) always refers to a right action. Denote by \(\CC_{p}\) the completion of \(\bar{\QQ}_{p}\) for the \(p\)-adic topology.
By Frobenius, we will always refer to the geometric Frobenius. 
Let \(\chi_{\mrm{cycl}}:\Gal_{\QQ}\to\ZZ_{p}^{\times}\)
denote the cyclotomic character, sending \(\Fr_{l}\) to \(l^{-1}\). By our convention, \(\chi_{\mrm{cycl}}\)
has Hodge-Tate weight \(-1\). Our convention for local class field theory takes the uniformizer to a lift of the geometric Frobenius. 

For a group \( G \) and a space \( X \), a \( G \)-action on \( X \) always refers to a \emph{right} action. 

\section{Some technical lemmas}
\label{subsec:technical-lem}
We start with some technical lemmas which will be used later. Readers are welcomed to jump to the next section, and refer back when necessary. 
We will work with schemes over \( \QQ \), which are realized as solid stacks 
via \cite[\cref{1-egScheme}]{Jiang2026Shla}.

\subsection{Na\"ive dual}
We will use the following notion of na\"ive dual \( \mbb{D}_{X} \).
\begin{notation}\label{notation:Naive-dual}
 For any solid stack \( U \), we denote \[
\mbb{D}_{U}:\QCoh(U)^{op}\to \QCoh(U),M\mapsto \mbb{D}_{U}(M):=M^{\vee}:=\unl{\RHom}_{U}(M,\mO_{U}).
\]
 \end{notation}
\begin{notation}[{\cite[\cref{1-dfnFnaturalFunctor}]{Jiang2026Shla}}]\label{dfnFnaturalFunctor}
For any cohomologically smooth map \(h:X\to Y\), we write \(h_{\natural}(-)\cong h_{!}(-\otimes h^{!}1)\)
for the left adjoint of \(h^{*}\). 

Note that by \cite[\cref{1-lemSmBaseChange}]{Jiang2026Shla}, we have \emph{smooth base change}: for any Cartesian diagram \[
\begin{tikzcd}
X'\arrow[r,"g_{X}"]\arrow[d,"h'"] & X\arrow[d,"h"]\\
Y'\arrow[r,"g_{Y}"] & Y,
\end{tikzcd}
\] where \( h \) is cohomologically smooth, we have natural isomorphisms \[
h'_{\natural}g_{X}^{*}\cong g_{Y}^{*}h_{\natural},\;h^{*}g_{Y,*}\cong g_{X,*}h^{\prime,*}.
\]
\end{notation}
\begin{eg}[{\cite[\cref{1-lemBasicFiltered} \& \cref{1-lemDRcomplexVSpushforward}]{Jiang2026Shla}}]\label{eg:D-mod-lower-natural}
Let \( X \) be a log smooth rigid variety over \( \cf \). Then \( h:X\to X^{\pdR/\cf}_{\log} \) is cohomogically smooth, and \( h_{\natural}(\mO_{X})\cong D_{X/\cf,\log} \) if we identify \( \QCoh(X^{\pdR/L}) \) with the category of \( D_{X/\cf,\log} \)-modules via \cite[\cref{1-lemDRcomplexVSpushforward}]{Jiang2026Shla}.
Moreover, by adjunction, for any vector bundle \( V \) on \( X \), \[h_{*}V\cong (h_{\natural}V^{\vee})^{\vee}\cong  \unl{\RHom}_{X}(D_{X/L,\log}\otimes_{\mO_{X}} V^{\vee},\mO_{X}), \]
where the RHS
is endowed with a (left) \( D_{X/\cf,\log} \)-module structure by letting \( D\in T_{X/L,\log} \)
act via \[
D\cdot f:=D\circ f-f\circ D.
\]
\end{eg}
\begin{lem}\label{lem:pullback-and-dual}
(1) Consider any morphism \( g:X\to Y \) between solid stacks. Then we have a natural transformation \[
g^{*}\circ \mbb{D}_{Y}\to \mbb{D}_{X}\circ g^{*},
\] 
which is an isomorphism if \( g \) is cohomologically smooth.
This natural transformation is compatible with compositions: for \( X\xrightarrow{g}Y\xrightarrow{f}Z \), we have a commutative diagram \[
\begin{tikzcd}
(f\circ g)^{*}\circ \mbb{D}_{Z}\arrow[r]\arrow[d] & 
f^{*}\circ g^{*}\circ \mbb{D}_{Z}\arrow[r]
& f^{*}\circ \mbb{D}_{Y}\circ g^{*}\arrow[d]\\\mbb{D}_{X}\circ (f\circ g)^{*}\arrow[rr]
&& \mbb{D}_{X}\circ f^{*}\circ g^{*}.
\end{tikzcd}
\]

(2) Assume that we are give a Cartesian diagram of solid stacks \[
\begin{tikzcd}
X'\arrow[r,"g_{X}"] \arrow[d,"f'"] & X\arrow[d,"f"] \\
Y'\arrow[r,"g_{Y}"] & Y,
\end{tikzcd}
\] 
where \( f \) is cohomologically smooth.
Then we have the following commutative diagram
\begin{equation}\label{eq:natu-transfo-commute}
\begin{tikzcd}
g_{Y}^{*}\circ \mbb{D}_{Y}\circ f_{\natural}
\arrow[r,"\eta_{1}"]\arrow[d,"\eta_{3}","\cong"']
& \mbb{D}_{Y'}\circ g^{*}_{Y}\circ f_{\natural}
\arrow[r,"\eta_{2}","\cong"']
& \mbb{D}_{Y'}\circ f'_{\natural}\circ g_{X}^{*}\arrow[d,"\eta_{6}","\cong"']
\\
g^{*}_{Y}\circ f_{*}\circ \mbb{D}_{X}
\arrow[r,"\eta_{4}"]
& 
f'_{*}\circ g_{X}^{*}\circ \mbb{D}_{X}
\arrow[r,"\eta_{5}"]
& f'_{*}\circ \mbb{D}_{X'}\circ g^{*}_{X},
\end{tikzcd}
\end{equation}
 where \( \eta_{1} \)
and \( \eta_{5} \)
are as in (1).
\end{lem}
\begin{proof}
(1) For \( M\in \QCoh(Y) \),
we  have a natural map \( M\otimes \mbb{D}_{Y}(M)\to \mO_{Y} \), and thus \( g^{*}M\otimes g^{*}\mbb{D}_{Y}(M)\to \mO_{Y'} \). By adjunction, this corresponds to a natural map \( g^{*}\mbb{D}_{Y}(M)\to \mbb{D}_{Y'}g^{*}(M) \). This construction is compatible with compositions as \( g^{*} \) is symmetric monoidal.
When \( g \) is cohomologically smooth, we have 
\begin{align*}
g^{*}\mbb{D}_{Y}(M)
\cong g^{!}\mbb{D}_{Y}(M)\otimes g^{!}(\mO_{Y})^{-1}
\cong \unl{\RHom}_{X}(g^{*}(M),g^{!}(\mO_{Y}))\otimes g^{!}(\mO_{Y})^{-1}
\cong \mbb{D}_{X}(g^{*}(M)).
\end{align*}

(2) Consider \( f\circ g_{X}\cong g_{Y}\circ f' \). Thus by (1), we have a commutative diagram
\begin{equation}\label{eq:nat-trans-commute1}
\begin{tikzcd}
f^{\prime,*}g_{Y}^{*}\mbb{D}_{Y}\arrow[r]\arrow[d] & f^{\prime,*}\mbb{D}_{Y'}g_{Y}^{*} \arrow[r] & 
\mbb{D}_{X'}f^{\prime,*}g_{Y}^{*}\arrow[d]
\\
g_{X}^{*}f^{*}\mbb{D}_{Y} \arrow[r]
& g_{X}^{*}\mbb{D}_{X}f^{*} \arrow[r] 
& \mbb{D}_{X'}g_{X}^{*}f^{*}.
\end{tikzcd}
\end{equation}
This gives an isomorphism between two natural transformations \[
f^{\prime,*}g_{Y}^{*}\mbb{D}_{Y}\to
\mbb{D}_{X'}g_{X}^{*}f^{*},
\]
which by adjunction gives an isomorphism between two natural transformations \[
g_{Y}^{*}\mbb{D}_{Y}f_{\natural}\to f^{\prime}_{*}\mbb{D}_{X'}g_{X}^{*}.
\]
It is straight-forward to show that the two natural transformations are isomorphic those appearing in (\ref{eq:natu-transfo-commute}). For example, the composition of the right vertical arrow and the upper horizontal arrows in (\ref{eq:nat-trans-commute1}) fits into the following commutative diagram 
\begin{equation}
	\begin{tikzcd}
	g_{Y}^{*}\mbb{D}_{Y}f_{\natural}
\arrow[r] & \arrow[d] \mbb{D}_{Y'}g^{*}_{Y}f_{\natural}
\arrow[r] & \arrow[d] 
f'_{*}\mbb{D}_{X'}f^{\prime,*}g_{Y}^{*}f_{\natural} \\ 
& \mbb{D}_{Y'}f'_{\natural}g_{X}^{*}\arrow[r] & 
 f'_{*}\mbb{D}_{X'}f^{\prime,*}f'_{\natural}g_{X}^{*}\arrow[r] & f'_{*}\mbb{D}_{X'}g_{X}^{*},
	\end{tikzcd}
\end{equation} which gives an isomorphism to the composition of the right vertical arrow and the upper horizontal arrows in (\ref{eq:natu-transfo-commute}).
\end{proof}

In \cite{Jiang2026Shla}, we construct an action of \( Z(\fg) \) on \( \QCoh(X^{\pdR/\cf}/G) \) as follows:
\begin{construction}[{\cite[\cref{1-lemInfiniteAction}]{Jiang2026Shla}}]\label{constr:infi-char-on-eqiv-D-mod}
Let \( Z(\fg):=Z(U(\fg)) \). Then by considering the following correspondence \[
\begin{tikzcd}
& (X^{\pdR/\cf}\times^{G} G^{\pdR/\cf})/G\arrow[r,"s"]\arrow[dl,"p_{1}"]\arrow[dr,"m"] & G\backslash G^{\pdR/\cf}/G
\\
X^{\pdR/\cf}/G && X^{\pdR/\cf}/G,
\end{tikzcd}
\] where \( m \) is induced by the action \( X^{\pdR/\cf}\times G^{\pdR/\cf}\to X^{\pdR/\cf} \),
we obtain a (right) action of \( \QCoh(G\backslash G^{\pdR/\cf}/G) \) (equipped with the convolution product)
on \( \QCoh(X^{\pdR/\cf}/G) \): \[
M*L :=m_{\natural}(p_{1}^{*}(M)\otimes s^{*}(L)),\;L\in\QCoh(G\backslash G^{\pdR/\cf}/G),M\in\QCoh(X^{\pdR/\cf}/G).
\]
The unit of \(  \QCoh(G\backslash G^{\pdR/\cf}/G) \) is given by \( \Delta_{\natural}(\mO_{G\backslash G/G}) \)
for \[
\Delta:G\backslash G/G\to G\backslash G^{\pdR/\cf}/G,
\] and we have an action of \( Z(\fg) \) on \( \QCoh(X^{\pdR/\cf}/G) \) induced by \( Z(\fg)\to \End(\Delta_{\natural}(\mO_{G\backslash G/G})) \).
\end{construction}
The following lemma gives an explicit example of this action of \( Z(\fg) \). 
\begin{lem}\label{lem:explicit-Zg-action-DX}
Consider the Cartesian diagram (\ref{eq:D-XdR/G}). Then \( h^{\prime,*}f^{*}h_{\natural}(\mO_{X/G})\cong D_{X/\cf} \), and the action of \( Z(\fg) \) on \( h_{\natural}(\mO_{X/G}) \) induces an action on \( D_{X/\cf} \), which can be described as follows: 

The action of \( G \) 
on \( X \) induces a morphism \( \fg\to T_{X} \), and \( U(\fg)^{op} \) acts on \( D_{X/\cf} \)
via right multiplication. Then \( Z(\fg) \)-action on \( D_{X/\cf} \) is induced by \( Z(\fg)\cong Z(U(\fg)^{op}) \), where the isomorphism is induced by \( A\mapsto A \) for \( A\in\fg \).
\end{lem}
\begin{proof}
The first isomorphism follows from 
Note that \( \hat{\fg}=\Ker(G\to G^{\pdR/\cf}) \) acts trivially on \( X^{\pdR/\cf} \). Thus (\ref{eq:D-XdR/G}) decomposes into two Cartesian squares \begin{equation}\label{eq:D-XdR/G-g}
\begin{tikzcd}
X\arrow[r,"h'"] \arrow[d,"f_{0}'"] & X^{\pdR/\cf}\arrow[d,"f_{0}"]\\
X/\hat{\fg} \arrow[r,"h_{1}"]\arrow[d]
& X^{\pdR/\cf}\times B\hat{\fg}\arrow[d,"u"]\\
X/G \arrow[r,"h"]
& X^{\pdR/\cf}/G.
\end{tikzcd} 
\end{equation} 
Moreover, by a similar construction as in \cref{constr:infi-char-on-eqiv-D-mod}, we have an action of \( \QCoh(\hat{\fg}\backslash */\hat{\fg}) \) on \( X^{\pdR/\cf}\times B\hat{\fg} \), such that by smooth base change, 
for any \( L\in \QCoh(G\backslash G^{\pdR/\cf}/G) \)
and \( M\in \QCoh(X^{\pdR/\cf}/G) \), we have \[
u^{*}(M*L)\cong u^{*}(M)*u_{0}^{*}(L),
\] where \( u_{0}:\hat{\fg}\backslash */\hat{\fg}\to G\backslash G^{\pdR/\cf}/G\). Thus it suffices to verify the similar statement for the upper square, that is, to show that the right \( U(\fg) \)-module structure on \( D_{X/\cf} \) (that corresponds to the descent of \( h'_{\natural}\mO_{X} \) to \( h_{1,\natural}\mO_{X/\hat{\fg}} \)) is induced by the morphism \( U(\fg)\to D_{X/\cf} \) and right multiplication.

By adjunction, \( \mbb{D}_{X^{\pdR/\cf}\times B\hat{\fg}}(h_{1,\natural}(\mO_{X/\hat{\fg}}))\cong h_{1,*}(\mO_{X/\hat{\fg}}) \).
As \( f_{0} \) is cohomologically smooth, \( f_{0}^{*} \)
commutes with the na\"ive dual by \cref{lem:pullback-and-dual} (1). 
By \cite[\cref{1-lemDiagonalLogProperty} (3)]{Jiang2026Shla}, \( \mbb{D}_{X^{\pdR/\cf}}(h'_{*}\mO_{X})\cong h'_{\natural}\mO_{X} \), 
and thus \( \mbb{D}_{X^{\pdR/\cf}\times B\hat{\fg}}(h_{1,*}\mO_{X/\hat{\fg}})\cong h_{1,\natural}\mO_{X/\hat{\fg}}\).
Thus it suffices to determine the left \( U(\fg) \)-module structure on \( h'_{*}(\mO_{X}) \) (that corresponds to its descent to \( h_{1,*}\mO_{X/\hat{\fg}} \)).
For this, we can localize to reduce to the case where \( X=\Spec(B) \), and then by definition, the \( U(\fg) \)-action on \( R\Gamma(X^{\pdR/\cf},h'_{*}\mO_{X})\cong B \)
is precisely induced by \( \fg\to T_{B/\cf} \). 
Moreover, we have a natural map \( h'_{*}\mO_{X}\to h'_{*}h^{\prime,*}h'_{*}\mO_{X} \), which after taking \( R\Gamma(X^{\pdR/\cf},-) \) can be described as \( B\to (B\otimes_{\cf} B)^{\wedge}_{I},b\mapsto 1\otimes b \),  where \( I:=\Ker(B\otimes_{\cf}B\to B) \).
Then the \( U(\fg) \)-action on \( h^{\prime,*}h'_{*}\mO_{X}\cong  (B\otimes_{\cf} B)^{\wedge}_{I} \) is the unique continuous action that is \( B\otimes 1 \)-linear, and 
acts on the second factor via  \( \fg\to T_{B/\cf} \).
Note that the pairing between \( D_{B/\cf} \)
and \( (B\otimes_{\cf} B)^{\wedge}_{I} \) is induced by \[
D_{B/\cf}\otimes (B\otimes_{\cf} B)^{\wedge}_{I}\to B, \;A\otimes (b_{1}\otimes b_{2})\mapsto b_{1}A(b_{2}).
\]  This implies the desired statement for the right \( U(\fg) \)-structure on \( D_{B/\cf} \).
\end{proof}

\subsection{Equivariant \texorpdfstring{\( D \)}{D}-modules}
Throughout this subsection, without mentioning otherwise, let \( \cf \) be a field over \( \QQ \), let \( G \) be an affine group scheme over \( \cf \), and let \( X \)
be a smooth (classical) scheme over \( \cf \) equipped with a right action of \( G \). Let \( \fg:=\Lie(G) \).
\begin{lem}\label{lem:t-str-D-mod-equiv}
Consider the natural maps 
\begin{equation}\label{eq:D-XdR/G}
\begin{tikzcd}
X\arrow[r,"h'"] \arrow[d,"f'"] & X^{\pdR/\cf}\arrow[d,"f"]\\
X/G \arrow[r,"h"]
& X^{\pdR/\cf}/G.
\end{tikzcd}
\end{equation}
Then \( \QCoh(X^{\pdR/\cf}/G) \) admits a t-structure such that \[ M\in \QCoh(X^{\pdR/\cf}/G)^{\ge 0}\;  (\text{resp. } \QCoh(X^{\pdR/\cf}/G)^{\le 0} ) \] if and only if \[ h^{\prime,*}f^{*}(M)\in \QCoh(X)^{\ge 0}\;  (\text{resp. } \QCoh(X)^{\le 0} ). \]
\end{lem}
\begin{proof}

By \cite[\cref{1-lemDRcomplexVSpushforward}]{Jiang2026Shla}, 
\( \QCoh(X^{\pdR/\cf}) \) is equivalent to the category of \( D_{X/\cf} \)-modules. In particular, it has a t-structure inherited from \( \QCoh(X) \) by \( h^{\prime,*} \).

By the construction in \cite[\cref{1-egScheme}]{Jiang2026Shla}, \( G\to \AnSpec(\cf) \) is proper, and so are the vertical maps in (\ref{eq:D-XdR/G}).
Thus by Barr-Beck-Lurie theorem (\cite[Theorem 4.7.3.5]{Lurie2017HA}; see also \cite[Proposition 3.1.27 (2)]{Juan2024analyticdeRham}), we have \[ \QCoh(X^{\pdR/\cf}/G)\cong \mrm{Comod}_{f^{*}f_{*}}(\QCoh(X^{\pdR/\cf})). \]

Consider the diagram \begin{equation*}
\begin{tikzcd}
X^{\pdR/\cf}\arrow[r,"f"]\arrow[d,"u"]
& X^{\pdR/\cf}/G\arrow[d,"u"]\\
\AnSpec(\cf)\arrow[r,"f_{0}"] & BG.
\end{tikzcd}
\end{equation*}
By proper base change and projection formula, \( f^{*}f_{*}(M)\cong M\otimes u^{*}(\mO(G))  \), so \( f^{*}f_{*} \) is t-exact, and \( \QCoh(X^{\pdR/\cf}/G) \) inherits a t-structure from \( \QCoh(X^{\pdR/\cf}) \) via \( f^{*} \).
\end{proof}

\begin{lem}\label{lem:Bt-gerbe-heart-fufaith}
Let \( G \) and \( X \) be as in \cref{lem:t-str-D-mod-equiv}. 
Let \( T\subset Z(G) \) be a torus with \( \mfk{t}:=\Lie(T) \), and assume that \( Y:=X/T \)
 is represented by a smooth (classical) scheme over \( \cf \). 
Let \( \bar{G}:=G/T \), which acts on \( Y \). Consider the natural morphism \[
f^{\dR}:X^{\pdR/\cf}/G\to Y^{\pdR/\cf}/\bar{G}.
\] Then \( (f^{\dR})^{*} \) is t-exact for the t-structure defined in \cref{lem:t-str-D-mod-equiv}, and restricting to the heart, \[ (f^{\dR})^{*}|_{\QCoh(Y^{\pdR/\cf}/\bar{G})^{\heartsuit}}:\QCoh(Y^{\pdR/\cf}/\bar{G})^{\heartsuit}\to \QCoh(X^{\pdR/\cf}/G)^{\heartsuit} \]
is an exact fully faithful functor of abelian categories, with the essential image consisting of \( M\in \QCoh(Y^{\pdR/\cf}/\bar{G})^{\heartsuit} \) on which the action of \( \mfk{t} \) (defined in \cref{constr:infi-char-on-eqiv-D-mod}) is trivial.
\end{lem}
\begin{proof}
The morphism 
\( f:X\to Y \) is flat, and thus \( f^{*} \) is t-exact. 
This implies that \( (f^{\dR})^{*} \) is t-exact, and \( (f^{\dR})^{*}|_{\QCoh(Y^{\pdR/\cf}/\bar{G})^{\heartsuit}} \) is exact. 

We now show that \( \mfk{t} \) acts on \( (f^{\dR})^{*}(M) \) trivially for any \( M\in \QCoh(Y^{\pdR/\cf}/\hat{G}) \). 
As in the proof of \cref{lem:explicit-Zg-action-DX}, the \( Z(\fg) \)-action can be studied on the level of \( X^{\pdR/\cf}/\hat{\fg}\cong X^{\pdR/\cf}\times B\hat{\fg} \), and thus \( Z(\mfk{t})  \)-action can be studied after pulling back to \( X^{\pdR/\cf}/\hat{\mfk{t}} \). But the composition \[ X^{\pdR/\cf}\times B\hat{\mfk{t}}\to X^{\pdR/\cf}/G\to Y^{\pdR/\cf}/\bar{G} \]
factors through \(X^{\pdR/\cf}\to Y^{\pdR/\cf} \). 


As
\( X\to Y \) is a \( T \)-torsor, so \( X^{\pdR/\cf}\to Y^{\pdR/\cf} \)
is a \( T^{\pdR/\cf} \)-torsor by \cite[\cref{1-thmAnDRStackGeneral} (1)]{Jiang2026Shla}. Thus \( f^{\dR} \) is a gerbe banded by \( B\hat{\mfk{t}} \) by \cite[\cref{1-lemRealizeLieRep}]{Jiang2026Shla}. In particular, \( f^{\dR} \) is cohomologically smooth by \cite[Proposition 4.2.5 (3)]{Juan2024analyticdeRham}. 

To show that \( (f^{\dR})^{*} \) is fully faithful, we claim that there is a natural isomorphism \[\id \cong 
\tau^{\ge 0}(f^{\dR})_{\natural}(f^{\dR})^{*}. 
\]
We have a natural morphism from the LHS to the RHS by adjunction. 
Considering the morphism \(u: B\hat{\mfk{t}}\to \AnSpec(L) \), we have \[ u_{\natural}u^{*}(\mO_{L})\cong \bigoplus_{i=0}^{\dim\mfk{t}}\wedge^{i}\mfk{t}[i]. \]
Thus we conclude that the natural morphism is an isomorphism by smooth base change and projection formula. 

Finally, we show that for any \( M\in\QCoh(X^{\pdR/\cf}/G)^{\heartsuit} \) such that \( \mfk{t} \) acts trivially, 
\( M \) lies in the essential image of \( (f^{\dR})^{*} \).
%
By localizing, we can assume that \( Y \) and \( X \) are affine. Then the category of such \( M \) is generated under colimits by \( h_{X,\natural}\mO_{X/G}\otimes_{U\mfk{t}}1 \),
where \( h_{X}:X/G\to X^{\pdR/\cf}/G \).
We claim that \[
h_{X,\natural}\mO_{X/G}\otimes_{U\mfk{t}}1\cong (f^{\dR})^{*}h_{Y,\natural}\mO_{Y/\bar{G}},
\] where \( h_{Y}:Y/\bar{G}\to Y^{\pdR/\cf}/\bar{G} \).
For this, we can consider the commutative diagram 
\begin{equation}\label{eq:rel-dR-Xg-YG}
\begin{tikzcd}
X/G\arrow[rd,"h_{X}"]\arrow[r,"h_{X/Y}"]
& X^{\pdR/Y}/G \arrow[r,"f'"]\arrow[d] 
& Y/\bar{G}\arrow[d,"h_{Y}"]\\
& X^{\pdR/\cf}/G\arrow[r,"f^{\dR}"]
& Y^{\pdR/\cf}/\bar{G},
\end{tikzcd}
\end{equation}
where the square  is Cartesian.
Then by smooth base change, it suffices to show that \[
h_{X/Y,\natural}\mO_{X}\otimes_{U\mfk{t}}1\cong \mO_{X^{\pdR/Y}/G}.
\] By adjunction, we have a natural map from the RHS to the LHS. To verify that it is an isomorphism, we can reduce to the case where \( G=T \), and \( X=T\times Y \), where the isomorphism translates into the isomorphism \( U(\mfk{t})\otimes_{U(\mfk{t})}1\cong 1 \).
\end{proof}
\begin{lem}\label{lem:BT-gerbe-fufaith}
Let \( G \) and \( X \) be as in \cref{lem:t-str-D-mod-equiv}. 
Let \( T\subset Z(G) \) be a torus with \( \mfk{t}:=\Lie(T) \), \( \bar{G}:=G/T \), and assume that the \( G \)-action on \( X \) factors through \( \bar{G} \).
Consider the natural morphism \[
f^{\dR}:X^{\pdR/\cf}/G\to X^{\pdR/\cf}/\bar{G}.
\] Then \( (f^{\dR})^{*} \) 
is fully faithful with the left inverse given by \( (f^{\dR})_{*} \).
\end{lem}
\begin{proof}
The map \( f^{\dR} \) is pulled back from \( s:BT\to \AnSpec(\cf) \). 
By \cite[Proposition 3.1.24 (2ii)]{Juan2024analyticdeRham}, \( s \) is cohomologically co-smooth. Thus by proper base change and projection formula, it suffices to show that \( s_{*}s^{*}\cong \id \), which is clear by the identification of \( \QCoh(B\GG_{m}) \)
with graded objects in \( D(\cf_{\square}) \) ({\cite{MoulinosTasos2021Tgof}, \cite[Proposition 2.2.6]{BhattLurie2022prismatic}}; see also \cite[\cref{1-corFilterdatA1modGm}]{Jiang2026Shla} for a translation of their results to the setting of solid stacks).
\end{proof}
\begin{lem}\label{lem:fix-t-action-on-torsor}
We continue with the notation of \cref{lem:Bt-gerbe-heart-fufaith}. Assume in addition that \( G\to \bar{G} \) admits a splitting as algebraic groups, that is, \( G\cong T\times \bar{G} \). 

For any algebraic character \( \lambda \) of \( T \), let \( \mcal{L}(\lambda) \)
be the line bundle on \( X/G\cong Y/\bar{G} \) defined by the Tannakian formalism. 
Then \[
h_{X,\natural}(\mO_{X/G})\otimes_{U(\mfk{t})}\lambda\cong (f^{\dR})^{*}(h_{Y,\natural}(\mcal{L}(\lambda^{-1})))\otimes s^{*}(\lambda),
\] where the morphisms are as in \[
\begin{tikzcd}
X/G\arrow[r,"f",equal]\arrow[d,"h_{X}"] & Y/\bar{G}\arrow[d,"h_{Y}"]&\\
X^{\pdR/\cf}/G\arrow[r,"f^{\dR}"] \arrow[rd,"s"]\arrow[d]
& Y^{\pdR/\cf}/\bar{G}\\
 BG\arrow[r]& BT
\end{tikzcd}
\]
\end{lem}
\begin{proof}
We consider the diagram (\ref{eq:rel-dR-Xg-YG}).
By smooth base change, it suffices to show that 
\[
h_{X/Y,\natural}(\mO_{X/G})\otimes_{U(\mfk{t})}\lambda \cong f^{\prime,*}(\mcal{L}(\lambda))\otimes s^{\prime,*}(\lambda^{-1}),
\] where \( s' \) denotes the composition \( X^{\pdR/Y}/G\to BG\to BT \).
Note that we have the  following diagram where both squares are Cartesian 
\[
\begin{tikzcd}
X/G\arrow[r,"f"]\arrow[d]
& Y/G\arrow[r]\arrow[d,"s''"]
& Y/\bar{G}\arrow[d]\\
BT
\arrow[r,"\Delta"]
& B(T\times T)\arrow[r,"p_{2}"]
\arrow[d,"p_{1}"]
& BT
\\
& BT 
\end{tikzcd}
\] and \( s'\cong p_{1}\circ s'' \),
where \( i_{1} \), \( p_{1} \)
and \( m \)
are induced by \[
\begin{tikzcd}
T\arrow[r,"{x\mapsto (x,x)}"',"\Delta"] & 
T\times T\arrow[d,"p_{1}"]
\arrow[r,"p_{2}"] & T \\
& T.
\end{tikzcd}
\]
Taking the relative de Rham stacks with respect to the horizontal maps in the square on the left, we have the following diagram where both squares are Cartesian by \cite[\cref{1-thmAnDRStackGeneral} (1)]{Jiang2026Shla}: \[
\begin{tikzcd}
X/G\arrow[r,"h_{X/Y}"]\arrow[d]
& X^{\pdR/Y}/G\arrow[r,"f'"]\arrow[d]
& Y/\bar{G}\arrow[d]\\
BT
\arrow[r,"\Delta"]
& B((\Delta(T)\subset T\times T)^{\wedge})\arrow[r,"p_{2}"]
& BT.
\end{tikzcd}
\] 

Therefore, by smooth base change, it suffices to show that \[
\Delta_{\natural}(\mO_{BT})\otimes_{U(\mfk{t})}\lambda\cong p_{1}^{*}(\lambda)\otimes p_{2}^{*}(\lambda^{-1}),
\] where the \( U(\mfk{t}) \)-action on \( \Delta_{\natural}(\mO_{BT}) \)
is induced by the first component.
For this, we have \[
(\Delta(T)\subset T\times T)^{\wedge}\cong T\times \mfk{t},\;(g_{1},g_{2})\mapsto (g_{2},g_{1}g_{2}^{-1}).
\] Then 
\begin{equation}\label{eq:transl-Delta-to-What}
\begin{tikzcd}
& B((\Delta(T)\subset T\times T)^{\wedge})\arrow[d,"\sim"]\arrow[dl,"p_{1}"]\arrow[dr,"p_{2}"]
\\
BT & B(T)\times B(\mfk{t})\arrow[r,"p_{1}"]\arrow[l,"m"] & BT,
 \end{tikzcd}
\end{equation}
where \( m \) is induced by \( m:T\times \mfk{t}\to T,\;(g_{1},g_{2})\mapsto g_{1}g_{2} \). Thus
\( \QCoh(B((\Delta(T)\subset T\times T)^{\wedge}))\cong \Rep(T)\otimes U(\mfk{t})-\Mod \),
and \( \Delta_{\natural}(\mO_{BT})\cong \mO_{BT}\boxtimes U(\mfk{t}) \) by \cite[\cref{1-lemRealizeLieRep}]{Jiang2026Shla} and smooth base change. Along this isomorphism, the \( U(\mfk{t}) \)-action on the LHS coincides with the action on \( U(\mfk{t})-\Mod \). Thus \[
\Delta_{\natural}(\mO_{BT})\otimes_{U(\mfk{t})}\lambda\cong \mO_{BT}\boxtimes \lambda
\cong m^{*}(\lambda)\otimes p_{1}^{*}(\lambda^{-1}). 
\]
The RHS translates back to \( p_{1}^{*}(\lambda)\otimes p_{2}(\lambda^{-1}) \) via (\ref{eq:transl-Delta-to-What}).
\end{proof}
\subsection{Pull-back of \texorpdfstring{\( D \)}{D}-modules along regular embeddings}
\begin{lem}\label{lem:pullback-D-mod-close-immer}
Consider an algebraic variety \( X \) over \( \cf \)
with a smooth morphism \(t: X\to \mA^{1} \), and let \( X_{0}:=X\times_{\mA^{1}}\{0\} \). Consider the diagram \[
\begin{tikzcd}
X_{0}
\arrow[r,"i"]\arrow[d,"h_{0}"]
& X\arrow[d,"h"]\\
X_{0}^{\pdR/L}
\arrow[r,"i^{\dR}"]
& X^{\pdR/L}.
\end{tikzcd}
\]
Then \( i^{\dR} \) is an open immersion, and the adjunction \( (i^{\dR}_{\natural},(i^{\dR})^{*}) \) coincides with the usual construction for \( D \)-modules. More precisely, for any \( D_{X} \)-module \( M \), \[
(i^{\dR})^{*}(M)\cong i^{*}D_{X} \otimes_{i^{-1}D_{X}}i^{-1}(M),
\] where \( i^{*}D_{X} \)
is a left \( D_{X_{0}} \)-module and a right \( i^{-1}D_{X} \)-module.

Furthermore, the right multiplication by \( t \) on \( i^{*}D_{X} \) is surjective, with kernel given by \( D_{X_{0}} \), where \( t \) is the coordinate on \( \mA^{1} \). Thus we have an ascending \( t^{\infty} \)-torsion filtration on \( (i^{\dR})^{*}(h_{0})_{\natural}(\mO_{X_{0}}) \) \[
\Fil_{i}^{t}((i^{\dR})^ {*}(h_{0})_{\natural}(\mO_{X_{0}})):=\Fib((i^{\dR})^{*}(h_{0})_{\natural}(\mO_{X_{0}})\xrightarrow{t^{i}}(i^{\dR})^{*}(h_{0})_{\natural}(\mO_{X_{0}})\cdot t^{i}),\;i\in\NN
\] with \( \gr_{i}^{t}\cong h_{\natural}(\mO_{X}) \) for \( i\in\NN \).
\end{lem}
\begin{proof}
By \cite[\cref{1-thmAnDRStackGeneral} (1)]{Jiang2026Shla}, \( i^{\dR} \)
is pulled back from \( \{0\}\to \mA^{1,\pdR/\cf} \), which can be pulled back to \( \hat{\mA}^{1}\to \mA^{1} \), which is an open immersion. Thus \( i^{\dR} \) is an open immersion. To calculate \( (i^{\dR})^{*} \), it suffices to calculate 
\( (i^{\dR})^{*}(D_{X}) \),
and \( (i^{\dR})^{*}(D_{X}) \)
has an endomorphism induced by right multiplication by \( t \).
Note that \[
h_{0}^{*}(i^{\dR})^{*}(D_{X})
\cong i^{*}h^{*}(h_{\natural}(\mO_{X}))
\cong i^{*}D_{X},
\] on which the right multiplication of \( t \) is nilpotent. 

The adjunction gives a natural morphism \( h_{0,\natural}(\mO_{X_{0}})\to (i^{\dR})^{*}h_{\natural}(\mO_{X}) \),
which corresponds to a natural \( D_{X_{0}} \)-linear map \( D_{X_{0}}\to i^{*}D_{X} \).
Moreover, \[ (i^{\dR})_{\natural}(D_{X_{0}})\cong h_{\natural}i_{\natural}(\mO_{X_{0}})\cong h_{\natural}(\mO_{X}/t)[-1]\cong D_{X}/(t)[-1]. \]
Applying \( (i^{\dR})^{*} \), this shows that the right multiplication by \( t \) on \( i^{*}D_{X} \) is surjective, with kernel given by \( D_{X_{0}} \).
\end{proof}
We will also need the following equivariant and multi-variable version:

\begin{notation}\label{notation:tuple-eg}
Let \( \Sigma \) be a finite set.

(1) For \( \unl{j}=(j_{\tau})\in\ZZ^{\Sigma} \), write \( |\unl{j}|=\sum_{\tau}j_{\tau} \).

(2) For \( \tau_{0}\in\Sigma \), let \( \delta_{\tau_{0}}\in \ZZ^{\Sigma} \), with \( (\delta_{\tau_{0}}):=\delta_{\tau_{0},\tau}\), which \( =1 \) if \( \tau_{0}=\tau \), and \( =0 \) otherwise.

(3)
For any \( w\in\ZZ \), write \( w^{\Sigma}:=\sum_{\tau\in\Sigma}w\cdot \delta_{\tau}\in\ZZ^{\Sigma} \).
\end{notation}
\begin{dfn}[Multi-filtration]\label{dfn:multi-filtr}
Let \( \mscr{C} \) be a stable category, and \( M\in\mscr{C} \). 
We endow \( \ZZ^{\Sigma} \) with a partial order \( (m_{\tau})\ge (n_{\tau}) \) if and only if \( m_{\tau}\ge n_{\tau} \) for all \( \tau \). We realize \( \ZZ^{\Sigma} \) as a category where \( \Hom((n_{\tau}),(m_{\tau})) \)
is a singleton when \( (m_{\tau})\ge (n_{\tau})\), and is empty otherwise. 

Then an \emph{ascending \( \ZZ^{\Sigma} \)-filtration} on \( M \) is the datum of a functor \( \Fil_{*}M:\ZZ^{\Sigma}\to \mscr{C} \) and an isomorphism \( \varinjlim_{\unl{j}\in\ZZ^{\Sigma}}\Fil_{\unl{j}}M\cong M \). 
For any \( \unl{j}\in\ZZ^{\Sigma} \), define \[
\gr_{\unl{j}}M:=\varinjlim_{\unl{d}\in \{-1,0\}^{\Sigma}}\Fil_{\unl{j}+\unl{d}}M,
\] where \( \{-1,0\}^{\Sigma} \)
denotes the corresponding full subcategory of \( \ZZ^{\Sigma} \).

Given an {ascending \( \ZZ^{\Sigma} \)-filtration} \( Fil_{*}M \) on \( M \), we define its total filtration (when the colimit exists) as \[
\Fil_{j}^{tot}M:=\varinjlim_{\unl{j}\in\ZZ^{\Sigma}:|\unl{j}|\le j}\Fil_{\unl{j}}M. 
\] Then the associated graded piece \[
\gr_{j}^{tot}M\cong \bigoplus_{|\unl{j}|=j}\gr_{\unl{j}}M.
\]

A \emph{descending \( \ZZ^{\Sigma} \)-filtration} is defined similarly by replacing \( \ZZ^{\Sigma} \)
with \( (\ZZ^{\Sigma})^{op} \). 
\end{dfn}
\begin{lem}\label{lem:equiv-pullback-D-mod-close-immer}
Consider an algebraic variety \( X \) over \( \cf \)
with a smooth morphism \(\mbf{t}=(t_{\tau})_{\tau\in\Sigma}^{d}: X\to (\mA^{1})^{\Sigma} \), and let \( X_{0}:=X\times_{(\mA^{1})^{\Sigma}}\{0\} \), and \( X':=X\times_{(\mA^{1})^{\Sigma}}\Delta(\mA^{1}) \), where \( \Delta:\mA^{1}\to (\mA^{1})^{\Sigma} \) is the diagonal embedding.
Write 
 \( t_{\tau} \) for the coordinate of the \( \tau \)-component of \( (\mA^{1})^{\Sigma} \), and \( t \) for that of \( \Delta(\mA^{1}) \). 

Let \( G \) be an affine group scheme over \( \cf \) that acts on \( X \). Assume that there is a character \( \chi:G\to \GG_{m} \) such that \( t:X\to (\mA^{1})^{\Sigma} \) is \( G \)-equivariant, where \( G \) acts on \( (\mA^{1})^{\Sigma} \) by composing \( \chi \) and the usual action of \( \GG_{m} \) on \( (\mA^{1})^{\Sigma} \). 

Consider the diagram \[
\begin{tikzcd}
X_{0}/G
\arrow[r,"i"]\arrow[d,"h_{0}"]\arrow[rr,bend left=12pt,"\ti{i}"] & X'/G\arrow[d,"h"]\arrow[r,"\Delta"]
& X/G\arrow[d,"\ti{h}"]\\
X_{0}^{\pdR/L}/G 
\arrow[r,"i^{\dR}"]\arrow[rr,bend right=12pt,"\ti{i}^{\dR}"] & (X')^{\pdR/\cf}/G\arrow[r,"\Delta^{\dR}"]
& X^{\pdR/L}/G.
\end{tikzcd}
\] Then 

(1) \( (\ti{i}^{\dR})^{*}\ti{h}_{\natural}(\mO_{X/G}) \) admits an ascending \( \ZZ^{\Sigma} \)-filtration \[
\Fil^{(t_{\tau})_{\Sigma}}_{\unl{j}}(\ti{i}^{\dR})^{*}\ti{h}_{\natural}(\mO_{X/G}):=
\begin{cases}
(\ti{i}^{\dR})^{*}\ti{h}_{\natural}\mrm{Kos}(t_{\tau}^{j_{\tau}+1}:\tau\in\Sigma),&\unl{j}=(j_{\tau})\in\NN^{\Sigma};\\
0,&\mrm{else},
\end{cases}
\] where 
\( \mrm{Kos}(t_{\tau}^{j_{\tau}+1}:\tau\in\Sigma) \) denotes the Koszul complex 
induced by 
\[
\bigoplus_{\tau\in\Sigma}t_{\tau}^{j_{\tau}+1}:
\mO_{X/G}\to \bigoplus_{\tau\in\Sigma} \mO_{X/G}\{j_{\tau}+1\},
\] which comes from pulling back 
\[\bigoplus_{\tau\in\Sigma}t_{\tau}^{j_{\tau}+1}:
\mO_{[\mA^{\Sigma}/\GG_{m}]}\to \bigoplus_{\tau\in\Sigma}\mO_{[\mA^{\Sigma}/\GG_{m}]}\{j_{\tau}+1\},
\]
where \( \{j\} \) denotes the twist by 
the pull-back of the \( j \)-th power of the tautological bundle on \( B\GG_{m} \).

The non-trivial graded pieces are \[
\gr^{(t_{\tau})_{\Sigma}}_{\unl{j}}(\ti{i}^{\dR})^{*}\ti{h}_{\natural}(\mO_{X/G})
\cong h_{0,\natural}(\mO_{X_{0}/G})\otimes \left\{|\unl{j}|\right\},\;\unl{j}=(j_{\tau})\in\NN^{\Sigma}.
\] 
For any \( \tau_{0}\in\Sigma \), \( t_{\tau_{0}} \) induces \[
t_{\tau_{0}}:(\ti{i}^{\dR})^{*}\ti{h}_{\natural}(\mO_{X/G})\to (\ti{i}^{\dR})^{*}\ti{h}_{\natural}(\mO_{X/G})\{1\},
\] which induces \( \gr^{(t_{\tau})_{\Sigma}}_{\unl{j}+\delta_{\tau_{0}}}\cong \gr^{(t_{\tau})_{\Sigma}}_{\unl{j}} \) for any \( \unl{j}\in\NN^{\Sigma} \). 

For any \( \unl{j}=(j_{\tau})\in\NN^{\Sigma} \), we have a Cartesian diagram 
\begin{equation}\label{eq:Fil-Cartesian-Fil0}
\begin{tikzcd}
\Fil^{(t_{\tau})_{\Sigma}}_{\unl{j}}(\ti{i}^{\dR})^{*}\ti{h}_{\natural}(\mO_{X/G})\arrow[r]\arrow[d,hook] &
h_{0,\natural}(\mO_{X_{0}/G})\{|\unl{j}|\}\arrow[d,hook,"\alpha"]
\\
(\ti{i}^{\dR})^{*}\ti{h}_{\natural}(\mO_{X/G})\arrow[r,"\unl{t}^{\unl{j}}"]
& (\ti{i}^{\dR})^{*}\ti{h}_{\natural}(\mO_{X/G})\{|\unl{j}|\},
\end{tikzcd}
\end{equation} where \( \alpha \) is induced by adjunction, and \( \unl{t}^{\unl{j}}:=\prod_{\tau}t_{\tau}^{j_{\tau}} \)
for \( \unl{j}=(j_{\tau}) \).

(2)
Consider the pro-system 
 \[
\ti{\mF}:=``\varprojlim_{\unl{m}=(m_{\tau})\in\ZZ^{\Sigma}}^{}"
(\ti{i}^{\dR})^{*}\ti{h}_{\natural}(\mO_{X/G})\left\{-|\unl{m}|\right\}
\] 
where
the transition map 
from \( \unl{m} \) to \( \unl{m}-\delta_{\tau} \) 
is induced by right multiplication by \( t_{\tau} \)as in (1).

Then \( \ti{\mF} \)
has an ascending \( t_{\tau}^{\bullet} \)-torsion \( \ZZ^{\Sigma} \)-filtration \( \{\Fil^{(t_{\tau})_{\Sigma}}_{\unl{j}}\}_{\unl{j}\in\ZZ^{\Sigma}} \), with respect to which \( \ti{\mF} \) is complete, and the graded pieces are  \[
\gr^{(t_{\tau})_{\Sigma}}_{\unl{j}}\ti{\mF}\cong h_{0,\natural}(\mO_{X_{0}/G})\left\{|\unl{j}|\right\},\;\unl{j}\in\ZZ^{\Sigma}.
\]
\end{lem}
\begin{rmk}
The filtration in (1) can be seen from smooth base change and the following Cartesian diagram \[
\begin{tikzcd}
X_{0}^{\pdR/X}/G\arrow[r]\arrow[d]
& X_{0}^{\pdR/X,+}/G\arrow[r] \arrow[d]
& X/G \times {[\mA^{1}/\GG_{m}]}
\\ X_{0}^{\pdR/\cf}/G\arrow[r]
& X_{0}^{\pdR/\cf}/G\times {[\mA^{1}/\GG_{m}]},
\end{tikzcd}
\] where the horizontal maps are induced by \( *=\GG_{m}/\GG_{m}\to \mA^{1}/\GG_{m} \). 
\end{rmk}
\begin{proof}
(1)
Consider \( i_{0}:B\GG_{m}\to \hat{\mA}^{1}/\GG_{m} \). Then we have a fiber sequence \[
i_{0,\natural}\mO_{B\GG_{m}}\to \mO_{\hat{\mA}^{1}/\GG_{m}}\xrightarrow{t} \mO_{\hat{\mA}^{1}/\GG_{m}}\{1\},
\] and thus in general, we have \[
i_{0,\natural}\mO_{B\GG_{m}}\{i\}\cong \Fib\left(\mO_{\hat{\mA}^{1}/\GG_{m}}\{i\}\xrightarrow{t} \mO_{\hat{\mA}^{1}/\GG_{m}}\{i+1\}\right),
\] which gives the description of the graded pieces by smooth base change. The more non-trivial part is to show that the natural morphism \[
\varinjlim_{j}\Fil_{j}^{t}(i^{\dR})^{*}h_{\natural}(\mO_{X/G})\to (i^{\dR})^{*}\ti{h}_{\natural}(\mO_{X/G})
\] is an isomorphism.
For this, by induction on \( \#(\Sigma) \), we can reduce to the case where \( \Sigma \) is a singleton, and the desired statement then follows from \cref{lem:pullback-D-mod-close-immer}.
(2) follows immediately from (1).
\end{proof}
\begin{lem}\label{lem:pullback-compare-diagonal}
We continue with the setting of \cref{lem:equiv-pullback-D-mod-close-immer}. Then

(1) Fix \( \tau_{0}\in\Sigma \). We have a natural injection \[
\delta:
(i^{\dR})^{*}h_{\natural}(\mO_{X'/G})\to (\ti{i}^{\dR})^{*}\ti{h}_{\natural}(\mO_{X'/G}),
\] in \( \QCoh(X_{0}^{\pdR/
\cf}/G)^{\heartsuit} \). 
Moreover, \( \delta \) identifies \( (i^{\dR})^{*}h_{\natural}(\mO_{X'/G}) \)
with \[
\mrm{Kos}(t_{\tau_{0}}-t_{\tau}:\tau\in\Sigma\backslash \{\tau\})=\left[
	(\ti{i}^{\dR})^{*}\ti{h}_{\natural}(\mO_{X'/G})
	\xrightarrow{(t_{\tau_{0}}-t_{\tau})}
	((\ti{i}^{\dR})^{*}\ti{h}_{\natural}(\mO_{X'/G})\{1\})^{\Sigma\backslash\{\tau_{0}\}}\to \cdots
\right].
\]


For any \( \tau_{0}\in\Sigma \), we have a commutative diagram \[
\begin{tikzcd}
(i^{\dR})^{*}h_{\natural}(\mO_{X'/G})\arrow[r,hook,"\delta"]\arrow[d,"t"] & (\ti{i}^{\dR})^{*}\ti{h}_{\natural}(\mO_{X'/G})\arrow[d,"t_{\tau_{0}}"]\\
(i^{\dR})^{*}h_{\natural}(\mO_{X'/G})\{1\}\arrow[r,hook,"\delta"] & (\ti{i}^{\dR})^{*}\ti{h}_{\natural}(\mO_{X'/G})\{1\}.
\end{tikzcd}
\]

(2) Consider the pro-system \[
\ti{\mF}^{0}:=``\varprojlim_{m\in\ZZ}^{}"
({i}^{\dR})^{*}{h}_{\natural}(\mO_{X'/G})\left\{-m\right\},
\] where
the transition map 
is induced by right multiplication by \( t \) as in \cref{lem:equiv-pullback-D-mod-close-immer} (1). Then we have a natural morphism \[
\delta:\ti{\mF}^{0}\to \ti{\mF},
\] 
where \( \ti{\mF} \) is as in \cref{lem:equiv-pullback-D-mod-close-immer} (2), which induces 
a natural \emph{filtered} isomorphism 
\begin{equation}\label{eq:mF0-to-Kos}
\ti{\mF}^{0}\cong \mrm{Kos}(t_{\tau}t_{\tau_{0}}^{-1}-1:\tau\in\Sigma\backslash\{\tau_{0}\}),
\end{equation}
where the LHS is endowed with \( \Fil^{t}_{j} \) defined as in \cref{lem:equiv-pullback-D-mod-close-immer} (1) for \( \Delta=\{*\} \), and \( t_{\tau}t_{\tau_{0}}^{-1} \) acts on \( \ti{\mF} \) and preserves the total filtration of \( \Fil^{(t_{\tau})_{\Sigma}}_{*} \), which induces a filtration on the RHS. 

(3)
Moreover, the action of \( t_{\tau_{1}}t_{\tau_{0}}^{-1} \) on \[
\gr^{(t_{\tau})_{\Sigma},tot}_{j}(\ti{\mF})\cong \bigoplus_{|\unl{j}|=j}\gr^{(t_{\tau})_{\Sigma}}_{\unl{j}}\ti{\mF}
\] induces \( \gr^{(t_{\tau})_{\Sigma}}_{\unl{j}}\cong \gr^{(t_{\tau})_{\Sigma}}_{\unl{j}-\delta_{\tau_{1}}+\delta_{\tau_{0}}} \). 
For any \( \unl{j}_{0}\in\ZZ^{\Sigma} \) with \( |\unl{j}_{0}|=j \), we have the following commutative diagram 
\[\begin{tikzcd}
\gr_{j}^{t}(\ti{\mF}^{0})
\arrow[r,hook,"\gr_{j}(\delta)"]
\arrow[d,"\pr_{\unl{j}_{0}}\circ\gr_{j}(\delta)"',"\cong"] & \gr_{j}^{(t_{\tau})_{\Sigma},tot}(\ti{\mF})
\cong \bigoplus_{|\unl{j}|=j}
 \gr_{\unl{j}}^{(t_{\tau})_{\Sigma}}(\ti{\mF}),\\
 \gr^{(t_{\tau})_{\Sigma}}_{\unl{j}_{0}}(\ti{\mF})\arrow[ru,"\Delta_{\unl{j}_{0}}"']
\end{tikzcd}
\] 
where the \( \unl{j} \)-component of \( \Delta_{\unl{j}_{0}} \) is induced by right multiplication by \( \unl{t}^{\unl{j}_{0}-\unl{j}}:=\prod_{\tau}t_{\tau}^{(\unl{j}_{0}-\unl{j})_{\tau}} \).
\end{lem}
\begin{rmk}
Note that \( \bigoplus_{i}M_{i}\cong ``{\prod_{i}^{}}"M_{i} \) in the pro-category, and thus 
\( \Delta_{\unl{j}_{0}} \) is well-defined.
\end{rmk}
\begin{proof}
(1): By pulling back to \( X_{0} \), and changing the coordinates on \( (\mA^{1})^{\Sigma} \) via \( (t_{\tau})\mapsto ((t_{\tau_{0}}-t_{\tau})_{\tau\ne \tau_{0}},t_{\tau_{0}}) \), we can see that \( \mrm{Kos}(t_{\tau_{0}}-t_{\tau}:\tau\in\Sigma)\in \QCoh(X_{0}^{\pdR/
\cf}/G)^{\heartsuit} \) by the surjection of right multiplication by \( t \) in \cref{lem:pullback-D-mod-close-immer}. 
The natural map \( \delta \) is induced by smooth base change, and
again by pulling back and changing the coordinates, we see that \( \delta \) factors through \( \mrm{Kos}(t_{\tau_{0}}-t_{\tau}:\tau\in\Sigma) \).
We check that it is an isomorphism by \cref{lem:pullback-D-mod-close-immer}. 
The claim on the commutative diagram can be similarly verified, which induces the map \( \delta:\ti{\mF}^{0}\to \ti{\mF} \). 

(2): Note that the action of \( t_{\tau} \) is invertible on the pro-system \( \ti{\mF} \), and thus 
\[
\delta:\ti{\mF}^{0}\cong \mrm{Kos}(t_{\tau}t_{\tau_{0}}^{-1}-1:\tau\in\Sigma\backslash\{\tau_{0}\})\hookrightarrow \ti{\mF}.
\] 
Note that the action of \( t_{\tau} \) on \( \ti{\mF} \) 
is invertible, and 
thus we have an isomorphism (\ref{eq:mF0-to-Kos}) without considering the filtration. By \cref{lem:equiv-pullback-D-mod-close-immer}, all the objects and the graded pieces lie in \( \mrm{Pro}(\QCoh(X^{\pdR/\cf}_{0}/G)^{\heartsuit}) \), so we will work in this abelian category. 

From (\ref{eq:Fil-Cartesian-Fil0}), 
we see that \( \delta \) induces \[
h_{0,\natural}(\mO_{X_{0}/G})\cong
\Fil_{0}^{t}(i^{\dR})^{*}h_{\natural}(\mO_{X'/G})\cong \Fil_{0^{\Sigma}}^{(t_{\tau})_{\Sigma}}(\ti{i}^{\dR})^{*}\ti{h}_{\natural}(\mO_{X/G}),
\] and \[
(i^{\dR})^{*}h_{\natural}(\mO_{X/G})\cap \Fil^{(t_{\tau})_{\Sigma}}_{\unl{j}}(\ti{i}^{\dR})^{*}\ti{h}_{\natural}(\mO_{X/G})=\Fil^{t}_{|\unl{j}|}(\ti{i}^{\dR})^{*}h_{\natural}(\mO_{X/G}),
\] which implies that \[
\ti{\mF}^{0}\cap \Fil^{(t_{\tau})_{\Sigma},tot}_{|\unl{j}|}\ti{\mF}=
\ti{\mF}^{0}\cap \Fil^{(t_{\tau})_{\Sigma}}_{\unl{j}}\ti{\mF}=\Fil^{t}_{|\unl{j}|}\ti{\mF}^{0}.
\] 
Then the rest follows formally from the simple fact that \[
\mrm{Kos}\left((t_{\tau}t_{\tau_{0}}^{-1}-1)|_{\gr^{(t_{\tau})_{\Sigma},tot}_{j}(\ti{\mF})}:\tau\in\Sigma\backslash\{\tau_{0}\}\right)
\] lies in \( \mrm{Pro}(\QCoh(X^{\pdR/\cf}/G)^{\heartsuit}) \), and is described as in (3) as the image of \( \Delta_{\unl{j}_{0}} \).
\end{proof}

\begin{lem}\label{lem:X0-X-epi}
Consider any solid stack \( X \) over \( \cf \)
with a morphism \(f: X\to \mA^{1} \), and let \( X_{0}:=X\times_{\mA^{1}}\{0\} \). Then 
the natural morphism \[ i':X_{0}\to X_{0}^{\pdR/X} \] is an epimorphism, and is cohomologically smooth.
\end{lem}
\begin{proof}
By \cite[\cref{1-thmAnDRStackGeneral} (1)]{Jiang2026Shla} and
\cite[Proposition 5.1.4]{Juan2024analyticdeRham},
\begin{align*}
X_{0}^{\pdR/X}
\cong X_{0}^{\pdR}\times_{X^{\pdR}}X
\cong \{0\}\times_{\mA^{1,\pdR}}X\cong 
\hat{\mA}^{1}\times_{\mA^{1}}X.
\end{align*}
Now it suffices to show that \( i:\{0\}\to \hat{\mA}^{1} \) is an epimorphism.
The natural map \( L[t]/t^{n}\to L \) has a nilpotent kernel, and thus \(
\{0\}\to \AnSpec(L[t]/t^{n})
\) is a descendable cover by \cite[Proposition. 3.35]{Mathew2016galois}, and thus a universal \( ! \)-cover by \cite[Lemma 4.7.4]{HeyerMann20246functorformalismssmoothrepresentations}. By the argument in \cite[Remark 4.2.2]{AnschützBoscoLeBrasCamargoScholze2025analyticrhamstacksfarguesfontaine}, we know that it is an epimorphism, and thus taking colimit, we know that \( i \) is an epimorphism. 
We conclude by noting that \( i \)
is cohomologically smooth, as \( \{0\}\to \mA^{1} \) is cohomologically smooth, and \( \hat{\mA}^{1}\to\mA^{1} \) is an open immersion. 
\end{proof}
\section{Vinberg semigroups and their subvarieties}\label{sec:Vinberg}
The Vinberg semigroup was introduced in \cite{Vinberg1995reductive}, \cite{Vinberg1995asymptotic}. A good summary can be found in \cite[\S 2.2]{FinkelbergKrylovMirkovicl2020DiGa}.

We work over an algebraically closed field \( C \)
of characteristic \( 0 \). Late we will take \( C:=\bar{\QQ}_{p} \). 
Let \( G \) be a reductive group over \( C \). Let \( B \) be a Borel subgroup, with a Levi subgroup \( T \) and the unipotent radical \( U \).
Let \( \Lambda:=X^{*}(T) \) and 
\(\Lambda^{+}:=X_{*}(T)^{+} \subset\Lambda \) 
the sub-monoid of \( B \)-dominant weights. Let \( \Delta:=\{\alpha_{i}:i=1,\ldots,n \} \) be the set of simple positive roots.
Let \( \Lambda^{pos}:=\sum_{i=1}^{n}\ZZ_{\ge0}\cdot \alpha_{i} \subset \Lambda \). Then \( \Lambda^{pos} \) induces a partial order \( \le \) on \( \Lambda \) such that \( \nu\le \nu' \) if and only if \( \nu'-\nu\in \Lambda^{pos} \). 
For \( \lambda\in\Lambda^{+} \), let \( V_{\lambda} \) be the irreducible algebraic representation of \( G \)
of highest weight \( \lambda \).

By Peter-Weyl Theorem, we have a \( G\times G \)-equivariant isomorphism \[
\mO(G)\cong \bigoplus_{\lambda\in \Lambda^{+}}V_{\lambda}\boxtimes V_{\lambda}^{\vee}.
\]
Then for any \( \nu\in\Lambda \),
let \[
\mO(G)_{\le\nu}:=\bigoplus_{\lambda\le \nu,\lambda\in\Lambda^{+}}V_{\lambda}\boxtimes V_{\lambda}^{\vee}\subset \mO(G).
\] 
The multiplicative structure of \( \mO(G) \) takes \[
\mO(G)_{\le \nu}\otimes \mO(G)_{\le \nu'}\to \mO(G)_{\le \nu+\nu'}.
\] In particular, \( \bigoplus_{\nu\in \Lambda}\mO(G)_{\le\nu}\cdot t^{\nu} \) is an algebra, where \( t^{\nu} \) are indeterminates satisfying the relation \( t^{\nu}\cdot t^{\nu'}=t^{\nu+\nu'} \).
\begin{dfn}[Vinberg semigroup \cite{Vinberg1995asymptotic}]\label{dfnVinbergSemigroup}
We define \[
\Vin_{G}:=\Spec\left(
  \bigoplus_{\nu\in \Lambda}\mO(G)_{\le\nu}\cdot t^{\nu}
\right).
\] 
There is a natural action of \( G\times G\times T \) on \( \Vin_{G} \), which we denote as \( *_{r}\times *_{l}\times \theta_{T} \). 
Let \[ \mA^{\Delta}:=(\mA^{1})^{\Delta}\cong\Spec\left(\bigoplus_{\nu\in\Lambda^{pos}} C\cdot t^{\nu}\right). \]
Then there is a natural map \[
\mbf{t}:\Vin_{G}\to \mA^{\Delta},
\] which is \( G\times G\times T \)-equivariant, where \( G\times G \) acts on \( \mA^{\Delta} \) trivially.

There is a non-degenerate locus \( \Vin_{G}^{*}\subset \Vin_{G} \) which is a \( G\times G\times T \)-stable Zariski open subscheme,
and which can be characterized by the property that for any \( x\in \mA^{\Delta}(C) \), \( \Vin_{G}^{*}\cap \Vin_{G}|_{x} \)
is the unique open \( G\times G \)-orbit.
\end{dfn}
Let \( Z(G) \) be the center of \( G \). By construction, \[
\Vin_{G}\times_{\mA^{\Delta}}\GG_{m}^{\Delta}\cong (G\times T)/Z(G).
\] Note that \( \GG_{m}^{\Delta}\cong T_{ad}:=T/Z(G) \)
Alternatively, there is a Tannakian definition as explained in {\cite[\S 2.2.8]{FinkelbergKrylovMirkovicl2020DiGa}}: 
\begin{lem}
    \label{lem:VinbergSemigroupTannakian}
For any (classical) \( C \)-algebra \( R \), \( \Vin_{G}(R) \) is isomorphic to the set of tuples
\[ ((g_{\lambda}\in \End_{R}(V_{\lambda}\otimes R))_{\lambda\in\Lambda^{+}},(\mbf{t}_{\nu}\in R)_{\nu\in \Lambda^{pos}}) \] satisfying the conditions:

(1) \( \mbf{t}_{0}=1 \), and \( \mbf{t}_{\bullet}:\Lambda^{pos}\to R \) is monoidal;

(2) \( g_{0}=\id \), and for any \( \lambda_{1},\lambda_{2},\lambda_{3}\in \Lambda^{+} \) with \( \lambda_{1}+\lambda_{2}\ge \lambda_{3} \), let \( i \) be the inclusion of the \( V_{\lambda_{3}} \)-isotypic subspace into \( V_{\lambda_{1}}\otimes V_{\lambda_{2}} \), and let \( \pr \) be the projection to that subspace, then \[
\mbf{t}_{\lambda_{1}+\lambda_{2}-\lambda_{3}}\cdot g_{\lambda_{3}}=g_{\lambda_{1}}\otimes g_{\lambda_{2}}.
\]
\end{lem}
\begin{proof}
We give a sketch of the argument as we failed to find a good reference. Assume that we are given \( f:\mO(\Vin_{G})\to R \).
Then for \( \lambda\in \Lambda^{+} \), consider the embedding \[ \iota_{\lambda}:V_{\lambda}\boxtimes V_{\lambda}\hookrightarrow \mO(G)_{\le \lambda}\hookrightarrow^{\cdot t^{\lambda}}\mO(\Vin_{G}),  \]
and then \( f\circ \iota_{\lambda}:V_{\lambda}\boxtimes V_{\lambda}\to R \) corresponds to a unique \( g_{\lambda}\in\End(V_{\lambda}\otimes R) \). For any \( \nu\in\Lambda^{pos} \), we put \( \mbf{t}_{\nu}:=f(t^{\nu}) \). It is straightforward to verify that the imposed conditions correspond to the Rees definition of \( \mO(\Vin_{G}) \).
\end{proof}

The construction of the following subvariety was explained to us by Arthur-C\'esar Le Bras. Let \( \mu\in X_{*}(T)^{+} \), which corresponds to maps \( \mu:\GG_{m}\to T \) and \[ \mu:\mA^{1}\to \mA^{\Delta},t\mapsto (t^{\langle \alpha,\mu\rangle})_{\alpha\in \Delta}. \]
Note that \( \langle \alpha,\mu\rangle \ge 0 \), so the map is well-defined.
\begin{dfn}\label{dfn:Vin-G-mu}
We define \[
\Vin_{G,\mu}:=\mA^{1}\times_{\mA^{\Delta}}\Vin_{G}\xrightarrow{\mbf{t}}\mA^{1}.
\] There is a natural action of \( G\times G\times \GG_{m} \) on \( \Vin_{G,\mu} \), which we denote as \( *_{r}\times *_{l}\times \theta_{\mu} \). Then \( \mbf{t} \) is \( G\times G\times \GG_{m} \)-equivariant for the trivial \( G\times G \)-action on \( \mA^{1} \).
We define 
\[ \Vin_{G,\mu}|_{0}:=\Vin_{G,\mu}\times_{\mA^{1}}\{0\},\;\Vin_{G,\mu}|_{\hat{\mA}^{1}}:=\Vin_{G,\mu}\times_{\mA^{1}}\hat{\mA}^{1}. 
\]
Then there is an action of 
\( G\times G\times \hat{\GG}_{m} \) on 
\( 
\Vin_{G,\mu}|_{\hat{\mA}^{1}} \). We denote the corresponding non-degenerate locus as 
\( \Vin_{G,\mu}^{*} \), \( \Vin_{G,\mu}^{*}|_{0} \)
and \( \Vin_{G,\mu}^{*}|_{\hat{\mA}^{1}} \).

We will refer to the action of \( \hat{\GG}_{m} \)
as the \emph{arithmetic Sen operator} \( \theta_{\mu} \).
\end{dfn}
\begin{lem}\label{lem:Vin-0-Per}
We have \( \Vin_{G,\mu}^{*}|_{0}\cong (P_{\mu}\times_{M_{\mu}}P^{\std}_{\mu})\backslash (G\times G) \);
in particular, it has a natural morphism
\[\pi^{\loc}_{\HT}\times \pi^{\loc}_{\GM}:
\Vin_{G,\mu}^{*}|_{0}\to \Fl_{G,\mu}\times \Fl^{\std}_{G,\mu},
\]
where \( \pi^{\loc}_{\HT} \) (resp. \( \pi^{\loc}_{\GM} \))
is \( G \)-equivariant for the action \( (G,*_{l}) \) (resp. \( (G,*_{r}) \)) on \( \Vin^{*}_{G,\mu}|_{0} \).
Here we are using freely the notation in \cite[Notation \ref{1-notationPartialFlagMinuscule}]{Jiang2026Shla}. 
\end{lem}
\begin{proof}
The natural map \( \mbf{t}:\Vin^{*}_{G,\mu}\to \mA^{1} \) admits a section \( \mfk{s} \) by \cite[\S 2.2.9]{FinkelbergKrylovMirkovicl2020DiGa}, and the action of \( G\times G \) is transitive on the fiber of \( \mbf{t} \). Therefore, it suffices to understand the stabilizer of the section. 
In terms of the Tannakian definition (\cref{lem:VinbergSemigroupTannakian}), for \( x\in \mA^{1} \) and any highest weight \( \lambda \), \( \mfk{s}(a) \) acts on the \( (\lambda-\nu) \)-weight space of  \( V_{\lambda} \) as \( x^{\langle \nu,\mu\rangle} \). Specializing to \( x=0 \), 
it is straightforward to see that  \[
\mrm{Stab}_{G\times G}(s(0))\subset P_{\mu}\times P_{\mu}^{\std}.
\] Moreover, the subspace
\( \im(s(0)|_{V_{\lambda}})\subset V_{\lambda} \) is a representation of \( M_{\mu} \), and \( (g_{1},g_{2})\in P_{\mu}\times P_{\mu}^{\std} \) lies in \( \mrm{Stab}_{G\times G}(s(0)) \) if and only if \[
\bar{g}_{1}|_{\im(s(0)|_{V_{\lambda}})}=\bar{g}_{2}|_{\im(s(0)|_{V_{\lambda}})},\forall\lambda\in X^{*}(T)^{+},
\] where \( \bar{g}_{i} \) denotes the image of \( g_{i} \) in \( M_{\mu} \). 
The map \( M_{\mu}\to \prod_{\lambda\in X^{*}(T)^{+}}\End(\im(s(0)|_{V_{\lambda}})) \) is injective, as one can verify by restricting to the maximal torus.  
This implies that \( \bar{g}_{1}=\bar{g}_{2} \), as desired.
\end{proof}

\begin{eg}
For \( G=\SL_{2} \), \( \Vin_{G}\cong \End(V_{(1,0)})\cong M_{2\times 2} \), the affine space of \( 2\times 2 \)-matrices, and \( \mbf{t}:\Vin_{G}\to \mA^{1} \) is induced by \( \det \). Then \( s:\mA^{1}\to \Vin_{G} \)
is given by \( t\mapsto \begin{pmatrix}
1&0\\0&t
\end{pmatrix} \). Then \( s(0)= \begin{pmatrix}
1&0\\0&0
\end{pmatrix}\), and \[
\mrm{Stab}_{G\times G}\left(\begin{pmatrix}
1&0\\0&0
\end{pmatrix}\right)=B\times_{T}\bar{B},
\] where \( B \) and \( \bar{B} \)
are the subgroups of upper-triangular and lower-triangular matrices respectively, and \( T=B\cap \bar{B} \).
\end{eg}

\begin{lem}\label{lem:Center-VinG}
(1) 
The actions \( (Z(G),*_{r}) \) and \( (Z(G),\theta_{T}) \) on \( \Vin_{G} \) coincide, and equal to the inverse of \( (Z(G),*_{l}) \). 
In other words, the action of \( (G\times G\times T,*_{l}\times *_{r}\times \theta_{T}) \) factors through \( G\times^{Z(G)}G\times^{Z(G)}T:=(G\times G\times T)/(Z(G)\times Z(G))  \), where 
\( Z(G)\times Z(G)\hookrightarrow G\times G\times T,\;(g_{1},g_{2})\mapsto (g_{1},g_{1}^{-1}g_{2},g_{2}) \).

(2)
If \( G\twoheadrightarrow G' \) is a central extension with the kernel \( Z_{0}\subset Z(G) \), then there is a natural morphism \(
\Vin_{G}/Z_{0}\to \Vin_{G'}
\), which induces \( \Vin_{G}^{*}/Z_{0}\cong \Vin_{G'}^{*} \). Here the action of \( Z_{0} \)
can be chosen to be either \( *_{l} \), \( *_{r} \), or \( \theta_{T} \), and the choice does not change the outcome by (1).
\end{lem}
\begin{proof}
(1) is straightforward from the formula in \cref{dfnVinbergSemigroup}. For (2), let us denote \( Z:=Z(G) \) and \( Z':=Z(G')=Z/Z_{0} \). Let \( d \) be the rank of \( G \) and \( G' \). Then  it suffices to show that the natural morphism \( p \) in the diagram \[
\begin{tikzcd}
\Vin_{G}^{*}/G\times^{Z}G\arrow[rr,"p"]\arrow[dr,"\mbf{t}"]
&& \Vin_{G'}^{*}/G'\times^{Z'}G'\arrow[dl,"\mbf{t}'"]\\
&\mA^{d}
\end{tikzcd}
\] is an isomorphism. 
As we see in the proof of \cref{lem:Vin-0-Per}, \( \mfk{t} \) is the classifying stack of a group scheme \( \ti{G} \) over \( \mA^{d} \), which is a closed subgroup of \( G\times^{Z} G\times \mA^{d} \) over \( \mA^{d} \). Similarly, \( \mfk{t}' \)
is the classifying stack of a group scheme \( \ti{G}\subset G'\times^{Z'}G'\times \mA^{d} \). The construction of the section in \cite[\S 2.2.9]{FinkelbergKrylovMirkovicl2020DiGa} is functorial along \( G\to G' \), and thus \( p \) is induced by a group morphism \(
\ti{G}\to \ti{G'}
\) that is compatible with the natural map \( G'\times^{Z'}G'\times \mA^{d}\twoheadrightarrow G'\times^{Z'}G'\times \mA^{d} \). 
By \cref{lem:Vin-0-Per}, both \(  \Ker(p) \)  and \(  \ti{G'}/\ti{G}  \) are schemes over \( \mA^{d} \) with a section, and become trivial when pulled back to any geometric point of \( \mA^{d} \). This implies that \( \Ker(p)\cong \ti{G'}/\ti{G}\cong \mA^{d} \) by considering the support of the augmentation ideal, and thus \( \ti{G}\cong \ti{G}' \).
\end{proof}

\section{\texorpdfstring{\( B_{\dR}^{+} \)}{BdR+}-thickenings of locally analytic Shimura varieties}\label{sec:BdR-thick}
In this section, we will show that the subvarieties \( \Vin_{G,\mu} \) of Vinberg semigroups are local models of the \( B_{\dR}^{+} \)-thickenings of locally analytic Shimura varieties introduced in \cite{Jiang2026Shla}. 
\subsection{Set-up}
Let us first recall the relevant constructions. 
We fix an identification \(\iota:\CC\cong\bar{\QQ}_{p}\).
\begin{notation}\label{notation:Shimura-SetUp}
Let \((G,X)\) be any Shimura datum, \(E\subset \CC\) be the reflex field, and \(E_{\p}
\) be the \(p\)-adic completion along the composition \(E\subset \CC\cong \bar{\QQ}_{p}\). 

For any neat open compact subgroup \(K\subset G(\mA_{f})\), let \(X_{K}\) 
be 
the associated Shimura variety defined over \(E\) (\cite{Milne1990canonical}).
We fix a cone decomposition, and let \( X_{K}^{\tor} \) be the associated toroidal compactification of \( X_{K} \) (\cite{Pink1989arithmetical}). 
We can choose the cone decomposition such that for a fixed level \(K^{p}K_{p,0}\), \(X_{K^{p}K_{p,0}}^{\tor}\) is smooth, and the boundary \(\partial X_{K}^{\tor}\) is a normal crossings divisor (\cite[Example 2.3.17]{DLLZ2019logarithmicFoundational}). 
Later, we will only consider \( X_{K}^{\tor} \) for \( K\subset K^{p}K_{p,0} \). 

For any \( p \)-adic field \(E_{\p}\subset L\subset \CC_{p} \), 
let \(\mX_{K,L}^{\tor}\)
be the analytification of \(X_{K}^{\tor}\times_{E}L\). We omit the subscript \( L \) when \( L=\CC_{p} \).
We regard them as solid stacks via \cite[Example \ref{1-egAnAdicSpace}]{Jiang2026Shla}. 
\end{notation}
From now on, we fix a neat open compact subgroup \( K^{p}\subset G(\mA_{f}^{p}) \). 
\begin{dfn}[{\cite[Definition \ref{1-dfnNonCompletedInfiniteLevel}]{Jiang2026Shla}}]\label{dfn:NonCompletedInfiniteLevel}
We define the \emph{smooth infinite-level Shimura varieties} \(\mX^{\tor,\sm}_{K^{p},L}\) of tame level \(K^{p}\) as
 \[\mX^{\tor,\sm}_{K^{p},L}:=\varprojlim_{K_{p}\subset G(\QQ_{p})}\mX^{\tor}_{\Kpp,L},\]
where the limits are taken in the category of solid stacks over \( L \) with \( K_{p} \) varying through all the open compact subgroups of \( G(\QQ_{p}) \). 
We omit the subscript \( L \) when \( L=\CC_{p} \).
Denote the natural map by \[
\pi_{K_{p}}^{\sm}:\mX^{\tor,\sm}_{K^{p},L}\to \mX^{\tor}_{K^{p}K_{p},L}.
\]
\end{dfn}

We will also consider the following construction:
\begin{construction}[{\cite[Definition \ref{1-dfnSMratConstruction}]{Jiang2026Shla}}]
For any \( E_{\p}\subset L\subset \CC_{p} \), and  \( ?\in\{(nothing),+,\hat{+}\} \)
we define the \emph{smooth log de Rham stack of Shimura varieties} as the solid stack \[(\mX^{\tor}_{K^{p},E_{\p}})^{\pdR/E_{\p},?,\sm}_{\log}:=\varprojlim_{K_{p}}(\mX^{\tor}_{\Kpp,E_{\p}})_{\log}^{\pdR/E_{\p},?},\]
and \((\mX^{\tor}_{K^{p},L})^{\pdR/L,?,\sm}_{\log}:= (\mX^{\tor}_{K^{p},E_{\p}})^{\pdR/E_{\p},?,\sm}_{\log}\times_{\AnSpec(E_{\p})}\AnSpec(L) \). 
\end{construction}

We could also consider \[\mX_{K^{p},\proket}^{\tor}:=\varprojlim_{K_{p}}\mX_{\Kpp}^{\tor},\]
where the limit is taken in the pro-Kummer \'etale site of \(\mX^{\tor}_{K^{p}K_{p,0}}\) (\cite[Definition 5.1.2]{DLLZ2019logarithmicFoundational}) with \( K_{p} \) varying through all the open compact subgroups of \( G(\QQ_{p}) \). Then \(\mX^{\tor}_{K^{p},\proket}\) is a pro-Kummer 
\'etale \(\widetilde{K_{p}}\)-torsor over \(\mX^{\tor}_{\Kpp}\). Consider the morphism between topoi: 
\[\nu:
(\mX^{\tor}_{K^{p}K_{p,0},\proket})_{/\mX^{\tor}_{K^{p},\proket}}\to |\mX^{\tor}_{K^{p},E_{\p}}|,
\] where the topological space \( |\mX^{\tor}_{K^{p},E_{\p}}|:=\varprojlim_{K_{p}}|\mX^{\tor}_{\Kpp,E_{\p}}| \). 

Let \( \hat{\mO}_{\mX_{\Kpp,\proket}^{\tor}} \) denote the completed structure sheaf of \( \mX_{\Kpp,\proket}^{\tor} \) (\cite[Definition 5.4.1 (2)]{DLLZ2019logarithmicFoundational}). 
Let \[
\hat{\mO}_{K^{p}}^{topo}:=R\nu_{*}(\hat{\mO}_{\mX^{\tor}_{K^{p},\proket}})\in D(|\mX^{\tor}_{K^{p},E_{\p}}|,\CC_{p}).
\] 
By \cite[Lemma \ref{1-lemUniformAffineProper}]{Jiang2026Shla}, \[ \QCoh(\mX^{\tor,\sm}_{K^{p},E_{\p}})\cong \varprojlim_{K_{p}}\QCoh(\mX^{\tor}_{K^{p}K_{p},E_{\p}}),\;\mF\mapsto ((\pi_{K_{p}}^{\sm})_{*}\mF)_{K_{p}}, \]
where the transition maps are given by \((-)_{*}\).
\begin{lem}\label{lem:Lift-Otopo-To-Xsm}
The sheaf \( \hat{\mO}_{K^{p}}^{topo}\in D(|\mX^{\tor}_{K^{p},E_{\p}}|,\CC_{p}) \) naturally lifts to a quasicoherent sheaf \( \hat{\mO}_{K^{p}} \) on \( \mX^{\sm}_{K^{p}} \).

More precisely, there exists \( \hat{\mO}_{K^{p}}\in \QCoh(\mX^{\sm}_{K^{p}}) \), such that for any \( K_{p}\subset G(\QQ_{p}) \) and \( U_{K_{p}}\subset \mX^{\tor}_{\Kpp} \), if \( U^{\sm} \) denotes its preimage in \( \mX^{\tor,\sm}_{K^{p}} \), and \( |U^{\sm}| \) denotes its preimage in \( |\mX^{\tor}_{K^{p},E_{\p}}| \), then \( R\Gamma(U^{\sm},\hat{\mO}_{K^{p}})\cong R\Gamma(|U^{\sm}|,\hat{\mO}^{topo}_{K^{p}}) \). 
\end{lem}
\begin{proof}
We take \( \hat{\mO}_{K^{p}} \) as in \cite[Notation \ref{1-notationOlaIla}]{Jiang2026Shla}. Then \( \hat{\mO}_{K^{p}} \) is a natural lift of \( \hat{\mO}_{K^{p}}^{topo} \) by \cite[Lemma \ref{1-lemAnalyticVSproket}]{Jiang2026Shla}. 
\end{proof}
\begin{notation}\label{notationWidetildeGroup}
We define \(Z_{c}\subset Z=Z(G)\) the maximal \(\QQ\)-anisotropic \(\RR\)-split torus, and  \(G^{c}:=G/Z_{c}\). For any subgroup \(H\subset G_{\bar{\QQ}}\) containing \(Z_{c}\), we denote \(H^{c}:=H/Z_{c}\). 

For any open subgroup \(K_{p}\subset G(\QQ_{p})\),
define \(\widetilde{K_{p}}:=K_{p}/\overline{(K^{p}K_{p}\cap Z_{c}(\QQ))}\).
Note that \(\widetilde{K_{p}}\) depends implicitly on \(K^{p}\).
\end{notation}
\begin{dfn}[{\cite[Definition \ref{1-dfnLocallyAnSV}]{Jiang2026Shla}}]\label{dfn:Ola}
We define \[ \mO^{\la}_{K^{p}}:=\varinjlim_{K_{p}\subset G(\QQ_{p})}R\Gamma(\widetilde{K_{p}\cap Z_{c}(\QQ_{p}}),\hat{\mO}_{K^{p}})^{R-K_{p}/(K_{p}\cap Z_{c}(\QQ_{p}))-\la}, \]
where \( \widetilde{K_{p}\cap Z_{c}(\QQ_{p}}) \) denotes the image of \( K_{p}\cap Z_{c}(\QQ_{p}) \) in \( \widetilde{K_{p}} \). 

Then \( \mO^{\la}_{K^{p}}\in \mrm{CAlg}(\QCoh(\mX^{\tor,\sm}_{K^{p}})) \) is analytic locally concentrated in degree \( 0 \) by \cite[Theorem \ref{1-thmOlaConcentrated0}]{Jiang2026Shla}, and we define \[
\mX^{\tor,\la}_{K^{p}}:=\unl{\AnSpec}_{\mX^{\tor,\sm}_{K^{p}}}(\mO^{\la}_{K^{p}}).
\]
\end{dfn}

\subsection{\texorpdfstring{\( \BdR^{+} \)}{BdR+}-thickenings of \texorpdfstring{\( \mX^{\tor,\la}_{K^{p}} \)}{XlaKp}}
We now construct a certain \( B_{\dR}^{+} \)-thickening \( \BdR^{+}(\mX^{\tor}_{K^{p}})^{\la} \) of \( \mX^{\tor,\la}_{K^{p}} \). This is essentially the construction in \cite[Lemma \ref{1-lemBdRXtorKp}]{Jiang2026Shla} except we are going to glue \( \AnSpec\BdR^{+,\la}(\ti{U}) \) instead of \( \AnSpec\BdR^{+}(\ti{U}) \).

Let \( E_{\p}\subset \cf\subset \CC_{p} \) be a finite extension. We will choose \( \cf \) to be large enough such that \( G \) splits over \( \cf \), and \( \mX^{\tor,\la}_{K^{p}K_{p,0},\cf} \) \'etale locally admits smooth toris charts.

We consider the morphism between topoi: 
\[\nu:
(\mX^{\tor}_{K^{p}K_{p,0},\proket})_{/\mX^{\tor}_{K^{p},\proket}}\to |\mX^{\tor}_{K^{p},\cf}|,
\]
and let \( \mbb{B}_{\dR,\mX^{\tor}_{\Kpp},\proket}^{+} \)
denote the period sheaf on \( \mX^{\tor}_{\Kpp,\proket} \) defined in \cite[Definition 6.1]{Scholze13} and \cite[Definition 2.2.3]{DLLZ2022logarithmicJAMS}.
Let 
\begin{align*}
\BdRKp^{+,topo}&:=R\nu_{*}(\mbb{B}^{+}_{\dR,\mX^{\tor}_{\Kpp},\proket})\in \widehat{\Fil}^{\NN}D(|\mX^{\tor}_{K^{p},\cf}|,\QQ_{p})
\end{align*}
which is \( \widetilde{K_{p}} \)-equivariant, and is equipped with the \( t \)-adic filtration \( t^{\bullet} \).
Here for any stable category \( \mscr{C} \), \( \widehat{\Fil}^{\NN}\mscr{C} \) denote the full subcategory of the functor category from the poset \( \NN \) to \( \mscr{C} \) consisting of those \( F \) such that \( \varprojlim_{i\in\NN} F(i)\cong 0 \).

\begin{notation}[{\cite[\cref{2-notationArithSenHT}]{Jiang2025HMF}}]\label{notation:Darith}
Let \( \cfi:=\cf(\zeta_{p^{\infty}}) \). Let \( H_{\cf}:=\Gal(\bar{\QQ}_{p}/\cfi) \),
and \( \Gamma_{\cf}:=\Gal(\cfi/\cf) \).
Let \( \mscr{C} \) be any stable category that is linear over \( D((\QQ_{p},\ZZ_{p})_{\square}) \). Let \( M\in \mscr{C} \) be a \( \Gal_{\cf} \)-equivariant object, i.e. a comodule over \( C^{0}(\Gal_{\cf},\QQ_{p}) \). Then we define \[
\Darithcf(M):=R\Gamma(H_{\cf},M)^{R-\Gamma_{\cf}-\la},
\] which is equipped with the action of the arithmetic Sen operator \( \theta\in \Lie(\Gamma_{\cf}) \).
\end{notation}
Following the convention of \cref{dfn:Ola}, let 
\begin{align*}
\BdRKp^{+,\la,topo}&:=\varinjlim_{K_{p}\subset G(\QQ_{p})}R\Gamma(\widetilde{K_{p}\cap Z_{c}(\QQ_{p}}),\BdRKp^{+,topo})^{R-K_{p}/(K_{p}\cap Z_{c}(\QQ_{p}))-\la},\\
\BdifKpcf^{+,\la,topo}&:=\Darithcf(\BdRKp^{+,\la,topo}),
\end{align*}
where the functors are taken in the filtered complete sense.

\begin{construction}\label{construction:BdR-Xla}
Let \( I \) be a finite set, and \( \{U_{K_{p,0},i}\}_{i\in I} \) be an open cover of \( \mX_{K^{p}K_{p,0},\cf}^{\tor} \) consisting of good chartable open subspaces (\cite[Definition \ref{1-dfnChartableAffOp}]{Jiang2026Shla}). 
Such cover exists up to enlarging \( \cf \) and a suitable choice of the cone decomposition by \cite[Lemma \ref{1-lemOpenSetsFromFulltoSM} (1)]{Jiang2026Shla}. 
Let \( \ti{U}_{i} \) be the preimage of \( U_{K_{p,0},i} \) in \( \mX^{\tor}_{K^{p},\cf} \).

For any \( J\subset I \), put \( \ti{U}_{J}:=\bigcap_{i\in J}\ti{U}_{j} \). Since \( \BdRKp^{+,topo} \) is filtered by Tate twists of \( \hat{\mO}_{K^{p}}^{topo} \), \( R\Gamma(\ti{U}_{J},\gr_{t^{\bullet}}^{i}\BdRKp^{+,\la,topo}) \) and \( R\Gamma(\ti{U}_{J},\gr^{i}_{t^{\bullet}}\BdifKpcf^{+,\la,topo}) \) (\( J\subset I,i\in\NN \)) are concentrated in degree \( 0 \)  by \cite[Theorem \ref{1-thmOlaConcentrated0}]{Jiang2026Shla} and \cite[Theorem 6.3.6]{Juan2022.09locallyShi}. 

We can thus apply the construction in \cite[Definition \ref{1-dfnSpfFunctor}]{Jiang2026Shla} and define \begin{align*}
\BdR^{+}(\mX^{\tor}_{K^{p}})^{\la}&:=\varinjlim_{J\subset I,J\ne \emptyset} \Spf(H^{0}(\ti{U}_{J},\BdRKp^{+,\la,topo})),\\
\Bdifcf^{+}(\mX^{\tor}_{K^{p}})^{\la}&:=\varinjlim_{J\subset I,J\ne \emptyset} \Spf(H^{0}(\ti{U}_{J},\BdifKpcf^{+,\la,topo})),\\
\mX^{\tor,\la}_{K^{p},\cfi}&:=\varinjlim_{J\subset I,J\ne \emptyset}\AnSpec(\Darithcf(H^{0}(\ti{U}_{J},\mO^{\la}_{K^{p}})))
\end{align*}
Note that by the proof of \cite[Lemma \ref{1-lemBdRTorsor} (2)]{Jiang2026Shla}, the transition maps are open immersions.
\end{construction}
\begin{lem}\label{lem:BdR-etc-Well-Defined}
The definition of the stacks \( \BdR^{+}(\mX^{\tor}_{K^{p}})^{\la} \) 
and \( \Bdifcf^{+}(\mX^{\tor}_{K^{p}})^{\la} \)
are independent of the choice of the open cover \( \{U_{K_{p,0},i}\} \).

Moreover, there are natural maps \[
\mX^{\tor,\la}_{K^{p}}\to \BdR^{+}(\mX^{\tor}_{K^{p}})^{\la}\to \Bdifcf^{+}(\mX^{\tor}_{K^{p}})^{\la}\to (\mX^{\tor}_{K^{p},\cfi})^{\pdR/\cfi,\sm}.
\]
\end{lem}
\begin{proof}
The proof is essentially the same as that of \cite[Lemma \ref{1-lemBdRXtorKp}]{Jiang2026Shla} but simpler. 
\end{proof}
The stacks \( \BdR^{+}(\mX^{\tor}_{K^{p}})^{\la} \) and \( \Bdifcf^{+}(\mX^{\tor}_{K^{p}})^{\la} \)
are related to the completed cohomology of Emerton (\cite{Emerton06}) as follows:
\begin{dfn}[\cite{Emerton06}]\label{dfn:CompleteCoh}
Define the \emph{completed cohomology} with tame level $K^{p}$ as  \[R\Gamma(K^{p},\ZZ/p^{n}):=\varinjlim_{K_{p}} R\Gamma_{\et}(X_{\Kpp,\bar{E}},\ZZ/p^{n}),
\] which carries the action of $G(\QQ_{p})\times\TT(K^{p})\times \Gal_{E}$. We also define \[R\Gamma(K^{p},\ZZ_{p}):=\Rlim_{n}R\Gamma(K^{p},\ZZ/{p^{n}}),\;R\Gamma(K^{p},\QQ_{p}):=R\Gamma(K^{p},\ZZ_{p})[\frac{1}{p}].\]
In general, for any Banach ring $C$ over $\ZZ_{p}$, $R\Gamma(K^{p},C):=R\Gamma(K^{p},\ZZ_{p})\hatotimes^{\cf}_{\ZZ_{p}}C$. We will write \(\ti{H}^{i}(K^{p},C):=H^{i}(R\Gamma(K^{p},C))\). 
	\end{dfn}
\begin{thm}\label{thm:cplt-coho-primitive-BdR}
There are natural filtered isomorphisms
\begin{align*}
(R\Gamma(K^{p},\QQ_{p})^{R-\la}\hatotimes_{\QQ_{p}}B_{\dR}^{+})^{\wedge}_{t}&\cong R\Gamma(\BdR^{+}(\mX^{\tor}_{K^{p}})^{\la}),\\
\Darithcf(R\Gamma(K^{p},\QQ_{p})^{R-\la}\hatotimes_{\QQ_{p}}B_{\dR}^{+})^{\wedge}_{t}&\cong R\Gamma(\Bdifcf^{+}(\mX^{\tor}_{K^{p}})^{\la}),
\end{align*} where \( (-)^{R-\la} \) is as in \cite[Notation \ref{1-notationRlaConvention}]{Jiang2026Shla}, \( (-)^{\wedge}_{t} \) denotes the \( t \)-adic completion, and \( \Darithcf \) is taken in the filtered complete sense.
\end{thm}
\begin{proof}
This follows from \cite[Proposition \ref{2-propPrimitiveComparisonBdR}]{Jiang2025HMF}.
\end{proof}
\begin{lem}[{\cite[Théorème 4.1]{Fontaine2004arithmetique}}]\label{lem:Fontaine0de-Rham}
	If \(V\) is a continuous Hodge-Tate representation of \(\Gal_{K}\)
	over \(\QQ_{p}\) of finite dimension.
	Let \( \cf\subset \bar{\QQ}_{p} \) be a finite extension of \( \QQ_{p} \) containing all embeddings of \( K \).
	Then \(V\) is de Rham
	if and only if the action of Fontaine operator \(\theta\in \Lie(\Gamma_{\cf})\)
	on \[\varinjlim_{i}E_{0}^{\theta}\left(\Darith(V\otimes_{\QQ_{p}}B_{\dR}^{+})^{\wedge}_{t}\cdot t^{-i}\right)\]
	is zero, where \( E_{0}^{\theta} \) as in \cref{dfn:E0-nilpotent} is taken in the filtered sense.
\end{lem} 
\begin{proof}
See for example \cite[\cref{2-lem:Fontaine0de Rham}]{Jiang2025HMF} for details.
\end{proof}



\subsection{Grothendieck-Messing theory for \texorpdfstring{\( \BdR^{+}(\mX^{\tor}_{K^{p}})^{\la} \)}{BdR+(XtorKp)la}}
In this subsection, we establish a Grothendieck-Messing theory for \( \BdR^{+}(\mX^{\tor}_{K^{p}})^{\la} \) and \( \Bdifcf^{+}(\mX^{\tor}_{K^{p}})^{\la} \). 

We take the field \( C:=\bar{\QQ}_{p} \) in \S \ref{sec:Vinberg}, and we have constructed \( \Vin_{G^{c},\mu} \) over \( \bar{\QQ}_{p} \), which descends to \( \Vin_{G^{c},\mu,\cf} \) over \( \cf \) for some large enough finite extension \( E_{\p}\subset \cf\subset \bar{\QQ}_{p} \). Recall that \( \cfi:=\cf(\zeta_{p^{\infty}}) \). In what follows, we will work over \( \cfi \): \( \mA^{1} \)
and \( \GG_{m} \) are implicitly  \( \mA^{1}_{\cfi} \)
and \( \GG_{m,\cfi} \), and absolute products are over \( (\cfi,\mO_{\cfi})_{\square} \).


\begin{thm}\label{thm:GM-theory-BdR}
Let \( Y_{\cf}:=(\Vin_{G^{c},\mu,\cfi}^{*}|_{0})^{\dR}\times_{\Fl^{\std,\dR}_{G^{c},\mu,\cfi}}
(\Fl^{\std}_{G^{c},\mu,\cfi})^{\pdR/\cfi} \).
We have the following commutative diagram of solid stacks over \( \cfi \)
\begin{equation}\label{eq:GM-theory-BdR}
\begin{tikzcd}
\mX^{\tor,\la}_{K^{p}}\arrow[d,"i'"] \arrow[r,"\pi_{\GM\HT}^{\la}"]& \Vin^{*}_{G^{c},\mu,\cfi}|_{0}/G^{c}\times \AnSpec(\CC_{p})\arrow[d,"i"]
\\
\BdR^{+}(\mX^{\tor}_{K^{p}})^{\la}
\arrow[d,"s'"] \arrow[r,"\pi_{\GM\HT}^{\BdR}"]& \Vin^{*}_{G^{c},\mu,\cfi}|_{\hat{\mA}^{1}}/G^{c}\times_{\mA^{1}}\Spf(B_{\dR}^{+})
\arrow[d,"s"]
\\
\Bdifcf^{+}(\mX^{\tor}_{K^{p}})^{\la}
\arrow[d] \arrow[r,"\pi_{\GM\HT}^{\Bdif}"]& \Vin^{*}_{G^{c},\mu,\cfi}|_{\hat{\mA}^{1}}/G^{c}
\arrow[d]
\\
(\mX^{\tor}_{K^{p},\cfi})^{\pdR/\cfi,\sm}_{\log}
\arrow[r,"\pi^{\pdR}_{\GM\HT}"]& Y_{\cf}/G^{c},
\end{tikzcd} 
\end{equation}
where all the squares are Cartesian, and the quotient by \( G^{c} \) on the RHS is taken with respect to \( *_{r} \). Here, \( \Spf \) is as in \cite[Definition \ref{1-dfnSpfFunctor}]{Jiang2026Shla}.

Moreover, the diagram is \( \Gal_{\cf}\times \hat{\fg}^{c} \)-equivariant, where \( \hat{\fg}^{c}:=(1\subset G^{c})^{\wedge} \) acts on  \( \Vin_{G^{c},\mu,\cf} \) via \( *_{l} \), and \( \Gal_{\cf} \) acts on \( \Vin_{G^{c},\mu,\cf} \) via the composition of \( \chi_{\cycl} \) with the \( \GG_{m} \)-action on \( \Vin_{G^{c},\mu,\cf} \).
\end{thm}
\begin{proof}
We start by constructing \[
\pi_{\GM\HT}^{\Bdif}:
\Bdifcf^{+}(\mX^{\tor}_{K^{p}})^{\la}\to \Vin_{G^{c},\mu,\cfi}^{*}|_{\hat{\mA}^{1}}/G^{c}. 
\]

By \cite[\cref{1-lemDRtorsorDescends}]{Jiang2026Shla}, the de Rham \( G^{c} \)-torsor \( \mcal{G}^{c}_{\dR,K} \) sends any \( V\in\Rep(G^{c}) \)
to the corresponding filtered automorphic vector bundle on \( \mX^{\tor}_{K,E_{\p}} \) with an integrable connection, and thus descends to a \( G^{c} \)-torsor \( \bar{\mcal{G}}^{c}_{\dR,K} \) on \( (\mX^{\tor}_{K,E_{\p}})^{\pdR,+}_{\log} \). 
We denote by \( \bar{\mcal{G}}^{c}_{\dR,K^{p}} \) its pull-back to \( (\mX^{\tor}_{K^{p},\cfi})^{\pdR/\cf,\sm}_{\log} \), and 
denote by \( \mcal{G}^{c,\Bdif}_{\dR,K^{p}} \) (resp. \( \mcal{G}^{c,\BdR}_{\dR,K^{p}} \), resp. \( \mcal{G}^{c,\la}_{\dR,K^{p}} \)) the pull-back of \( \bar{\mcal{G}}^{c}_{\dR,K^{p}} \)
to \( \Bdifcf^{+}(\mX^{\tor}_{K^{p}})^{\la} \) (resp. \( \BdR^{+}(\mX^{\tor}_{K^{p}})^{\la} \), resp. \( \mX^{\tor,\la}_{K^{p},\cfi} \)).
Thus it suffices to construct a \( G^{c} \)-equivariant map \[\pi^{\Bdif}_{\GM\HT}:
\mcal{G}^{c,\Bdif}_{\dR,K^{p}}\to \Vin_{G^{c},\mu,\cfi}|_{\hat{\mA}^{1}}.
\]

For this, one uses the Riemman-Hilbert correspondence (\cite{LiuZhu2017rigidity}, \cite{DLLZ2022logarithmicJAMS}): for any dominant \( \lambda \), let \( V_{\lambda}\in\Rep(G^{c}) \) be the corresponding highest weight representation, and then there is a filtered isomorphism on \( \mX^{\tor}_{K,\proket} \) \[
V_{\lambda}\otimes \mO\mbb{B}_{\dR,\mX^{\tor}_{K},\proket}\cong 
\mcal{G}_{\dR,K}^{c}(V_{\lambda})
\otimes_{\mO_{\mX^{\tor}_{K}}} \mO\mbb{B}_{\dR,\mX^{\tor}_{K},\proket}, 
\] 
that is compatible with the logarithmic integrable connections on both sides. 

Locally, restricting ourselves to \( \ti{U}_{i} \) as in \cref{construction:BdR-Xla},
and taking the (log) de Rham complexes on both sides for the connections,
we obtain a canonical modification between two \( \BdRKp^{+,\la}(\ti{U}_{i}) \)-lattices \[
R\Gamma(\nabla,\mcal{G}^{c}_{\dR,K}(V_{\lambda})\otimes_{\mO_{\mX_{K}^{\tor}}} \OBdRKp^{+,\la,topo}(\ti{U}_{i}))
\cong 
\mcal{G}^{c,\BdR}_{\dR,K^{p}}(V_{\lambda})\dashrightarrow 
V_{\lambda}\otimes \BdRKp^{+,\la,topo}(\ti{U}_{i})
\] which is bounded by \( t^{\langle \mu,\lambda\rangle} \).
Thus multiplying this morphism by \( t^{\langle \mu,\lambda\rangle} \) induces a morphism between \( \BdR^{+,\la}(\mS_{K^{p}}(G)) \)-lattices \[\mcal{G}_{\dR,K^{p}}^{c,\BdR}(V_{\lambda})\to 
V_{\lambda}\otimes \BdR^{+,\la}(\mS_{K^{p}}(G)).
\] 
Taking \( \Darithcf \) on both sides, we obtain a morphism between \( \Darithcf(\BdR^{+,\la}(\mS_{K^{p}}(G))) \)-lattices \[
\mcal{G}_{\dR,K^{p}}^{c,\Bdif}(V_{\lambda})\to V_{\lambda}\otimes \Darithcf(\BdR^{+,\la}(\mS_{K^{p}}(G))).
\] 

After pulling back to \( \mcal{G}^{c,\Bdif}_{\dR,K^{p}} \), both sides are trivialized, and then we obtain a morphism \[
\mcal{G}^{c,\Bdif}_{\dR,K^{p}}\to \Vin_{G^{c},\mu,\cfi}
\] using the Tannakian definition (\cref{lem:VinbergSemigroupTannakian}), such that the composition to \( \mA^{1} \)
is induced by \( t\in \BdR^{+}(\SGTKp{}) \). Since \( \Darithcf(\BdR^{+}(\SGTKp{})) \) is \( t \)-complete, we know that the map factors canonically through \[ \pi_{\GM\HT}^{\Bdif}: \mcal{G}^{c,\Bdif}_{\dR,K^{p}}\to \Vin_{G^{c},\mu,\cfi}|_{\hat{\mA}^{1}}. \]

By construction, restricting the morphism \( \pi^{\Bdif}_{\dR,K^{p}} \) to \( \mcal{G}^{c,\la}_{\dR,K^{p}} \), we recover the morphism \( \pi^{\la}_{\GM\HT} \) in \cite[\cref{1-thmLAstrucGMHTper}]{Jiang2026Shla}, which in particular factors through \( \Vin^{*}_{G^{c},\mu,\cfi}|_{0} \). Since the reduction of \( \mcal{G}^{c,\Bdif}_{\dR,K^{p}} \) coincides with that of \( \mcal{G}^{c,\la}_{\dR,K^{p}} \), we know that \( \pi^{\Bdif}_{\GM\HT} \) also factors through \( \Vin^{*}_{G^{c},\mu,\cfi}|_{\hat{\mA}^{1}} \).
We thus obtain the desired morphism. 

To show the commutativity, it suffices to show that the following diagram is 
commutative: 
\begin{equation}\label{eq:G-to-V}
\begin{tikzcd}
\mcal{G}^{c,\la}_{\dR,K^{p}}
\arrow[r,"\pi^{\la}_{\GM\HT}"]\arrow[d,"i'"]
& \Vin^{*}_{G^{c},\mu,\cfi}|_{0}\arrow[d,"i"]\\
\mcal{G}^{c,\BdR}_{\dR,K^{p}}
\arrow[d] \arrow[r,"\pi_{\GM\HT}^{\BdR}"]& \Vin^{*}_{G^{c},\mu,\cfi}|_{\hat{\mA}^{1}}
\arrow[d]
\\
\bar{\mcal{G}}^{c}_{\dR,K^{p}}
\arrow[r,"\pi^{\pdR}_{\GM\HT}"]& Y_{L},
\end{tikzcd}
\end{equation}

By \cref{lem:Vin-0-Per}, and \cite[\cref{1-thmLAstrucGMHTper} (1)]{Jiang2026Shla}, 
we have a Cartesian commutative diagram \[
\begin{tikzcd}
\mX^{\tor,\la}_{K^{p},\cfi}
\arrow[r]\arrow[d] & \Vin^{*}_{G^{c},\mu,\cfi}|_{0}/G^{c}\arrow[d]\\
(\mX^{\tor,\la}_{K^{p},\cfi})^{\pdR/\cfi,\sm}_{\log}
\arrow[r] & Y_{L}/G^{c}.
\end{tikzcd}
\]
In particular,
the larger diagram in (\ref{eq:G-to-V}) is commutative. 

On the other hand, the first row and the second row in (\ref{eq:G-to-V}) coincide after applying the functor \( (-)^{\pdR} \).
Since the object at the right corner coincides with its own algebraic de Rham stack,
we know that the lower square is also commutative.

It remains to show that the diagram (\ref{eq:GM-theory-BdR}) is Cartesian. By \cref{lem:Vin-0-Per} and \cite[\cref{1-thmLAstrucGMHTper} (1)]{Jiang2026Shla}, the bigger square is Cartesian. By \cref{lem:X0-X-epi}, the upper square is Cartesian with \( ! \)-surjective vertical maps.
Thus we conclude that the lower square is also Cartesian.

Finally, we note that the multiplication by \( t^{\langle \mu,\lambda\rangle} \) accounts for the twist of the Galois action on \( \Vin^{*}_{G^{c},\mu,\cfi} \). 
\end{proof}

\begin{eg}
Let us make the construction explicit in the case of modular curves. We consider the \emph{homological convention} with \( \mu=(1,0) \). With this normalization, the standard representation of \( G=\GL_{2,\QQ} \) corresponds to the Tate module of the universal elliptic curves \( \mcal{E} \). 
By Riemann-Hilbert correspondence (\cite{LiuZhu2017rigidity}, \cite{DLLZ2022logarithmicJAMS}), we have a natural map \[
T_{p}\mcal{E}\otimes \BdRKp^{+}\hookrightarrow 
\mcal{G}^{\BdR}_{K^{p}}(V_{(1,0)}):=
R\Gamma(\nabla,D_{K_{p}}\otimes_{\mO_{\mX_{\Kpp}^{\tor}}} \OBdR^{+})
\] where \( D_{K_{p}} \) is the first relative de Rham homology of \( \mcal{E} \) with respect to \( \mX_{\Kpp}^{\tor} \). Multiplying the ``inverse'' of this map with \( t=t^{\langle \lambda,\mu\rangle } \) induces a morphism \[
\mcal{G}^{\BdR}_{K^{p}}(V_{(1,0)})\hookrightarrow 
T_{p}\mcal{E}\otimes \BdRKp^{+}.
\] 
As we will explain in \S \ref{subsec:eg-GL2}, \( \Vin_{G,\mu}\cong M_{2\times 2}\times \GG_{m} \), and this construction defines a map from \( \BdR^{+}(\mX^{\tor}_{K^{p}})^{\la} \) into the \( M_{2\times 2} \)-component. One replace \( V_{(1,0)} \) with \( V_{(1,1)} \) to obtain a map into the \( \GG_{m} \)-component.
\end{eg}
 This enables us to study the period sheaf \( \BdRKp^{+} \) using differential operators on \( \Vin_{G^{c},\mu} \).
 
\begin{cor}\label{cor:BdR-Base-Change}
Put \[
\bar{\mX}^{\tor,\la}_{K^{p},\cfi}:=(\mX^{\tor}_{K^{p},\cf})^{\pdR/\cf,\sm}\times_{Y_{\cf}}(\Vin^{*}_{G^{c},\mu,\cfi}|_{0})^{\pdR/\cfi}/G^{c}.
\]
Then there is a Cartesian diagram over \( \cfi \) \begin{equation}\label{eq:Alg-DR-Bdif}
\begin{tikzcd}
\Bdifcf^{+}(\mX^{\tor}_{K^{p}})^{\la}\arrow[d,"h^{\Bdif}"]\arrow[r,"\pi_{\GM\HT}^{\Bdif}"] & \Vin^{*}_{G^{c},\mu,\cfi}|_{\hat{\mA}^{1}}/G^{c}\arrow[d,"\hat{h}_{0}"]\arrow[r,"\hat{i}"]
& \Vin^{*}_{G^{c},\mu,\cfi}/G^{c}\arrow[d,"h"]\\
\bar{\mX}^{\tor,\la}_{K^{p},\cfi}\arrow[r,"\bar{\pi}_{\GM\HT}"] & (\Vin^{*}_{G^{c},\mu,\cfi}|_{0})^{\pdR/\cfi}/G^{c}
\arrow[r,"i^{\dR}"] & (\Vin^{*}_{G^{c},\mu,\cfi})^{\pdR/\cfi}/G^{c},
\end{tikzcd}
\end{equation}
where \( \hat{h}_{0} \) is cohomologically smooth, and we have a natural 
\( \hat{\fg}^{c} \)-equivariant filtered
isomorphism \[
i:h^{\Bdif}_{\natural}(\mO_{\Bdifcf^{+}(\mX^{\tor}_{K^{p}})^{\la}})\cong \bar{\pi}_{\GM\HT}^{*}\circ \hat{h}_{0,\natural}(\mO_{\Vin^{*}_{G^{c},\mu,\cfi}|_{\hat{\mA}^{1}}/G^{c}})\cong \bar{\pi}_{\GM\HT}^{*}(i^{\dR})^{*}h_{\natural}(\mO_{\Vin_{G^{c},\mu,\cfi}^{*}/G^{c}}),
\] where the LHS is endowed with a descending \( t \)-adic filtration, and the RHS is endowed with a descending filtration induced by the ascending \( t^{\bullet} \)-torsion filtration on \( (i^{\dR})^{*}h_{\natural}
(\mO_{\Vin_{G^{c},\mu,\cfi}^{*}/G^{c}}) \) in \cref{lem:equiv-pullback-D-mod-close-immer}.

Moreover, along \( i \),
the action of \( \theta\in\Lie(\Gamma_{\cf}) \)
on the LHS coincides with the action of \( \theta_{\mu} \) (as in \cref{dfn:Vin-G-mu}) on the RHS.
\end{cor}
\begin{proof}
It suffices to show that \( \hat{h}_{0} \) is cohomologically smooth. The rest follows from the smooth base change.
By \cite[\cref{1-lemBasicFiltered}]{Jiang2026Shla},
\[ \Vin^{*}_{G^{c},\mu,\cfi}|_{0}\to(\Vin^{*}_{G^{c},\mu,\cfi}|_{0})^{\pdR/\cfi} \] is cohomologically smooth. 
By \cref{lem:X0-X-epi}, \[
\Vin^{*}_{G^{c},\mu,\cfi}|_{0}\to \Vin^{*}_{G^{c},\mu,\cfi}|_{\hat{\mA}^{1}}
\] is a suave \( ! \)-cover. Thus \( \hat{h}_{0} \) is also cohomologically smooth by \cite[Corollary 8.6]{Mann20226functorformalismmathbbzell}.
\end{proof}
\begin{notation}
Denote \[\BdifKpcf^{+,\la}:=h^{\Bdif}_{*}(\mO_{\Bdifcf^{+}(\mX^{\tor}_{K^{p}})^{\la}})\in\QCoh(\bar{\mX}^{\tor,\la}_{K^{p},\cfi}).\]
\end{notation}
We think of \( \BdifKpcf^{+,\la} \)
as a lift of \( \BdifKpcf^{+,\la,topo}\in D(|\mX^{\tor}_{K^{p},\cf}|) \) to \( \bar{\mX}^{\tor,\la}_{K^{p},\cfi} \); the latter can be thought of as a suitable version of the algebraic de Rham stack of \( \mX^{\tor,\la}_{K^{p}} \). 

Now by adjunction and smooth base change, we have  
\begin{align}\label{align:Bdif-D-pi-i-h}
\BdifKpcf^{+,\la}&\cong \mbb{D}_{\bar{\mX}^{\tor,\la}_{K^{p}}}\circ h^{\Bdif}_{\natural}(\mO_{\Bdifcf^{+}(\mX^{\tor}_{K^{p}})^{\la}})\\
&\cong \mbb{D}_{\bar{\mX}^{\tor,\la}_{K^{p}}}\circ \bar{\pi}_{\GM\HT}^{*}\circ \hat{h}_{0,\natural}(\mO_{\Vin^{*}_{G^{c},\mu,\cfi}|_{\hat{\mA}^{1}}/G^{c}})\\
&\cong \mbb{D}_{\bar{\mX}^{\tor,\la}_{K^{p}}}\circ \bar{\pi}_{\GM\HT}^{*}\circ (i^{\dR})^{*}\circ h_{\natural}(\mO_{\Vin_{G^{c},\mu,\cfi}^{*}/G^{c}}).
\end{align}

The isomorphism (\ref{align:Bdif-D-pi-i-h}) is equivariant for the action of \( \fg^{c} \). We then obtain the following corollary.
\begin{cor}\label{cor:BdR-lambda-D-pi-GMHT}
Let \( Z(\fg^{c}):=Z(U(\fg^{c})) \). Then
for any infinitesimal character \( \lambda:Z(\fg^{c})\to \cfi \), we have a filtered isomorphism 
\begin{align}\label{align:BdR-lambda-D-pi-GMHT}
i:\BdifKpcf^{+,\la,\lambda}
\cong \mbb{D}_{\bar{\mX}^{\tor,\la}_{K^{p}}}\circ \bar{\pi}_{\GM\HT}^{*}\circ (i^{\dR})^{*}\left(h_{\natural}
(\mO_{\Vin_{G^{c},\mu,\cfi}^{*}/G^{c}})\otimes_{Z(\fg^{c})}\lambda\right),
\end{align}
where the LHS is endowed with a descending \( t \)-adic filtration, and the RHS is endowed with a descending filtration induced by the ascending \( t^{\bullet} \)-torsion filtration on \( (i^{\dR})^{*}h_{\natural}
(\mO_{\Vin_{G^{c},\mu,\cfi}^{*}/G^{c}}) \) in \cref{lem:equiv-pullback-D-mod-close-immer}, the subscript \( \lambda \) refers to \( \RHom_{Z(\fg^{c})}(\lambda,-) \),
and \( Z(\fg^{c}) \)-action on \( h_{\natural}
(\mO_{\Vin_{G^{c},\mu,\cfi}^{*}/G^{c}}) \)
is induced by \( (G^{c},*_{l}) \) via \cref{constr:infi-char-on-eqiv-D-mod}.

Moreover, \( i \) is equivariant for the action of \( \theta\in\Lie(\Gamma_{L}) \)
on the LHS and the action of \( \theta_{\mu} \) on the RHS.
\end{cor}
\begin{dfn}\label{dfn:E0-nilpotent}
Let \( \mscr{C} \) be a stable category.
For any \( M\in\mscr{C} \) equipped with an endomorphism \( \theta \), denote \[
E^{\theta}_{0}(M):=\varinjlim_{n}\Fib(M\xrightarrow{\theta^{n}}M),\;E^{0}_{\theta}(M):=\varprojlim_{n} \mrm{Cofib}(M\xrightarrow{\theta^{n}}M),
\] where the transition maps in the two systems are induced by \[
\begin{tikzcd}
(M\arrow[r,"\theta^{n}"]\arrow[d,"\id"] & M)\arrow[d,"\theta"]\\
(M\arrow[r,"\theta^{n+1}"] & M)
\end{tikzcd}\text{ and }
\begin{tikzcd}
(M\arrow[r,"\theta^{n+1}"]\arrow[d,"\theta"] & M)\arrow[d,"\id"]\\
(M\arrow[r,"\theta^{n}"] & M)
\end{tikzcd}
\] respectively. In particular, we have natural transformations \[
E_{0}^{\theta}\to \id \to E^{\theta}_{0}.
\]

More generally, if \( M\in\mscr{C} \) is equipped with commuting endomorphism \( (\theta_{\tau})_{\tau\in\Sigma} \) for a finite set \( \Sigma=\{\tau_{1},\ldots,\tau_{m}\} \), we define \[
E_{0}^{(\theta_{\tau})}=E_{0}^{\theta_{\tau_{1}}}\circ E_{0}^{\theta_{\tau_{2}}}\circ \cdots\circ E_{0}^{\theta_{\tau_{m}}},\;
E^{0}_{(\theta_{\tau})}=E^{0}_{\theta_{\tau_{1}}}\circ E^{0}_{\theta_{\tau_{2}}}\circ \cdots\circ E^{0}_{\theta_{\tau_{m}}}.
\] It is clear that the resulting functors are independent of choice of the ordering of \( \Sigma \). 
\end{dfn}
\begin{rmk}
When \( \mscr{C}=D(R) \), \( (M,\theta) \)
corresponds to a quasicoherent sheaf on \( \mA^{1}_{R} \). Consider \( j:\hat{\mA}^{1}_{R}\hookrightarrow \mA^{1}_{R} \). Then the natural transformations above correspond to \[
j_{!}j^{*}M\to M\to j_{*}j^{*}M.
\] 
\end{rmk}
\begin{cor}\label{cor:BdR-lambda-E0D-pi-GMHT}
For any infinitesimal character \( \lambda:Z(\fg)\to \cfi \), we have an isomorphism 
\begin{align*}
i:E_{0}^{\theta}(\BdifKpcf^{+,\la,\lambda})
\cong \mbb{D}_{\bar{\mX}^{\tor,\la}_{K^{p},\cfi}}\circ \bar{\pi}_{\GM\HT}^{*}\circ E^{0}_{\theta_{\mu}}\circ (i^{\dR})^{*}\left(h_{\natural}
(\mO_{\Vin_{G^{c},\mu,\cfi}^{*}/G^{c}})\otimes_{Z(\fg)}\lambda\right),
\end{align*}
where both functors \( E_{0}^{\theta} \)
and \( E^{0}_{\theta_{\mu}} \) are taken in the filtered category, the subscript \( \lambda \) refers to \( \RHom_{Z(\fg)}(\lambda,-) \), and \( Z(\fg^{c}) \)-action on \( h_{\natural}
(\mO_{\Vin_{G^{c},\mu,\cfi}^{*}/G^{c}}) \)
is induced by \( (G^{c},*_{l}) \).

Moreover, \( i \) is equivariant for the action of \( \theta\in\Lie(\Gamma_{L}) \)
on the LHS and the action of \( \theta_{\mu} \) on the RHS. Here 
\end{cor}
\begin{proof}
By \cref{cor:BdR-lambda-D-pi-GMHT},
we have 
\begin{align*}
&\Fib(\BdifKpcf^{+,\la,\lambda}\xrightarrow{\theta^{n}}\BdifKpcf^{+,\la,\lambda})\\
&\cong \mbb{D}_{\bar{\mX}^{\tor,\la}_{K^{p},\cfi}}\circ \bar{\pi}_{\GM\HT}^{*}\circ \left((i^{\dR})^{*}\left(h_{\natural}
(\mO_{\Vin_{G^{c},\mu,\cfi}^{*}/G^{c}})\otimes_{Z(\fg)}\lambda\right)/^{L}\theta_{\mu}^{n}\right),
\end{align*}
and we obtain the desired isomorphism by taking inverse limit. 
\end{proof}
The study of the nilpotent action of \( \theta \) on \( E^{0}_{\theta}(\BdifKpcf^{+,\la,\lambda}) \) is closely related to Fontaine operators, and thus related to classicality of classes appearing in the completed cohomology via \cref{thm:cplt-coho-primitive-BdR}. 
\cref{cor:BdR-lambda-E0D-pi-GMHT} allows us to reduce the problem to a calculation on the local model \( \Vin^{*}_{G^{c},\mu} \).
\section{Fontaine operators on Vinberg semigroups}\label{sec:Font-op-on-Vin}
In this section, we will study in the cases of \( \GL_{2} \) and of \( \mrm{Res}_{F/\QQ}\GL_{2} \) the action of \( \theta_{\mu} \) on \[
E_{0}^{\theta_{\mu}}\circ (i^{\dR})^{*}\left(h_{\natural}
(\mO_{\Vin_{G^{c},\mu,\cfi}^{*}/G^{c}})\otimes_{Z(\fg)}\lambda\right).
\]
By \cref{cor:BdR-lambda-E0D-pi-GMHT}, this is related to Fontaine operators on Shimura varieties. For this reason, we will refer to this action as Fontaine operators on Vinberg semigroups.

\subsection{Example: \texorpdfstring{\( G=\GL_{2} \)}{GL2}}\label{subsec:eg-GL2}
In this section, we specify to the case of \( G=\GL_{2,\cf} \), and \( \mu:t\mapsto \mrm{diag}(t,1) \), and study the action of \( \theta_{\mu} \) on the \( D_{\Vin_{G,\mu,\cf}^{*}|_{0}/L} \)-module \[ E_{0}^{\theta_{\mu}}\circ (i^{\dR})^{*}\left(h_{\natural}
(\mO_{\Vin_{G^{c},\mu,\cfi}^{*}/G^{c}})\otimes_{Z(\fg)}\lambda\right) \] for a (regular algebraic) infinitesimal character \( \lambda:Z(\fg)\to \cf \). As we will explain in \S \ref{sec:Fontaine-Ope-HMV}, by \cref{cor:BdR-lambda-E0D-pi-GMHT}, this recovers the formula \( N_{k}=c_{k}\cdot d^{k}\circ \bar{d}^{k} \) in \cite[Theorem 6.2.6]{Pan2209.06II}.

In this section, we will work with schemes, and work over \( \cf=\QQ \), and we will omit \( \cf \) from the subscripts for simplicity.

Put \( G:=\GL_{2,\QQ} \), and \( \mu:t\mapsto \mrm{diag}(t,1) \). We fix the Borel \( B\subset G \) consisting of upper-triangular matrices, and its Levi subgroup \( T \) consisting of diagonal matrices.
In this case, we have \[
\Vin_{G}\cong \Vin_{G,\mu}\cong M_{2\times 2}\times \GG_{m}\supset \Vin_{G}^{*}\cong M_{2\times 2}^{*}\times \GG_{m},
\]
where \( M_{2\times 2} \) is the affine space of \( 2\times2 \)-matrices,
corresponds to the endomorphism space of the standard representation \(\mrm{Std}= V_{(1,0)} \), \( M_{2\times 2} \) is the subspace of non-zero matrices,
and \( \GG_{m} \) corresponds to that of the determinant \( \det=V_{(1,1)} \). 
The natural map to \( \mA^{1}=\mA^{\Delta} \)
is given by
\[
\mbf{t}:M_{2\times 2}\times \GG_{m}\to \mA^{1},(X,a)\mapsto a^{-1}\det(X). 
\] Thus \[
\Vin_{G}^{*}|_{0}\cong M_{2\times 2}^{\rk=1}\times \GG_{m}.
\]

The action of \( (G\times G\times T ,*_{l}\times *_{r}\times \theta_{T}) \) on \( \Vin_{G} \)
is given by \begin{align*}
G\times G\ni (g_{1},g_{2}):&(X,a)\mapsto (g_{1}^{-1}Xg_{2},a\det(g_{1})^{-1}\det(g_{2}))\\
T\ni \begin{pmatrix}
t_{1} & 0\\0 & t_{2}
\end{pmatrix}:&(X,a)\mapsto (t_{1}X,t_{1}t_{2}a).
\end{align*}
Let \( x_{ij} \)
for \( 1\le i,j\le 2 \)
be the coordinate for \( M_{2\times 2} \), 
and \( a \)
be that for \( \GG_{m} \).

Let \[ z=
\begin{pmatrix}
1 & 0\\0& 1
\end{pmatrix},
h=\begin{pmatrix}
1 & 0\\0& -1
\end{pmatrix},e=\begin{pmatrix}
0&1\\0&0
\end{pmatrix},f=\begin{pmatrix}
0&0\\1&0
\end{pmatrix}\in\mfk{gl}_{2}. \]
\( z \)
and \( h \)
also form a basis for \( \mfk{t} \), and we will denote by \( \theta \)
the action of \( \mfk{t} \)
induced by that of \( T \).

We will consider the following elements in \( D_{\Vin_{G}} \): 
\begin{align*}
\theta_{\mu}&:=\theta_{T}(d\mu),\\
\Omega&=\frac{1}{2}h^{2}+h+2fe\in Z(\fg),\\
\mrm{Eu}&:=\sum_{1\le i,j\le 2}x_{ij}\partial_{x_{ij}}=\theta(h),\\
\det(\partial)&:=\det \begin{pmatrix}
\partial_{x_{11}} & \partial_{x_{12}}\\
\partial_{x_{21}} & \partial_{x_{22}}
\end{pmatrix}.
\end{align*}
Then we have by the following \emph{Capelli identity} in \( D_{\Vin_{G}} \) (\cite{Capelli1887zuruckfuhrung}) 
\begin{align}\label{align:capelli}
\Omega=\frac{1}{2}\mrm{Eu}^{2}+\mrm{Eu}-2\det(X)\det(\partial).
\end{align}
See also for example \cite[(1.12)]{HoweUmeda1991TCit} or \cite[Lemma 4.14]{BenZviGanev2022wonderful}.
In particular, \( \Omega \), \( \theta_{\mu} \),
\( \mrm{Eu} \)
and \( \det(X)\det(\partial) \)
commute with each other.



Note that
\( P_{\mu}=B \) (resp. \( P_{\mu}^{\std}=\bar{B} \)) the subgroup of upper-triangular (resp. lower-triangular) matrices.
We have natural maps 
\[
\Fl^{std}_{G,\mu}:=\bar{B}\backslash G\cong \mP^{1}\xleftarrow{\pi_{\GM}^{\loc}}\Vin_{G}^{*}|_{0}\xrightarrow{\pi_{\HT}^{\loc}}\Fl_{G,\mu}:=B\backslash G\cong \mP^{1}.
\]

We will following the convention in \cite{Jiang2025HMF}:
\begin{notation}[{\cite[\cref{2-notationSomeAlgRep}]{Jiang2025HMF}}]\label{notation:weight-rep}
We define a \emph{weight} as a pair \((k,w)\in\ZZ^{\oplus2}\) with \(k\equiv w\pmod{2}\). We say that \( (k,w) \) is a \emph{regular} weight if \( k\ne 1 \).

(1) If  \((k,w)\) is regular,
denote by \(V^{(k,w)}\) the algebraic representation of \(G\) \[V^{(k,w)}:=\mrm{Sym}^{\max\{-k,k-2\}}\mrm{Std}\otimes\det{}^{\frac{-w-\max\{-k,k-2\}}{2}},\] 
and denote by \(\lambda^{(k,w)}:Z(\fg)\to \QQ\) its infinitesimal character.
Concretely, \( \lambda^{(k,w)}(z)=-w \)
and \( \lambda^{(k,w)}(\Omega)=\frac{1}{2}k(k-2). \)

(2)
Denote by \(\chi^{(k,w)}:=\left(\frac{w-k}{2},\frac{w+k}{2}\right)\) the algebraic representation of \(H\); that is,
 \(\chi^{(k,w)}(h):=-k\),
and \(\chi^{(k,w)}(z):=w\).
We then define \[
\omega^{(k,w)}_{\Fl^{\std}_{G,\mu}}:=\mcal{P}^{\std}_{\mu}(\chi^{(k,w)})\in\mrm{Pic}(\Fl^{\std}_{G,\mu}/G),\;\omega^{(k,w)}_{\Fl_{G,\mu}}:=\mcal{P}_{\mu}(\chi^{(-k,w)})\in\mrm{Pic}(\Fl_{G,\mu}/G),
\] where \( \mcal{P}_{\mu}^{\std} \) and \( \mcal{P}_{\mu} \) are torsors as in \cite[\cref{1-notationDRtorsor}]{Jiang2026Shla}.
See \cite[\cref{2-rmk:normalization-kw}]{Jiang2025HMF} for a discussion of the normalization here. 
\end{notation}
The following is the main theorems of this section:
\begin{notation}\label{notation:twist-T-action}
For any solid stack \( X \) over \( BT \), \( \chi\in\QCoh(BT) \), and \( M\in \QCoh(X) \), we will write \[
M\{\chi\}:=M\otimes s^{*}(\chi),
\] where \( s \) denotes the structure map \( X\to BT \).
\end{notation}
\begin{thm}\label{thm:Vin-GL2-Gr-t}
Let \( \ti{G}:=G\times G\times T \).
Consider the diagram \[
\begin{tikzcd}[column sep=huge]
\Vin_{G}^{*}/\ti{G}\arrow[d,"h"]
& \Vin_{G}^{*}|_{0}/\ti{G}\arrow[l,"i"']\arrow[r,"{(\pi_{\GM}^{\loc}, \pi_{\HT}^{\loc})}"]\arrow[d,"h_{0}"] & \Fl_{G,\mu}^{\std}/G\times \Fl_{G,\mu}/G\arrow[d,"{(h_{\Fl}^{\std},h_{\Fl})}"]\\
(\Vin_{G}^{*})^{\pdR}/\ti{G}
& \arrow[l,"i^{\dR}"](\Vin^{*}_{G}|_{0})^{\pdR}/\ti{G}
\arrow[r,"{(\pi_{\GM}^{\loc,\dR},\pi_{\HT}^{\loc,\dR})}"] & \Fl_{G,\mu}^{\std,\pdR}/G\times \Fl_{G,\mu}^{\pdR}/G.
\end{tikzcd}
\]
Let \( (k,w) \) be a regular weight with \( k\ge 2 \). We will make free use of \cref{notation:twist-T-action}. 

(1)
Consider the pro-system 
 \[
\ti{\mF}_{(k,w)}:=``\varprojlim_{m\in\ZZ}^{}"
(\ti{i}^{\dR})^{*}\ti{h}_{\natural}(\mO_{\Vin_{G,\cf}^{*}/\ti{G}})\otimes_{Z(\fg)}\lambda^{(k,w)}\left\{\chi^{(2m,0)}\right\},
\] where 
the transition maps are induced by right multiplication by \( t \) as in \cref{lem:equiv-pullback-D-mod-close-immer} (2), and the action of \( Z(\fg) \) is induced by \( (\fg,*_{l}) \) via \cref{constr:infi-char-on-eqiv-D-mod}.

Then \( \ti{\mF}_{(k,w)} \)
has an ascending \( t^{\bullet} \)-torsion filtration \( \{\Fil^{t}_{j}\} \), with respect to which \( \ti{\mF}^{0}_{\kw} \) is complete, and the graded pieces are  \[
\gr^{t}_{j}\ti{\mF}_{(k,w)}\cong h_{0,\natural}(\mO_{\Vin_{G,\cf}^{*}|_{0}/\ti{G}})\otimes_{Z(\fg)}\lambda^{(k,w)}\left\{\chi^{(-2j,0)}\right\},\;j\in\ZZ,
\] 
and 
\( (\mfk{t},\theta_{T}) \) acts on \( \gr^{t}_{*}(\ti{\mF}_{(k,w)}) \) semisimply with integral eigen-characters.

(2)
The natural morphism 
\[ E_{0}^{\theta_{\mu}}\left(\ti{\mF}_{\kw}\right)\to E^{0}_{\theta_{\mu}}\left(\ti{\mF}_{\kw}\right) \]
is an isomorphism, 
where \( E^{\theta_{\mu}}_{0} \) and \( E_{\theta_{\mu}}^{0} \) are as in \cref{dfn:E0-nilpotent} taken in the filtered sense with respect to \( \Fil^{t}_{*} \).

Moreover, the pro-system \( E_{0}^{\theta_{\mu}}\left(\ti{\mF}_{\kw}\right) \)
is eventually constant for sufficiently small \( m \).
We denote this limit as \( \mF_{(k,w)} \). Then
\( \mF_{(k,w)} \) has an ascending \( t^{\bullet} \)-torsion filtration \( \Fil^{t}_{*} \), with \[
\gr^{t}_{j}\cong 
\begin{cases}
	(\pi^{\loc,\dR}_{\GM})^{*} (h_{\Fl}^{\std})_{\natural}(\omega^{(-k,-w)}_{\Fl^{\std}_{G,\mu}})\otimes (\pi^{\loc,\dR}_{\HT})^{*} (h_{\Fl})_{\natural}(\omega^{(-k,w)}_{\Fl_{G,\mu}})\{\chi^{(w,w)}\},&j=\frac{k-w}{2};\\
	(\pi^{\loc,\dR}_{\GM})^{*} (h_{\Fl}^{\std})_{\natural}(\omega^{(k-2,-w)}_{\Fl^{\std}_{G,\mu}})\otimes (\pi^{\loc,\dR}_{\HT})^{*} (h_{\Fl})_{\natural}(\omega^{(k-2,w)}_{\Fl_{G,\mu}})\{\chi^{(w,w)}\},&j=\frac{2-k-w}{2};\\
	0,&\text{else}.
\end{cases}
\] 

(3) The (right) action of \( \theta_{\mu} \) preserves the filtration,
is zero on \( \gr_{j}^{t}(\mF_{(k,w)}) \), and its action on \( \mF_{(k,w)} \)
induces a monodromy operator \[
\ti{N}_{k}:\gr^{t}_{\frac{k-w}{2}}(\mF_{(k,w)})\to \gr^{t}_{\frac{2-k-w}{2}}(\mF_{(k,w)}),
\]
which after twisted by \( \{\chi^{-(w,w)}\} \) naturally descends along \( \pi^{\loc,\dR}_{\GM}\times \pi^{\loc,\dR}_{\HT} \)
to \[
N_{k}:
(h^{\std}_{\Fl})_{\natural}(\omega_{\Fl^{\std}_{G,\mu}}^{(-k,-w)})
\boxtimes
(h_{\Fl})_{\natural}(\omega^{(-k,w)}_{\Fl_{G,\mu}})
\to 
(h^{\std}_{\Fl})_{\natural}(\omega_{\Fl^{\std}_{G,\mu}}^{(k-2,-w)})\boxtimes
(h_{\Fl})_{\natural}(\omega^{(k-2,w)}_{\Fl_{G,\mu}}).
\] 
\end{thm}
\begin{thm}\label{thm:Vin-GL2-Fontaine=Diff}We continue with the notation in \cref{thm:Vin-GL2-Gr-t}.
Then we have an equality \[
N_{k}=c_{k}\cdot \theta^{k-1}_{\Fl^{\std}_{G,\mu}}
\boxtimes \theta^{k-1}_{\Fl_{G,\mu}},
\] where \( c_{k}\in\QQ^{\times} \), and
\begin{align*}
\theta^{k-1}_{\Fl^{\std}_{G,\mu}}&:(h^{\std}_{\Fl})_{\natural}(\omega_{\Fl^{\std}_{G,\mu}}^{(-k,-w)})\to (h^{\std}_{\Fl})_{\natural}(\omega_{\Fl^{\std}_{G,\mu}}^{(k-2,-w)}),\\
\theta^{k-1}_{\Fl_{G,\mu}}&:(h_{\Fl})_{\natural}(\omega_{\Fl_{G,\mu}}^{(-k,w)})\to (h_{\Fl})_{\natural}(\omega_{\Fl_{G,\mu}}^{(k-2,w)}),
\end{align*}
correspond via \cite[\cref{1-lemDiffOpeToLogStack} (1)]{Jiang2026Shla} to the theta operators on \( \Fl^{\std}_{G,\mu} \)
and \( \Fl_{G,\mu} \) respectively.
\end{thm}
\begin{rmk}
It is possible to calculate the constant \( c_{k} \)
and compare with that of \cite[Theorem 1.3]{Qiu2026fontaineoperatorcuspsmodular}.
\end{rmk}
The proof requires some preparation:

\begin{prop}\label{prop:Zg-gr_t-VinG}
Let \( (k,w) \) be a regular weight and \( \ti{G}:=G\times G\times T \). Then \[
h_{0,\natural}(\mO_{\Vin_{G}^{*}|_{0}/\ti{G}})\otimes_{Z(\fg)}\lambda^{(k,w)}
\cong N_{-k}\oplus N_{k-2},
\]
where the \( Z(\fg) \)-action is induced by \( (\fg,*_{l}) \) by \cref{constr:infi-char-on-eqiv-D-mod} and \[
N_{m}\cong 
(\pi^{\loc,\dR}_{\GM})^{*} (h_{\Fl}^{\std})_{\natural}(\omega^{(m,-w)}_{\Fl^{\std}_{G,\mu}})\otimes (\pi^{\loc,\dR}_{\HT})^{*} (h_{\Fl})_{\natural}(\omega^{(m,w)}_{\Fl_{G,\mu}})\left\{\chi^{(m,-w)}\right\},
\]
 and \( (\mfk{t},\theta) \)
acts on \( N_{m} \) via \( \chi^{(-m,w)} \) for \( m\in\{-k,k-2\} \). 
\end{prop}
\begin{proof}
By \cref{lem:pullback-D-mod-close-immer}, \( D_{\Vin_{G}^{*}|_{0}}\subset i^{*} D_{\Vin_{G}^{*}} \) is the subsheaf consisting of elements annihilated by right multiplication by \( t \).
Then by the Capelli identity (\ref{align:capelli}), the right action of \( \Omega \) on \( D_{\Vin_{G}^{*}|_{0}} \) equals that of \( \frac{1}{2}\Eu(\Eu+2) \). 
By \cref{lem:explicit-Zg-action-DX}, this implies that \[
\Omega=\frac{1}{2}\Eu(\Eu+2)\in H^{0}\End(h_{0,\natural}(\mO_{\Vin^{*}_{G}|_{0}})),
\] where the actions of \( \Omega \) and \( \Eu:=\theta_{T}(1,-1) \) on \( h_{0,\natural}(\mO_{\Vin^{*}_{G}|_{0}}) \) are defined as in \cref{constr:infi-char-on-eqiv-D-mod}. 
Thus
\begin{align*}
h_{0,\natural}(\mO_{\Vin_{G}^{*}|_{0}/\ti{G}})\otimes_{Z(\fg)}\lambda^{(k,w)}&\cong D_{\Vin_{G}^{*}|_{0}}/(\theta(z)-w,(\Eu+k)(\Eu-k+2))
\\
&\cong N_{-k}
\oplus N_{k-2},
\end{align*}
where 
\[
N_{m}:= 
h_{0,\natural}(\mO_{\Vin_{G}^{*}|_{0}/\ti{G}})\otimes_{U(\mfk{t})}\chi^{(-m,w)},
\; m\in\{-k,k-2\}.
\]
Note that \( (\mfk{t},\theta) \)
comes from an algebraic action \( (T,\theta) \)
on \( \Vin_{G}^{*}|_{0} \), with \( \Vin_{G}^{*}|_{0}/T\cong \Fl_{G,\mu}^{\std}\times \Fl_{G,\mu} \).
Concretely, we have \( T \)-torsors \( \mcal{P}_{\mu}\to \Fl_{G,\mu} \)
and \( \mcal{P}^{\std}_{\mu}\to \Fl_{G,\mu}^{\std} \), 
and 
\begin{align}\label{align:VinG0-torsor-Mod-T}
\Vin_{G}^{*}|_{0}\cong (\mcal{P}_{\mu}^{\std}\times \mcal{P}_{\mu})/T,
\end{align} 
where \( T \) acts on \( \mcal{P}_{\mu}^{\std}\times \mcal{P}_{\mu} \)
diagonally, and along (\ref{align:VinG0-torsor-Mod-T}),
\( (T,\theta) \) is induced by the action of \( T \) on \( \mcal{P}_{\mu}^{\std} \). Thus applying \cref{lem:fix-t-action-on-torsor}, we obtain 
\begin{align*}
N_{m}
&\cong (\pi^{\loc,\dR}_{\GM}\times \pi^{\loc,\dR}_{\HT})^{*}\circ (h^{\std}_{\Fl}\times h_{\Fl})_{\natural}(\omega_{\Fl^{\std}_{G,\mu}}^{(m,-w)}\boxtimes \omega_{\Fl_{G,\mu}}^{(m,w)})\{\chi^{(-m,w)}\}\\
&\cong 
(\pi^{\loc,\dR}_{\GM})^{*} (h_{\Fl}^{\std})_{\natural}(\omega^{(m,-w)}_{\Fl^{\std}_{G,\mu}})\otimes (\pi^{\loc,\dR}_{\HT})^{*} (h_{\Fl})_{\natural}(\omega^{(m,w)}_{\Fl_{G,\mu}})\{\chi^{(-m,w)}\},
\end{align*}
where we have used the K\"unneth formula  (\cite[\cref{1-lemKunnethProper}]{Jiang2026Shla}), and \( \{\chi^{(-m,w)}\} \) is as in \cref{notation:twist-T-action}.
\end{proof}
\begin{proof}[Proof of \cref{thm:Vin-GL2-Gr-t}]
The existence of such a filtration follows immediately from \cref{lem:equiv-pullback-D-mod-close-immer} (2). 
By \cref{prop:Zg-gr_t-VinG},
\[
\gr_{j}^{t}\ti{\mF}_{(k,w)}\cong h_{0,\natural}(\mO_{\Vin_{G}^{*}|_{0}/\ti{G}})\otimes_{Z(\fg)}\lambda^{(k,w)}\left\{\chi^{(-2j,0)}\right\}\cong N_{-k}^{(j)}\oplus N_{k-2}^{(j)},\;j\in\NN,
\] with \[
N_{m}^{(j)}\cong 
(\pi^{\loc,\dR}_{\GM})^{*} (h_{\Fl}^{\std})_{\natural}(\omega^{(m,-w)}_{\Fl^{\std}_{G,\mu}})\otimes (\pi^{\loc,\dR}_{\HT})^{*} (h_{\Fl})_{\natural}(\omega^{(m,w)}_{\Fl_{G,\mu}})\left\{\chi^{(-m-2j,w)}\right\},\;m\in\{-k,k-2\}.
\]
Thus by \cref{lem:Bt-gerbe-heart-fufaith}, \( (\mfk{t},\theta) \)
acts on \( N_{m}^{(j)} \) via \( \chi^{(-m-2j,w)} \), and \( E_{0}^{\theta_{\mu}}(N_{m}^{(j)})\cong 0 \) unless \( j\in\{\frac{k-w}{2},\frac{2-k-w}{2}\} \). Thus \( E_{0}^{\theta_{\mu}}(N_{m}^{(j)}) \) is eventually constant for \( m<\min(\frac{k-w}{2},\frac{2-k-w}{2}) \),
and \( \mF_{(k,w)} \) is equipped with an induced filtration with graded pieces described as in the statement of \cref{thm:Vin-GL2-Gr-t}.

Finally, the existence of the monodromy operator \( \ti{N}_{k} \) is provided by the t-structure defined in \cref{lem:t-str-D-mod-equiv}, and the descent of \( \ti{N}_{k} \) along \( \pi^{\loc,\dR}_{\GM}\times \pi^{\loc,\dR}_{\HT} \) follows from \cref{lem:Bt-gerbe-heart-fufaith}.
\end{proof}
\begin{proof}[Proof of \cref{thm:Vin-GL2-Fontaine=Diff}]
By \cite[\cref{1-lemDiffOpeToLogStack} (1)]{Jiang2026Shla}, 
\( N_{k} \) corresponds to a \( G\times G \)-equivariant differential operator from 
\( \omega^{(2-k,w)}_{\Fl_{G,\mu}^{\std}}\boxtimes \omega^{(2-k,-w)}_{\Fl_{G,\mu}} \)
to 
\( \omega^{(k,w)}_{\Fl_{G,\mu}^{\std}}\boxtimes \omega^{(k,-w)}_{\Fl_{G,\mu}} \) on \( \Fl^{\std}_{G,\mu}\times \Fl_{G,\mu} \). The space of such differential operators is of dimension \( 1 \) over \( \QQ \), generated by \( \theta^{k-1}_{\Fl_{G,\mu}^{\std}}\boxtimes \theta^{k-1}_{\Fl_{G,\mu}} \). 
Thus there exists \( c_{k}\in\QQ \), such that \[
N_{k}=c_{k}\cdot \theta^{k-1}_{\Fl^{\std}_{G,\mu}}\boxtimes \theta^{k-1}_{\Fl_{G,\mu}}.
\]
Now it remains to show that \( c_{k}\ne 0 \). 

For this, we pull back the construction along \( (\Vin_{G}^{*})^{\pdR}\to( \Vin_{G}^{*})^{\pdR}/\ti{G} \), base change from \( \QQ \) to \( \QQ_{p} \), and localize to the dagger neighborhood of the closed locus \( \{x_{21}=1,a=1\} \). Put 
\begin{align*}
B:=\QQ_{p}\{a-1,x_{21}-1\}^{\dagger}[x_{11},x_{22},x_{21}]\twoheadrightarrow 
A:=\QQ_{p}\{a-1,x_{21}-1\}^{\dagger}[z_{1},z_{2}]\cong B/t,
\\
x_{21}\mapsto x_{21},
x_{11}\mapsto z_{1}x_{21},x_{22}\mapsto z_{2}x_{21},x_{12}\mapsto z_{1}z_{2}x_{21},a\mapsto a,
\end{align*} with \( t=\det(x_{ij}) \), and \( \QQ_{p}\{x\}^{\dagger}:=\varinjlim_{m}\QQ_{p}\left\langle\frac{1}{p^{m}}x\right\rangle \).
Then \[X_{0}:=
\AnSpec\left(A
\right)\cong 
X\times_{\Vin^{*}_{G}}\Vin^{*}_{G}|_{0}
\hookrightarrow X:=\AnSpec(B)\subset \Vin^{*}_{G}.
\] Put \( \hat{B}^{+}:=B^{\wedge}_{t} \), and \( \hat{B}:=\hat{B}^{+}[1/t] \). Then \( \hat{B}^{+} \)
and \( \hat{B} \)
are equipped with \( t \)-adic filtrations \( \Fil^{\bullet}_{t} \).

We have maps from \( (\Vin_{G}^{*}|_{0})_{x_{21}} \) to \( U_{1}:=\Spec(\QQ_{p}[z_{1}])\subset\Fl_{G,\mu}^{\std}\cong \mP^{1} \)
and \( U_{2}:=\Spec(\QQ_{p}[z_{2}])\subset \Fl_{G,\mu}\cong \mP^{1} \) induced by \( z_{1} \)
and \( z_{2} \) 
respectively. Then we have \( (\Vin_{G}^{*}|_{0})_{x_{21}}\cong U_{1}\times U_{2}\times T_{0} \), for \( T_{0}:=\Spec(\QQ_{p}\{x_{21}-1,a-1\}^{\dagger}) \).

Consider the extension in \cref{thm:Vin-GL2-Gr-t} \[
0\to \gr^{t}_{\frac{2-k-w}{2}}(\mF_{(k,w)})
\to \mF_{(k,w)}
\to \gr^{t}_{\frac{k-w}{2}}(\mF_{(k,w)})\to 0.
\] 
Note that by smooth base change, \[
(i^{\dR})^{*}h_{\natural}(\mO_{\Vin_{G}^{*}})\cong \hat{h}_{0,\natural}(\mO_{\Vin_{G}^{*}|_{\hat{\mA}^{1}}}), 
\] where \( \hat{h}_{0} \)
is the natural morphism \( \hat{h}_{0}:\Vin_{G}^{*}|_{\hat{\mA}^{1}}\to (\Vin_{G}^{*}|_{0})^{\pdR} \).
Thus by adjunction, \[ (\hat{h}_{0,\natural}(\mO_{\Vin_{G}^{*}|_{\hat{\mA}^{1}}}))^{\vee}\cong \hat{h}_{0,*}(\mO_{\Vin_{G}^{*}|_{\hat{\mA}^{1}}}), \] and 
\begin{align*}
\mF_{(k,w)}^{\vee}\cong \varinjlim_{m\in\ZZ} E^{0}_{\theta_{\mu}}(t^{-m}\cdot \hat{h}_{0,*}(\mO_{\Vin_{G}^{*}|_{\hat{\mA}^{1}}})^{\lambda^{(k,w)}}),
\end{align*}
 where \( (-)^{\lambda^{(k,w)}}:=\RHom_{Z(\fg)}(\lambda^{(k,w)},-) \), and \( (-)^{\vee} \) is as in \cref{notation:Naive-dual}.  
\( \mF_{(k,w)}^{\vee} \) has an induced descending \( t \)-adic filtration \( \Fil_{t^{\bullet}}^{*} \), with graded pieces 
\[
\gr_{t^{\bullet}}^{j}\cong 
\begin{cases}
	(\pi^{\loc,\dR}_{\GM}\times \pi^{\loc,\dR}_{\HT})^{*}(h^{\std}_{\Fl}\times h_{\Fl})_{*}(\omega^{(k,w)}_{\Fl^{\std}_{G,\mu}}\boxtimes \omega^{(k,-w)}_{\Fl_{G,\mu}}),&j=\frac{k-w}{2};\\
	(\pi^{\loc,\dR}_{\GM}\times \pi^{\loc,\dR}_{\HT})^{*}(h^{\std}_{\Fl}\times h_{\Fl})_{*}(\omega^{(2-k,w)}_{\Fl^{\std}_{G,\mu}}\boxtimes \omega^{(2-k,-w)}_{\Fl_{G,\mu}}),&j=\frac{2-k-w}{2};\\
	0,&\text{else}.
\end{cases}
\]
Thus we have \[ R\Gamma((X^{\pdR},\mF_{(k,w)}^{\vee})\cong E^{0}_{\theta_{\mu}}(\hat{B}^{\lambda^{(k,w)}}):=``\varinjlim_{m}" E^{0}_{\theta_{\mu}}(t^{-m}\cdot \hat{B}^{+,\lambda^{(k,w)}}), \]
with the direct limit system being eventually constant.
Then \( E^{0}_{\theta_{\mu}}(\hat{B}^{\lambda^{(k,w)}}) \)
has an induced \( t \)-adic filtration with the only non-trivial graded pieces being \[
\gr_{t^{\bullet}}^{\frac{i-w}{2}}\cong R\Gamma(U_{1}\times U_{2},\omega^{(i,w)}_{\Fl^{\std}_{G,\mu}}\boxtimes \omega^{(i,-w)}_{\Fl_{G,\mu}})\cong H^{0}(U_{1}\times U_{2},\omega^{(i,w)}_{\Fl^{\std}_{G,\mu}}\boxtimes \omega^{(i,-w)}_{\Fl_{G,\mu}})[0]
\] for \( i\in\{k,2-k\} \).
Thus \[ E^{0}_{\theta_{\mu}}(\hat{B}^{\lambda^{(k,w)}})\cong E^{0}_{\theta_{\mu}}((t^{\frac{2-k-w}{2}}\hat{B}^{+}/t^{\frac{2+k-w}{2}}\hat{B}^{+})^{\lambda^{(k,w)}} ),\]
and is concentrated in degree \( 0 \). It suffices to show that the action of \( \theta_{\mu} \) on \( E^{0}_{\theta_{\mu}}(\hat{B}^{\lambda^{(k,w)}}) \) is non-zero.


We put \( u(x):=\frac{1}{2}(x+k)(x-k+2) \), and \[
\Omega_{k}:=\Omega-\frac{1}{2}k(k-2)=u(\Eu)-2t\det(\partial).
\] We have \[ \gr^{\frac{2-k-w}{2}}_{t^{\bullet}}\cong\left \{\bar{F}=t^{\frac{2-k-w}{2}}a^{\frac{2-k+w}{2}}x_{21}^{k-2}f(z_{1})g(z_{2}):f(z_{1})\in \QQ_{p}[z_{1}],g\in\QQ_{p}[z_{2}] \right\}.\]
Any such function \( \bar{F} \) lifts to \[
F_{0}:=t^{\frac{2-k-w}{2}}G_{0}\in t^{\frac{2-k-w}{2}}\hat{B}^{+},\;\text{with }G_{0}:=a^{\frac{2-k+w}{2}}x_{21}^{k-2}f(z_{1})g(z_{2})\in\hat{B}^{+}.
\] Then \( \Eu(G_{0})=(k-2)G_{0} \), and \[
\Omega_{k}(F_{0})=t^{\frac{2-k-w}{2}}\Omega_{k}(G_{0})=2t^{\frac{4-k-w}{2}}a^{\frac{2-k+w}{2}}x_{21}^{k-2}f'(z_{1})g'(z_{2}).
\]
We remedy this defect by considering \[
F_{1}:=\sum_{i=0}^{k-2}\frac{(-2)^{i}}{\prod_{j=1}^{i}u(k-2-2j)}t^{\frac{2-k-w}{2}+i}a^{\frac{2-k+w}{2}}x_{21}^{k-2}f^{(i)}(z_{1})g^{(i)}(z_{2}).
\]
Then \( F_{1} \) is another lift of \( \bar{F} \), \( \Omega_{k}(F_{1})=
t^{\frac{k-w}{2}}G_{1} \), with \[G_{1}
:= \frac{(-2)^{k-1}}{\prod_{j=1}^{k-2}u(k-2-2j)}a^{\frac{2-k+w}{2}}x_{21}^{k-2}f^{(k-1)}(z_{1})g^{(k-1)}(z_{2})\in\hat{B}^{+}.
\] We have \( \Eu(G_{1})\cong (k-2)G_{1} \) and \( (\mfk{t},\theta) \) acts on \( F_{1} \) by \( (0,w) \).

Consider \( \ln(x_{21}):=\sum_{n=}^{\infty}(1-x_{21})^{n}/n\in\QQ_{p}\{x_{21}-1\}^{\dagger} \).
Then 
\( \Eu(\ln(x_{21}))=1 \), and \( (\Eu-2+k)(G_{1}\ln(x_{21}))=G_{1} \). Thus
\begin{align*}
\Omega_{k}(G_{1}\ln(x_{21}))&\equiv u(\Eu)(G_{1}\ln(x_{21}))\\&\equiv \frac{1}{2}(\Eu+k)G_{1}\equiv (k-1)G_{1}
\pmod{t}.
\end{align*}
Now put \[
F:=F_{1}-\frac{1}{k-1}t^{\frac{k-w}{2}}G_{1}\ln(x_{21}),
\] and then \( F \) is another lift of \( \bar{F} \) such that \( \Omega_{k}(F)\in t^{\frac{2+k-w}{2}}\hat{B}^{+} \). In other words, \( F \)
is a lift of \( \bar{F} \)
in \( E^{0}_{\theta_{\mu}}((t^{\frac{2-k-w}{2}}\hat{B}^{+}/t^{\frac{2+k-w}{2}}\hat{B}^{+})^{\lambda^{(k,w)}} ) \). Thus 
\begin{align*}
N_{k}(\bar{F})&=\theta_{\mu}(F)=-\frac{1}{k-1}t^{\frac{k-w}{2}}G_{1}\theta_{\mu}(\ln(x_{21}))=-\frac{1}{k-1}t^{\frac{k-w}{2}}G_{1}\\
&=\frac{2}{(k-1)((k-2)!)^{2}}t^{\frac{k-w}{2}}a^{\frac{2-k+w}{2}}x_{21}^{k-2}f^{(k-1)}(z_{1})g^{(k-1)}(z_{2})\in\hat{B}^{+}.
\end{align*} In particular, \( N_{k}\ne 0 \), as desired.
\end{proof}

\subsection{Example: \texorpdfstring{\( G=\Res_{F/\QQ}\GL_{2} \)}{ResF/QGL2} and \texorpdfstring{\( G^{c} \)}{Gc}}\label{subsec:eg-ResF-GL2}
In this section, we will take \( F \) to be a number field, and \( G:=\Res_{F/\QQ}\GL_{2} \). Put \( Z_{c}:=\{g\in\Res_{F/\QQ}Z(\GL_{2}):\Nm_{F/\QQ}(g)=1\} \), where \( Z(\GL_{2})\cong \GG_{m} \) is the center of \( \GL_{2} \), and \( \Nm_{F/\QQ}:\Res_{F/\QQ}Z(\GL_{2})\to Z(\GL_{2}) \), and \( G^{c}:=G/Z_{c} \). 

We take  \( \cf\subset \bar{\QQ} \) to be the Galois closure of \( F \) over \( \QQ \),
and let \( \Sigma:=\{\tau:F\hookrightarrow \cf\} \).
Then \( G_{\cf}\cong \prod_{\tau\in\Sigma}\GL_{2,\cf} \), and fix a Borel \( B=\prod_{\tau}B_{\tau} \) (resp. a maximal torus \( T=\prod_{\tau}T_{\tau} \), resp. \( \bar{B}=\prod_{\tau}\bar{B}_{\tau} \)), where \( B_{\tau}\subset\GL_{2,\cf} \) (resp. \( T_{\tau} \), resp. \( \bar{B}_{\tau} \)) is the subgroup of upper-triangular (resp. diagonal, resp. lower-triangular) matrices. 

For \( \tau_{0}\in\Sigma \), we put  \[
\mu_{\tau_{0}}:\GG_{m,\cf}\to G_{\cf},\;t\mapsto \left(
	\begin{pmatrix}
	t & 0\\0 &1
	\end{pmatrix}_{\tau=\tau_{0}},
	\begin{pmatrix}
	1 & 0\\0 &1
	\end{pmatrix}_{\tau\in\Sigma\backslash \{\tau_{0}\}}\right).
\] Let \( \mu:=\prod_{\tau\in\Sigma}\mu_{\tau} \).
By abuse of notation, we also denote by \( \mu \) the composition of \( \mu \) with \( G\to G^{c} \). We will study Fontaine operators on \( \Vin_{G,\mu} \), and \( \Vin_{G^{c},\mu} \),
which corresponds to the Hilbert modular varieties via \cref{cor:BdR-lambda-E0D-pi-GMHT}. 
\begin{rmk}
It is also possible to consider a variant \( \mu_{T}:=\prod_{\tau\in T^{c}}\mu_{\tau} \), which relates to the quaternionic Shimura varieties considered in \cite{Jiang2025HMF}. 
\end{rmk}

Since we have \( G_{\cf}\cong \prod_{\tau}\GL_{2,\cf} \), we have \[
\Vin_{G,\cf}\cong \prod_{\tau\in\Sigma}\Vin_{\GL_{2},\cf}\supset \Vin_{G,\cf}^{*}\cong \prod_{\tau\in\Sigma}\Vin_{\GL_{2},\cf}^{*}. 
\] The set of positive simple roots 
is \( \Delta=\{\alpha_{\tau}\}_{\tau\in\Sigma} \), where \( \alpha_{\tau} \) is \( (1,-1) \) at the \( \tau \)-component and \( (0,0) \) at the other components. 
We will identify \( \Delta\cong \Sigma \).
The natural map \( \mbf{t}:\Vin_{G,\cf}\to (\mA^{1})^{\Sigma} \) is induced by \[
\mbf{t}_{\tau}:\Vin_{\GL_{2},\cf}\to \mA^{1},\;\tau\in \Sigma,
\]  and \( \mbf{t}_{\tau} \) is as described in \S \ref{subsec:eg-GL2}. Then \( \mu:\mA^{1}\to (\mA^{1})^{\Sigma} \) as considered in \cref{dfn:Vin-G-mu} is precisely the diagonal embedding, and we have a Cartesian diagram over \( \cf \) \[
\begin{tikzcd}
\Vin_{G,\cf}^{*}|_{0}\arrow[r,"i"]\arrow[d,"\mbf{t}"] &
\Vin_{G,\mu,\cf}^{*}\arrow[r,"\Delta"]\arrow[d,"\mbf{t}"]
& \Vin^{*}_{G,\cf}\cong (\Vin^{*}_{\GL_{2},\cf})^{\Sigma}\arrow[d,"\mbf{t}=(\mbf{t}_{\tau})_{\tau}"]\\
\{0\}\arrow[r,"i"] &
\mA^{1}\arrow[r,"\Delta"]
& (\mA^{1})^{\Sigma}.
\end{tikzcd}
\]
We realize \( \GG_{m} \) as a subgroup of \( T \) via \( \mu \). 
Let \( \ti{G}:=G\times G\times T\supset \ti{G}_{0}=G\times G\times \GG_{m} \).
We consider the following commutative diagram \begin{equation}\label{eq:Pull-ti-G0}
\begin{tikzcd}
\Vin_{G,\cf}^{*}|_{0}/\ti{G}_{0}\arrow[r,"i"]\arrow[rr,"\ti{i}",bend left=12pt]\arrow[d,"h_{0}"] &
\Vin_{G,\mu,\cf}^{*}/\ti{G}_{0}\arrow[r,"\Delta"]\arrow[d,"h"]
& \Vin_{G,\cf}^{*}/\ti{G}_{0}
\arrow[d,"\ti{h}"]\\
(\Vin_{G,\cf}^{*}|_{0})^{\pdR}/\ti{G}_{0}\arrow[r,"i^{\dR}"]\arrow[rr,"\ti{i}^{\dR}",bend right=12pt] &
(\Vin_{G,\mu,\cf}^{*})^{\pdR}/\ti{G}_{0}\arrow[r,"\Delta^{\dR}"]
& (\Vin_{G,\cf}^{*})^{\pdR}/\ti{G}_{0},
\end{tikzcd}
\end{equation}
where the larger square is pulled back from \begin{equation}\label{eq:Res-Vin-G-basic}
\begin{tikzcd}
\Vin_{G,\cf}^{*}|_{0}/\ti{G}\arrow[r,"\ti{i}"]\arrow[d,"h_{0}"]
& \Vin_{G,\cf}^{*}/\ti{G}
\arrow[d,"\ti{h}"]\\
(\Vin_{G,\cf}^{*}|_{0})^{\pdR}/\ti{G}\arrow[r,"\ti{i}^{\dR}"]
& (\Vin_{G,\cf}^{*})^{\pdR}/\ti{G},
\end{tikzcd}
\end{equation} along \( B\ti{G}_{0}\to B\ti{G} \). Now the diagram (\ref{eq:Res-Vin-G-basic}) is a product over \( \tau\in\Sigma \) of the diagram in \cref{thm:Vin-GL2-Gr-t}.
Note that (\ref{eq:Pull-ti-G0}) is in the setting of \cref{lem:equiv-pullback-D-mod-close-immer}. 

\begin{dfn}[{\cite[\cref{2-dfn:multiweight}\& \cref{2-notationSomeAlgRep}]{Jiang2025HMF}}]
	\label{dfn:weight-HMF}
We say that a tuple \(\kw\)
	is a \emph{multiweight}, if \(\Bbbk=(k_{\tau})_{\tau\in\Sigma}\) 
	with \(k_{\tau},w\in\ZZ\),
	and \(k_{\tau}\equiv w\pmod{2}\). We say that \( \kw \) is \emph{regular} if \( k_{\tau}\ne 1 \) for all \( \tau\in\Sigma \).

Denote by \( \pr_{\tau} \) the projection of \( G_{\cf} \) to the \( \tau \)-th factor \( G_{\tau}\cong\GL_{2,\cf} \). Write \( \mu_{\tau}:=\pr_{\tau}\circ \mu_{\tau} \). 
We have \[ \Fl_{G,\mu,\cf}:=B\backslash G_{\cf}\cong \prod_{\tau\in\Sigma}\Fl_{G_{\tau},\mu_{\tau},\cf}, \Fl_{G,\mu,\cf}^{\std}:=\bar{B}\backslash G_{\cf} \cong \prod_{\tau\in\Sigma}\Fl_{G_{\tau},\mu_{\tau},\cf}^{\std},\] and we again denote the projection to the \( \tau \)-th factor as \( \pr_{\tau} \).

Given a multiweight \( \kw \) and \( \tau\in\Sigma \), \( (k_{\tau},w) \) is a weight in the sense of \cref{notation:weight-rep} for \( G_{\tau} \), 
and we have defined \( V^{(k_{\tau},w)} \), \( \lambda^{(k_{\tau},w)} \), \( \chi^{(k_{\tau},w)} \), \( \omega^{(k_{\tau},w)}_{\Fl^{\std}_{G_{\tau},\mu_{\tau},\cf}} \), and \( \omega^{(k_{\tau},w)}_{\Fl_{G_{\tau},\mu_{\tau},\cf}} \)
in \cref{notation:weight-rep}.
We then define \[
V^{\kw}:=\bigotimes_{\tau\in\Sigma}\pr_{\tau}^{*}(V^{(k_{\tau},w)})\in\Rep(G^{c}),\;
\chi^{\kw}:=\bigotimes_{\tau\in\Sigma}\pr_{\tau}^{*}(\chi^{(k_{\tau},w)}),
\] and define \( \lambda^{\kw}:Z(\fg^{c})\to \cf \) to be the infinitesimal character of \( V^{\kw} \). 

We define \[
\omega^{\kw}_{\Fl_{G,\mu,\cf}^{\std}}:=
\bigotimes_{\tau\in\Sigma}\pr_{\tau}^{*}(\omega^{(k_{\tau},w)}_{\Fl^{\std}_{G_{\tau},\mu_{\tau},\cf}}),\;
\omega^{\kw}_{\Fl_{G,\mu,\cf}}:=
\bigotimes_{\tau\in\Sigma}\pr_{\tau}^{*}(\omega^{(k_{\tau},w)}_{\Fl_{G_{\tau},\mu_{\tau},\cf}}).
\]
\end{dfn}

\begin{notation}[Weyl group actions]\label{notation:Weyl-action}
Let \(W\) denote the absolute Weyl group of \(G\), and let \(s_{\tau}\) denote the reflection induced by the simple positive root \( \alpha_{\tau} \) of \(G_{\tau}\). For any \( T\subset \Sigma \),
we denote \(s_{T}:=\prod_{\tau\in T}s_{\tau}\in W\cong (\ZZ/2)^{\Sigma_{\infty}}\),
and denote \(s_{\tau}:\Bbbk\mapsto s_{\tau}\Bbbk\) (resp. \(\Bbbk\mapsto s_{\tau}\cdot\Bbbk\)) acting on \(\ZZ^{\oplus \Sigma_{\infty}}\)
by only changing at \(\tau\)-component via \(k_{\tau}\mapsto -k_{\tau}\) (resp. \(k_{\tau}\mapsto 2-k_{\tau}\)). 
\end{notation}
We have the following corollary of \cref{thm:Vin-GL2-Gr-t} and \cref{thm:Vin-GL2-Fontaine=Diff}:
\begin{cor}\label{cor:plectic-Vin-Fontaine}
Consider the diagram (\ref{eq:Res-Vin-G-basic}). Let \( \theta_{\tau}:=\theta_{T}(d\mu_{\tau}) \) for \( \tau\in\Sigma \). Let \( \kw \) be a regular weight with \( k_{\tau}\ge 2 \) for all \( \tau\in\Sigma \). 
We will make free use of \cref{notation:twist-T-action} and \cref{notation:tuple-eg}. 

(1)
Consider the pro-system 
 \[
\ti{\mF}_{\kw}:=``\varprojlim_{\unl{m}=(m_{\tau})\in\ZZ^{\Sigma}}^{}"
(\ti{i}^{\dR})^{*}\ti{h}_{\natural}(\mO_{\Vin_{G,\cf}^{*}/\ti{G}})\otimes_{Z(\fg)}\lambda^{\kw}\left\{\chi^{(2\unl{m},0)}\right\}
\] 
where
the transition map 
from \( \unl{m} \) to \( \unl{m}-\delta_{\tau} \) 
is induced by right multiplication by \( t_{\tau} \) (with \( t_{\tau}:= \) the coordinate of the \( \tau \)-component of \( (\mA^{1})^{\Sigma} \)) as in \cref{lem:equiv-pullback-D-mod-close-immer} (1), and the action of \( Z(\fg) \) is induced by \( (\fg,*_{l}) \) via \cref{constr:infi-char-on-eqiv-D-mod}.

Then \( \ti{\mF}_{\kw} \)
has an ascending \( t_{\tau}^{\bullet} \)-torsion \( \ZZ^{\Sigma} \)-filtration \( \{\Fil^{(t_{\tau})_{\Sigma}}_{\unl{j}}\}_{\unl{j}\in\ZZ^{\Sigma}} \), with respect to which \( \ti{\mF}_{\kw} \) is complete, and the graded pieces are  \[
\gr^{(t_{\tau})_{\Sigma}}_{\unl{j}}\ti{\mF}_{\kw}\cong h_{0,\natural}(\mO_{\Vin_{G,\cf}^{*}|_{0}/\ti{G}})\otimes_{Z(\fg)}\lambda^{\kw}\left\{\chi^{(-2\unl{j},0)}\right\},\;\unl{j}\in\ZZ^{\Sigma},
\] and
\( (\mfk{t},\theta_{T}) \) acts on \( \gr^{(t_{\tau})_{\Sigma}}_{*}(\ti{\mF}_{\kw}) \) semisimply with integral eigen-characters.

(2) 
The natural morphism 
\[ E_{0}^{(\theta_{\tau})}\left(\ti{\mF}_{\kw}\right)\to E^{0}_{(\theta_{\tau})}\left(\ti{\mF}_{\kw}\right) \]
is an isomorphism, 
where \( E^{(\theta_{\tau})}_{0}:=E^{(\theta_{\tau})_{\tau\in\Sigma}}_{0} \) and \( E_{(\theta_{\tau})}^{0}:=E_{(\theta_{\tau})_{\tau\in\Sigma}}^{0} \) are as in \cref{dfn:E0-nilpotent} taken in the filtered sense with respect to \( \Fil^{(t_{\tau})_{\Sigma}}_{*} \).

Moreover, the pro-system \( E_{0}^{(\theta_{\tau})}\left(\ti{\mF}_{\kw}\right) \)
is eventually constant for sufficiently small \( \unl{m} \).
We denote this limit as \( \mF_{\kw} \). Then
\( \mF_{\kw} \) has an ascending \( t_{\tau}^{\bullet} \)-torsion \( \ZZ^{\Sigma} \)-filtration \( \{\Fil^{(t_{\tau})_{\Sigma}}_{\unl{j}}\}_{\unl{j}\in\ZZ^{\Sigma}} \), with the \emph{non-trivial} graded pieces being 
\begin{equation}\label{eq:filtr-E0-lambda-H}
\gr^{(t_{\tau})_{\Sigma}}_{\unl{j}_{T}}\cong 
	(\pi^{\loc,\dR}_{\GM})^{*} (h_{\Fl}^{\std})_{\natural}(\omega^{-(s_{T}\cdot\Bbbk,w)}_{\Fl^{\std}_{G,\mu,\cf}})\otimes (\pi^{\loc,\dR}_{\HT})^{*} (h_{\Fl})_{\natural}(\omega^{-(s_{T}\cdot\Bbbk,-w)}_{\Fl_{G,\mu,\cf}})
	\left\{\chi^{(w^{\Sigma},w)}\right\},
\end{equation}
for \( T\subset\Sigma \), and \( \unl{j}_{T}:=(-\chi^{(s_{T}\cdot\Bbbk,w)}(d\mu_{\tau}))_{\tau} \).

(3) Denote by \( \Fil^{t_{\tau_{0}}}_{j}\mF_{\kw}:=\varinjlim_{\unl{j}:j_{\tau}=j}\Fil^{(t_{\tau})_{\Sigma}}_{\unl{j}}\mF_{\kw} \) for \( j\in\ZZ \), and denote by \( \gr^{t_{\tau_{0}}}_{j} \) the associated graded piceces, which still have induced \( \ZZ^{\Sigma\backslash \{\tau_{0}\}} \)-filtrations \( \Fil^{(t_{\tau})_{\Sigma\backslash\{\tau_{0}\}}}_{*} \), with \[
\gr^{(t_{\tau})_{\Sigma\backslash\{\tau_{0}\}}}_{(j_{\tau})_{\tau\in\Sigma\backslash\{\tau_{0}\}}}\gr^{t_{\tau_{0}}}_{j_{\tau_{0}}}\mF_{\kw}
\cong \gr^{(t_{\tau})_{\Sigma}}_{(j_{\tau})_{\tau\in\Sigma}}\mF_{\kw}.
\]

For \( \tau_{0}\in\Sigma \), the (right) action of \( \theta_{\tau_{0}} \) preserves the multi-filtration \( \Fil^{(t_{\tau})_{\Sigma}} \), is zero on \( \gr_{*}^{t_{\tau_{0}}}(\mF_{\kw} ) \), and its action on \( \mF_{\kw} \)
induces a monodromy operator \[
\ti{N}_{\tau_{0},k_{\tau_{0}}}:\gr^{t_{\tau_{0}}}_{\frac{k_{\tau_{0}}-w}{2}}(\mF_{\kw})\to \gr^{t_{\tau_{0}}}_{\frac{2-k_{\tau_{0}}-w}{2}}(\mF_{\kw}),
\] preserving \( \Fil^{(t_{\tau})_{\Sigma\backslash\{\tau_{0}\}}}_{*} \) and thus inducing 
\[
\ti{N}_{\tau_{0},k_{\tau_{0}}}:
\gr^{(t_{\tau})_{\Sigma}}_{\unl{j}_{T}}(\mF_{\kw})
\to 
\gr^{(t_{\tau})_{\Sigma}}_{\unl{j}_{T\cup\{\tau_{0}\}}}(\mF_{\kw}),
\] for \( \tau_{0}\notin T\subset \Sigma \).

(4) Along the isomorphism (\ref{eq:filtr-E0-lambda-H}), 
after twisted by \( \left\{\chi^{-(w^{\Sigma},w)}\right\} \), \( \ti{N}_{\tau_{0},k_{\tau_{0}}} \) naturally descends along \( \pi^{\loc,\dR}_{\GM}\times \pi^{\loc,\dR}_{\HT}:(\Vin_{G,\cf}^{*}|_{0})^{\pdR}/\ti{G}\to (\Fl_{G,\mu,\cf}^{\std})^{\pdR}/G \times (\Fl_{G,\mu,\cf})^{\pdR}/G \)
to 
\begin{align*}
\ti{N}_{\tau_{0},k_{\tau_{0}}}:&(h^{\std}_{\Fl})_{\natural}(\omega_{\Fl^{\std}_{G,\mu,\cf}}^{-(s_{T}\cdot\Bbbk,w)})\boxtimes
(h_{\Fl})_{\natural}(\omega^{-(s_{T}\cdot\Bbbk,-w)}_{\Fl_{G,\mu,\cf}})\\&
\to (h^{\std}_{\Fl})_{\natural}(\omega_{\Fl^{\std}_{G,\mu,\cf}}^{-(s_{T\cup\{\tau_{0}\}}\cdot\Bbbk,w)})\boxtimes
(h_{\Fl})_{\natural}(\omega^{-(s_{T\cup\{\tau_{0}\}}\cdot\Bbbk,-w)}_{\Fl_{G,\mu,\cf}}),
\end{align*}
and \[
\ti{N}_{\tau_{0},k_{\tau_{0}}}=c_{k_{\tau}}\cdot \theta^{k_{\tau_{0}}-1}_{\Fl^{\std}_{G,\mu},\tau_{0}}\boxtimes \theta^{k_{\tau_{0}}-1}_{\Fl_{G,\mu},\tau_{0}}
\] where 
\( c_{k_{\tau}}\in\QQ^{\times} \) is as in \cref{thm:Vin-GL2-Fontaine=Diff}, and
\begin{align*}
\theta^{k_{\tau_{0}}-1}_{\Fl^{\std}_{G,\mu},\tau_{0}}&:(h^{\std}_{\Fl})_{\natural}(\omega_{\Fl^{\std}_{G,\mu,\cf}}^{-(s_{T}\cdot\Bbbk,w)})
\to (h^{\std}_{\Fl})_{\natural}(\omega_{\Fl^{\std}_{G,\mu,\cf}}^{-(s_{T\cup\{\tau_{0}\}}\cdot\Bbbk,w)}),\\
\theta^{k_{\tau_{0}}-1}_{\Fl_{G,\mu},\tau_{0}}&:(h_{\Fl})_{\natural}(\omega_{\Fl_{G,\mu,\cf}}^{-(s_{T}\cdot\Bbbk,-w)})
\to (h_{\Fl})_{\natural}(\omega_{\Fl_{G,\mu,\cf}}^{-(s_{T\cup\{\tau_{0}\}}\cdot\Bbbk,-w)}),
\end{align*}
are
the \( \tau_{0} \)-theta operators on \( \Fl^{\std}_{G,\mu} \)
and \( \Fl_{G,\mu} \) respectively, that is, the differentials in the BGG complexes as defined in \cite[\cref{1-propFilterBGG}]{Jiang2026Shla}.
\end{cor}
\begin{proof}
This follows immediately from \cref{thm:Vin-GL2-Gr-t}, \cref{thm:Vin-GL2-Fontaine=Diff} and smooth base change. Note that the \( \tau_{0} \)-theta operators are induced by the theta operators on \( \tau_{0} \)-components by \cite[\cref{1-propPartialBGGDual}]{Jiang2026Shla}.
\end{proof} 
\begin{prop}\label{prop:BdR-to-Plectic}
Consider the diagram (\ref{eq:Pull-ti-G0}).
Let \( \theta_{\tau}:=\theta_{T}(d\mu_{\tau}) \), and \( \theta_{\mu}:=\theta_{T}(d\mu)=\sum_{\tau\in\Sigma}\theta_{\tau} \), and let \( \kw \) be a regular weight with \( k_{\tau}\ge 2 \) for all \( \tau\in\Sigma \). 
Consider the pro-system \[
\ti{\mF}_{\kw}^{0}:=
``\varprojlim_{m\in\ZZ}^{}"
(i^{\dR})^{*}h_{\natural}(\mO_{\Vin_{G,\mu,\cf}^{*}/\ti{G}_{0}})\otimes_{Z(\fg)}\lambda^{\kw}\left\{-m\right\},
\] where the transition maps
are induced by right multiplication by \( t \) (with \( t:= \) the coordinate of \( \mA^{1} \)) as in \cref{lem:equiv-pullback-D-mod-close-immer} (1), and the action of \( Z(\fg) \) is induced by \( (\fg,*_{l}) \) via \cref{constr:infi-char-on-eqiv-D-mod}.

(1) \( \ti{\mF}^{0}_{\kw} \) has an ascending \( t^{\bullet} \)-torsion filtraiton \( \{\Fil^{t}_{j}\}_{j\in\ZZ} \)
with respect to which \( \ti{\mF}^{0}_{\kw} \) is complete, and the graded pieces are 
\[
\gr^{t}_{j}\ti{\mF}^{0}_{\kw}\cong h_{0,\natural}(\mO_{\Vin_{G,\cf}^{*}|_{0}/\ti{G}_{0}})\otimes_{Z(\fg)}\lambda^{\kw}\left\{j\right\},\;j\in\ZZ.
\] 
Moreover, the natural morphism 
\[ E_{0}^{\theta_{\mu}}\left(\ti{\mF}_{\kw}^{0}\right)\to E^{0}_{\theta_{\mu}}\left(\ti{\mF}_{\kw}^{0}\right) \]
is an isomorphism, 
where \( E^{\theta_{\mu}}_{0} \) and \( E_{\theta_{\mu}}^{0} \) are as in \cref{dfn:E0-nilpotent} taken in the filtered sense with respect to \( \Fil^{t}_{*} \).

(2)
Denote by \( \mfk{s} \) the natural map \( (\Vin_{G}^{*}|_{0})^{\pdR}/\ti{G}_{0}\to (\Vin_{G}^{*}|_{0})^{\pdR}/\ti{G} \).
Then we have a natural filtered map \[ \delta:\ti{\mF}_{\kw}^{0}\to \mfk{s}^{*}\ti{\mF}_{\kw}, \]
where \( \ti{\mF}_{\kw} \) is as in \cref{cor:plectic-Vin-Fontaine},
and we endow it with the total filtration \( \Fil^{(t_{\tau})_{\Sigma},tot}_{*} \) of \( \Fil^{(t_{\tau})_{\Sigma}}_{*} \).

(3)
The composition \[
E^{0}_{\theta_{\mu}}(\ti{\mF}_{\kw}^{0})
\xrightarrow{E^{0}_{\theta_{\mu}}(\delta)} \mfk{s}^{*}E^{0}_{\theta_{\mu}}(\ti{\mF}_{\kw})
\to \mfk{s}^{*}E^{0}_{(\theta_{\tau})}(\ti{\mF}_{\kw})
\] is an isomorphism,
where the functors \( E^{0}_{\theta_{\mu}} \)
and \( E^{0}_{(\theta_{\tau})} \)
are taken in the filtered sense.
\end{prop}
\begin{rmk}
The isomorphism \( E^{0}_{\theta_{\mu}}(\ti{\mF}^{0}_{\kw})\cong \mfk{s}^{*}E^{0}_{(\theta_{\tau})}(\ti{\mF}_{\kw}) \)
allows us to translate the structure of the partial filtrations and ``partial Fontaine operators'' on  \( E^{0}_{(\theta_{\tau})}(\ti{\mF}_{\kw}) \) in \cref{cor:plectic-Vin-Fontaine} to \( E^{0}_{\theta_{\mu}}(\ti{\mF}^{0}_{\kw}) \). 
\end{rmk}
\begin{proof}
(1) follows immediately from \cref{lem:equiv-pullback-D-mod-close-immer}. (2) is provided by \cref{lem:pullback-compare-diagonal}.
For (3), as both sides are complete with respect to their filtrations, we can pass to the graded pieces. Then by \cref{lem:pullback-compare-diagonal} (3), \( \gr_{j}(\delta) \) can be identified with \[
\Delta_{\unl{j}_{0}}=(\unl{t}^{\unl{j}_{0}-\unl{j}}):\gr^{(t_{\tau})_{\Sigma}}_{\unl{j}_{0}}(\ti{\mF}_{\kw})\hookrightarrow \gr^{(t_{\tau})_{\Sigma},tot}_{j}(\ti{\mF}_{\kw})\cong \bigoplus_{|\unl{j}|=j}\gr^{(t_{\tau})_{\Sigma}}_{\unl{j}}(\ti{\mF}_{\kw})
\] for any \( |\unl{j}_{0}|=j \). 
By \cref{cor:plectic-Vin-Fontaine} (1), \( \theta_{\tau} \)
acts on both sides semisimply with integral eigen-characters.
Moreover, multiplication by \( \unl{t}^{\unl{j}_{0}-\unl{j}} \) changes the action of \( \theta_{\tau} \) by \( j_{\tau}-j_{0,\tau} \). It follows easily that the composition \[
E^{0}_{\theta_{\mu}}(\gr^{(t_{\tau})_{\Sigma}}_{\unl{j}_{0}}(\ti{\mF}_{\kw}))\hookrightarrow{\Delta_{\unl{j}_{0}}}
E^{0}_{\theta_{\mu}}(\gr^{(t_{\tau})_{\Sigma},tot}_{j}(\ti{\mF}_{\kw}))
\twoheadrightarrow 
E^{0}_{(\theta_{\tau})}(\gr^{(t_{\tau})_{\Sigma},tot}_{j}(\ti{\mF}_{\kw}))
\] is an isomorphism.
\end{proof}

Now we pass from \( G \) to \( G^{c} \).
Let \( \ti{G}^{c}:=G^{c}\times G^{c}\times T^{c}\supset \ti{G}_{0}^{c}=G^{c}\times G^{c}\times \GG_{m} \).
We consider the following commutative diagram \begin{equation}\label{eq:Pull-ti-G0-Gc}
\begin{tikzcd}
\Vin_{G^{c},\cf}^{*}|_{0}/\ti{G}^{c}_{0}\arrow[r,"i"]\arrow[rr,"\ti{i}",bend left=12pt]\arrow[d,"h_{0}"] &
\Vin_{G^{c},\mu,\cf}^{*}/\ti{G}^{c}_{0}\arrow[r,"\Delta"]\arrow[d,"h"]
& \Vin_{G^{c},\cf}^{*}/\ti{G}^{c}_{0}
\arrow[d,"\ti{h}"]\\
(\Vin_{G^{c},\cf}^{*}|_{0})^{\pdR}/\ti{G}^{c}_{0}\arrow[r,"i^{\dR}"]\arrow[rr,"\ti{i}^{\dR}",bend right=12pt] 
\arrow[d,"\pi_{\GM}^{\loc,\dR}\times \pi_{\HT}^{\loc,\dR}"]
&
(\Vin_{G^{c},\mu,\cf}^{*})^{\pdR}/\ti{G}^{c}_{0}\arrow[r,"\Delta^{\dR}"]
& (\Vin_{G^{c},\cf}^{*})^{\pdR}/\ti{G}^{c}_{0},\\
(\Fl_{G,\mu}^{\std})^{\pdR}/G^{c}
\times \Fl_{G,\mu}^{\pdR}/G^{c} & 
\Fl_{G,\mu}^{\std}/G^{c}
\times \Fl_{G,\mu}/G^{c}\arrow[l,"h_{\Fl}^{\std}\times h_{\Fl}"]
\end{tikzcd}
\end{equation}
where the larger square is pulled back from \begin{equation}\label{eq:Res-Vin-G-basic-Gc}
\begin{tikzcd}
\Vin_{G^{c},\cf}^{*}|_{0}/\ti{G}^{c}\arrow[r,"\ti{i}"]\arrow[d,"h_{0}"]
& \Vin_{G^{c},\cf}^{*}/\ti{G}^{c}
\arrow[d,"\ti{h}"]\\
(\Vin_{G^{c},\cf}^{*}|_{0})^{\pdR}/\ti{G}^{c}\arrow[r,"\ti{i}^{\dR}"]
& (\Vin_{G^{c},\cf}^{*})^{\pdR}/\ti{G}^{c},
\end{tikzcd}
\end{equation} along \( \mfk{s}:B\ti{G}^{c}_{0}\to B\ti{G}^{c} \).
\begin{cor}\label{cor:from-G-to-Gc}
The analogues of \cref{cor:plectic-Vin-Fontaine}
and \cref{prop:BdR-to-Plectic} hold if we replace \( G \) with \( G^{c} \), 
the diagram (\ref{eq:Pull-ti-G0}) with (\ref{eq:Pull-ti-G0-Gc}), and
the diagram (\ref{eq:Res-Vin-G-basic}) with (\ref{eq:Res-Vin-G-basic-Gc}). 

More precisely, consider the diagram (\ref{eq:Pull-ti-G0-Gc}) and (\ref{eq:Res-Vin-G-basic-Gc}). Let \( \theta_{\tau}:=\theta_{T^{c}}(d\mu_{\tau}) \) for \( \tau\in\Sigma \), and \( \theta:=\theta_{T^{c}}(d\mu)=\sum_{\tau\in\Sigma}\theta_{\tau} \). Let \( \kw \) be a regular weight with \( k_{\tau}\ge 2 \) for all \( \tau\in\Sigma \). 

(1) Consider the pro-system \[
\ti{\mF}_{\kw}^{c}:=``\varprojlim_{\unl{m}=(m_{\tau})\in\ZZ^{\Sigma}}^{}"
(\ti{i}^{\dR})^{*}\ti{h}_{\natural}(\mO_{\Vin_{G^{c},\cf}^{*}/\ti{G}^{c}})\otimes_{Z(\fg^{c})}\lambda^{\kw}\left\{\chi^{(2\unl{m},0)}\right\}
\] 
where
the transition map 
from \( \unl{m} \) to \( \unl{m}-\delta_{\tau} \) 
is induced by right multiplication by \( t_{\tau} \) as in \cref{lem:equiv-pullback-D-mod-close-immer} (1), and the action of \( Z(\fg^{c}) \) is induced by \( (\fg^{c},*_{l}) \) via \cref{constr:infi-char-on-eqiv-D-mod}.

Then \( \ti{\mF}_{\kw}^{c} \)
has an ascending complete \( t_{\tau}^{\bullet} \)-torsion \( \ZZ^{\Sigma} \)-filtration \( \{\Fil^{(t_{\tau})_{\Sigma}}_{\unl{j}}\}_{\unl{j}\in\ZZ^{\Sigma}} \), 
with  \[
\gr^{(t_{\tau})_{\Sigma}}_{\unl{j}}\ti{\mF}_{\kw}^{c}\cong h_{0,\natural}(\mO_{\Vin_{G^{c},\cf}^{*}|_{0}/\ti{G}^{c}})\otimes_{Z(\fg^{c})}\lambda^{\kw}\left\{\chi^{(-2\unl{j},0)}\right\},\;\unl{j}\in\ZZ^{\Sigma},
\] and
\( (\mfk{t}^{c},\theta_{T^{c}}) \) acts on \( \gr^{(t_{\tau})_{\Sigma}}_{*}(\ti{\mF}_{\kw}^{c}) \) semisimply with integral eigen-characters.

(2) The pro-system \( E_{0}^{(\theta_{\tau})}\left(\ti{\mF}_{\kw}^{c}\right)\cong E^{0}_{(\theta_{\tau})}\left(\ti{\mF}_{\kw}^{c}\right) \)
is eventually constant, where \( E^{0}_{(\theta_{\tau})} \) and \( E_{0}^{(\theta_{\tau})} \) are taken in the filtered sense.
We denote this limit as \( \mF_{\kw}^{c} \). Then
\( \mF_{\kw}^{c} \) has an ascending \( t_{\tau}^{\bullet} \)-torsion \( \ZZ^{\Sigma} \)-filtration \( \{\Fil^{(t_{\tau})_{\Sigma}}_{\unl{j}}\}_{\unl{j}\in\ZZ^{\Sigma}} \), with the \emph{non-trivial} graded pieces being 
\begin{equation}\label{eq:filtr-E0-lambda-H-new}
\gr^{(t_{\tau})_{\Sigma}}_{\unl{j}_{T}}(\mF^{c}_{\kw})\cong 
	(\pi^{\loc,\dR}_{\GM})^{*} (h_{\Fl}^{\std})_{\natural}(\omega^{-(s_{T}\cdot\Bbbk,w)}_{\Fl^{\std}_{G,\mu,\cf}})\otimes (\pi^{\loc,\dR}_{\HT})^{*} (h_{\Fl})_{\natural}(\omega^{-(s_{T}\cdot\Bbbk,-w)}_{\Fl_{G,\mu,\cf}})
	\left\{\chi^{(w^{\Sigma},w)}\right\},
\end{equation}
for \( T\subset\Sigma \), and \( \unl{j}_{T}:=(-\chi^{(s_{T}\cdot\Bbbk,w)}(d\mu_{\tau}))_{\tau} \).

(3)
For \( \tau_{0}\in\Sigma \), the (right) action of \( \theta_{\tau_{0}} \) preserves the multi-filtration \( \Fil^{(t_{\tau})_{\Sigma}} \), is zero on \( \gr_{*}^{t_{\tau_{0}}}(\mF_{\kw}^{c} ) \), and its action on \( \mF_{\kw}^{c} \)
induces a monodromy operator \[
\ti{N}_{\tau_{0},k_{\tau_{0}}}:\gr^{t_{\tau_{0}}}_{\frac{k_{\tau_{0}}-w}{2}}(\mF_{\kw}^{c})\to \gr^{t_{\tau_{0}}}_{\frac{2-k_{\tau_{0}}-w}{2}}(\mF_{\kw}^{c}),
\] preserving \( \Fil^{(t_{\tau})_{\Sigma\backslash\{\tau_{0}\}}}_{*} \) and thus inducing 
\[
\ti{N}_{\tau_{0},k_{\tau_{0}}}:
\gr^{(t_{\tau})_{\Sigma}}_{\unl{j}_{T}}(\mF_{\kw}^{c})
\to 
\gr^{(t_{\tau})_{\Sigma}}_{\unl{j}_{T\cup\{\tau_{0}\}}}(\mF_{\kw}^{c}),
\] for \( \tau_{0}\notin T\subset \Sigma \).

(4)
Along the isomorphism (\ref{eq:filtr-E0-lambda-H-new}), 
after twisted by \( \left\{\chi^{-(w^{\Sigma},w)}\right\} \), \( \ti{N}_{\tau_{0},k_{\tau_{0}}} \) naturally descends along \( \pi^{\loc.\dR}_{\GM}\times \pi^{\loc.\dR}_{\HT} \)
to 
\begin{align*}
\ti{N}_{\tau_{0},k_{\tau_{0}}}:&(h^{\std}_{\Fl})_{\natural}(\omega_{\Fl^{\std}_{G,\mu,\cf}}^{-(s_{T}\cdot\Bbbk,w)})\boxtimes
(h_{\Fl})_{\natural}(\omega^{-(s_{T}\cdot\Bbbk,-w)}_{\Fl_{G,\mu,\cf}})\\&
\to (h^{\std}_{\Fl})_{\natural}(\omega_{\Fl^{\std}_{G,\mu,\cf}}^{-(s_{T\cup\{\tau_{0}\}}\cdot\Bbbk,w)})\boxtimes
(h_{\Fl})_{\natural}(\omega^{-(s_{T\cup\{\tau_{0}\}}\cdot\Bbbk,-w)}_{\Fl_{G,\mu,\cf}}),
\end{align*}
and \[
\ti{N}_{\tau_{0},k_{\tau_{0}}}=c_{k_{\tau}}\cdot \theta^{k_{\tau_{0}}-1}_{\Fl^{\std}_{G,\mu},\tau_{0}}\boxtimes \theta^{k_{\tau_{0}}-1}_{\Fl_{G,\mu},\tau_{0}}
\] where the RHS is defined as in \cref{cor:plectic-Vin-Fontaine} (4). 

(5) Consider the pro-system \[
\ti{\mF}_{\kw}^{c,0}:=
``\varprojlim_{m\in\ZZ}^{}"
(i^{\dR})^{*}h_{\natural}(\mO_{\Vin_{G^{c},\mu,\cf}^{*}/\ti{G}^{c}_{0}})\otimes_{Z(\fg^{c})}\lambda^{\kw}\left\{-m\right\},
\] where the transition maps
are induced by right multiplication by \( t \) as in \cref{lem:equiv-pullback-D-mod-close-immer} (1), which is equipped with an ascending \( t^{\bullet} \)-torsion filtraiton \( \{\Fil^{t}_{j}\}_{j\in\ZZ} \) and an action of \( \theta \).

The pro-system \( E_{0}^{\theta}\left(\ti{\mF}_{\kw}^{c,0}\right)\cong E^{0}_{\theta}\left(\ti{\mF}_{\kw}^{c,0}\right) \)
is eventually constant, where \( E^{0}_{\theta} \) is taken in the filtered sense, and there exists a natural filtered isomorphism \[
E^{0}_{\theta}\left(\ti{\mF}_{\kw}^{c,0}\right) \cong E^{0}_{(\theta_{\tau})}\left(\ti{\mF}_{\kw}^{c}\right)
\] that is compatible with \( \theta \)-actions, where the LHS is endowed with the induced \( t^{\bullet} \)-torsion filtration \( \Fil^{t}_{*} \)
and the RHS is endowed with the total filtration of the \( \ZZ^{\Sigma} \)-filtration \( \Fil^{(t_{\tau})_{\Sigma}}_{*} \).
\end{cor}
\begin{proof}
Consider the diagram \[
\begin{tikzcd}
\Vin^{*}_{G,\cf}/\ti{G}\arrow[r,"p"]\arrow[d,"\ti{h}"]&
\Vin^{*}_{G^{c},\cf}/\ti{G}^{c}\arrow[r,"\mbf{t}"]\arrow[d,"\ti{h}"] & (\mA^{1})^{\Sigma}\\
(\Vin_{G,\cf}^{*})^{\pdR}/\ti{G}\arrow[r,"p"]&
(\Vin_{G^{c},\cf}^{*})^{\pdR}/\ti{G}^{c}\\
(\Vin_{G,\cf}^{*}|_{0})^{\pdR}/\ti{G}\arrow[r,"p"]\arrow[u,"i^{\dR}"]\arrow[r,"p"]\arrow[d,"\pi_{\GM}^{\loc,\dR}\times \pi_{\HT}^{\loc,\dR}"]&
(\Vin_{G^{c},\cf}^{*}|_{0})^{\pdR}/\ti{G}^{c}\arrow[u,"i^{\dR}"]
\arrow[d,"\pi_{\GM}^{\loc,\dR}\times \pi_{\HT}^{\loc,\dR}"]\\
(\Fl_{G,\mu}^{\std})^{\pdR}/G
\times \Fl_{G,\mu}^{\pdR}/G\arrow[r,"p"] & 
(\Fl_{G,\mu}^{\std})^{\pdR}/G^{c}
\times \Fl_{G,\mu}^{\pdR}/G^{c}
\end{tikzcd}
\]
Since \( \lambda^{\kw}|_{\Lie(Z_{c})}=0 \), it follows from \cref{lem:fix-t-action-on-torsor} that 
\[
\ti{h}_{\natural}(\mO_{\Vin^{*}_{G,\cf}/\ti{G}})\otimes_{Z(\fg)}\lambda^{\kw}
\cong p^{*}\left(h_{0,\natural}(\mO_{\Vin_{G^{c},\cf}^{*}|_{0}/\ti{G}^{c}})\otimes_{Z(\fg^{c})}\lambda^{\kw}\right).
\]
Thus by the construction of the \( \ZZ^{\Sigma} \)-filtration in \cref{lem:equiv-pullback-D-mod-close-immer}, we have a \( \ZZ^{\Sigma} \)-filtered isomorphism \[
p^{*}(\ti{\mF}^{c}_{\kw})\cong \ti{\mF}_{\kw}.
\]
A similar argument shows that we have a filtered isomorphism \[
p^{*}(\ti{\mF}^{c,0}_{\kw})\cong \ti{\mF}_{\kw}^{0}.
\]
The description of their graded pieces follows immediately from \cref{cor:plectic-Vin-Fontaine}
and \cref{prop:BdR-to-Plectic} via \cref{lem:Bt-gerbe-heart-fufaith}
and \cref{lem:BT-gerbe-fufaith}. More precisely, we decompose the natural map \( p \) into a composition \[
p:(\Vin^{*}_{G,\cf}|_{0})^{\pdR}/\ti{G}\xrightarrow{p_{1}} (\Vin^{*}_{G^{c},\cf}|_{0})^{\pdR}/(G\times G\times T^{c})\xrightarrow{\mfk{s}_{2}} (\Vin^{*}_{G^{c},\cf}|_{0})^{\pdR}/\ti{G}^{c}, 
\] and we apply \cref{lem:Bt-gerbe-heart-fufaith} to \( p_{1}^{*} \)
and \cref{lem:BT-gerbe-fufaith} to \( p_{2}^{*} \).

In particular, we see that \( (\mfk{t},\theta_{T^{c}}) \) acts on the graded pieces semisimply with integral eigen-characters. This implies that 
\[
E_{0}^{(\theta_{\tau})}\left(\ti{\mF}_{\kw}^{c}\right)\cong E^{0}_{(\theta_{\tau})}\left(\ti{\mF}_{\kw}^{c}\right),\;
E_{0}^{\theta}\left(\ti{\mF}_{\kw}^{c,0}\right)\cong E^{0}_{\theta}\left(\ti{\mF}_{\kw}^{c,0}\right).
\] Note that \( E^{0}_{(\theta_{\tau})} \) and \( E^{0}_{\theta} \)
commutes with \( p^{*} \). Thus we have \[
p^{*}(E^{0}_{(\theta_{\tau})}\left(\ti{\mF}_{\kw}^{c}\right))\cong E^{0}_{(\theta_{\tau})}\left(\ti{\mF}_{\kw}\right),
p^{*}(E^{0}_{\theta}\left(\ti{\mF}_{\kw}^{c,0}\right))
\cong E^{0}_{\theta}\left(\ti{\mF}_{\kw}^{0}\right).
\] Now we can again apply \cref{lem:Bt-gerbe-heart-fufaith}
and \cref{lem:BT-gerbe-fufaith}, and reduce the problem to \cref{cor:plectic-Vin-Fontaine}
and \cref{prop:BdR-to-Plectic}.
\end{proof}

\section{Fontaine operators on Shimura varieties}\label{sec:Fontaine-Ope-HMV}
\subsection{Differential operators}\label{subsec:diff-operator}
In \cite[\S \ref{1-subsec:Diff-ope}]{Jiang2026Shla}, we gave a general construction
of differential operators using filtered complete de Rham stacks, which we will recall in \cref{constr:Diff-op-old}.  
In this subsection, we will give another construction (\cref{constr:Diff-op}) using only the algebraic de Rham stacks, which is better suited for our purpose. 
We will compare the two constructions (\cref{prop:const2-coincide}), so that we could apply our previous results on the classicality of partial de Rham cohomology (\cite[\cref{2-thmClassicalityCohomoDR}]{Jiang2025HMF}).

Consider the diagram (\ref{eq:Alg-DR-Bdif}). Pulling back along \( \Vin_{G^{c},\mu,\CC_{p}}|_{0}\to \Vin_{G^{c},\mu,\cfi}|_{\hat{\mA}^{1}} \), we obtain a Cartesian diagram \[
\begin{tikzcd}
\mX^{\tor,\la}_{K^{p}}
\arrow[r,"\pi^{\la}_{\GM\HT}"]\arrow[d,"h^{\Bdif}_{0}"]
& \Vin_{G^{c},\mu,\CC_{p}}|_{0}/G^{c}\arrow[d,"h_{0}"]\\
\bar{\mX}^{\tor,\la}_{K^{p}}\arrow[r,"\bar{\pi}_{\GM\HT}"] & (\Vin^{*}_{G^{c},\mu,\CC_{p}}|_{0})^{\pdR/\CC_{p}}/G^{c},
\end{tikzcd}
\] where \( \bar{\mX}^{\tor,\la}_{K^{p}}:=\bar{\mX}^{\tor,\la}_{K^{p},\cfi}\hatotimes_{\cfi}\CC_{p} \).

Composing \( \bar{\pi}_{\GM\HT} \)
with \[
\pi_{\GM}^{\loc,\dR}\times \pi_{\HT}^{\loc,\dR}:
(\Vin^{*}_{G^{c},\mu,\CC_{p}}|_{0})^{\pdR/\CC_{p}}/G^{c}\to 
\Fl_{G,\mu,\CC_{p}}^{\std,\pdR/\CC_{p}}/G^{c}\times \Fl_{G,\mu,\CC_{p}}^{\pdR/\CC_{p}},
\]
and pulling back along \[
h_{\Fl}^{\std}\times h_{\Fl}:\Fl_{G,\mu}^{\std}\times \Fl_{G,\mu}
\to \Fl_{G,\mu}^{\std,\pdR/\CC_{p}}\times \Fl_{G,\mu}^{\pdR/\CC_{p}},
\] we obtain another Cartesian diagram
\begin{equation}\label{eq:Cartesian-Ola0-Cp}
\begin{tikzcd}
\mX^{\tor,\la,0}_{K^{p}}
\arrow[r]\arrow[d,"h^{\la,0}"] &[-2.1em] 
\mX^{\tor,\sm}_{K^{p}}\times \Fl_{G,\mu} \arrow[r,"\pi^{\sm}_{\GM}\times 1"]\arrow[d,"h^{\sm}\times h_{\Fl}"]
&[0.8em]\Fl_{G,\mu}^{\std}/G^{c}\times \Fl_{G,\mu}\arrow[d,"h_{\Fl}^{\std}\times h_{\Fl}"]\\
\bar{\mX}^{\tor,\la}_{K^{p}}
\arrow[r]& 
(\mX^{\tor}_{K^{p}})^{\pdR/\CC_{p},\sm}_{\log}\times \Fl_{G,\mu}^{\pdR/\CC_{p}} 
\arrow[r,"\pi^{\sm,\pdR}_{\GM}\times 1"] & 
\Fl_{G,\mu}^{\std,\pdR/\CC_{p}}/G^{c}\times \Fl_{G,\mu}^{\pdR/\CC_{p}}
\end{tikzcd}
\end{equation}
by \cite[\cref{1-corCartesianHilbertFlandSm}]{Jiang2026Shla}.

By construction, the square on the left is actually the restriction of the following square over \( [\mA^{1}/\GG_{m}]^{2} \)  along \( [\GG_{m}/\GG_{m}]^{2}\subset [\mA^{1}/\GG_{m}]^{2} \)
\begin{equation}\label{eq:la-0-cartesian}
\begin{tikzcd}
{\mX^{\tor,\la,0}_{K^{p}}\times [{\mA}^{1}/\GG_{m}]^{2}}
\arrow[r," \pi^{\la,0}_{\sm}\times\pi^{\la,0}_{\HT}"]
\arrow[d,"h^{\la,0,++}"] &
{(\mX^{\tor,\sm}_{K^{p}}\times [\mA^{1}/\GG_{m}])}
\times  
{(\Fl_{G,\mu}\times [\mA^{1}/\GG_{m}])}
\arrow[d," h^{\sm,+}\times h_{\Fl}^{+}"]
\\ 
(\mX^{\tor,\la,0}_{K^{p}})^{\pdR/\CC_{p},++}_{\log}
\arrow[r,"\pi^{\pdR,++}_{\GM\HT}"] 
&  (\mX^{\tor}_{K^{p}})^{\pdR/\CC_{p},{+},\sm}_{\log}\times\Fl_{G,\mu}^{\pdR/\CC_{p},+} ,
\end{tikzcd}
\end{equation}
constructed 
in \cite[\cref{1-CordfnBiFilAgLa0}]{Jiang2026Shla}.

\begin{construction}\label{constr:Diff-op}
Assume that we are given two vector bundles \( \mF_{1},\mF_{2} \) on \( \Fl^{\std}_{G,\mu}/G^{c} \), with a \( G^{c} \)-equivariant differential operator \( \delta \) from \( \mF_{1} \) to \( \mF_{2} \), and 
two vector bundles \( \mF'_{1},\mF_{2}' \) on \( \Fl_{G,\mu} \), with a differential operator \( \delta' \) from \( \mF'_{1} \) to \( \mF'_{2} \). 
Then \( \delta \)
and \( \delta' \)
corresponds to morphisms \[
\delta^{\vee}:(h_{\Fl}^{\std})_{\natural}(\mF_{2}^{\vee})
\to (h_{\Fl}^{\std})_{\natural}(\mF_{1}^{\vee}),\;
\delta^{\prime,\vee}:(h_{\Fl})_{\natural}(\mF_{2}^{\prime,\vee})
\to (h_{\Fl})_{\natural}(\mF_{1}^{\prime,\vee}).
\]
By smooth base change, we have natural isomorphisms 
\begin{align*}
h^{\la,0}_{\natural}
(\pi^{\la,0}_{\GM}\times \pi^{\la,0}_{\HT})^{*}(\mF_{1}^{\vee}\boxtimes \mF_{1}^{\prime,\vee})&\cong 
(\pi^{\dR}_{\GM}\times \pi^{\dR}_{\HT})^{*}((h_{\Fl}^{\std})_{\natural}(\mF_{1}^{\vee})\otimes (h_{\Fl})_{\natural}(\mF_{1}^{\prime,\vee})),\\
h^{\la,0}_{\natural}
(\pi^{\la,0}_{\GM}\times \pi^{\la,0}_{\HT})^{*}(\mF_{2}^{\vee}\boxtimes \mF_{2}^{\prime,\vee})&\cong 
(\pi^{\dR}_{\GM}\times \pi^{\dR}_{\HT})^{*}((h_{\Fl}^{\std})_{\natural}(\mF_{2}^{\vee})\otimes (h_{\Fl})_{\natural}(\mF_{2}^{\prime,\vee})),
\end{align*}
where \( \pi^{\la,0}_{\GM}:=\pi^{\sm}_{\GM}\circ \pi^{\la,0}_{\sm} \), and \( \pi^{\dR}_{\GM}\times \pi^{\dR}_{\HT}:=(\pi^{\loc,\dR}_{\GM}\times \pi^{\loc,\dR}_{\HT})\circ \bar{\pi}_{\GM\HT} \).
Thus \( \delta^{\vee}\boxtimes \delta^{\prime,\vee} \)
induces \[
h^{\la,0}_{\natural}
(\pi^{\la,0}_{\GM}\times \pi^{\la,0}_{\HT})^{*}(\mF_{2}^{\vee}\boxtimes \mF_{2}^{\prime,\vee})
\to
h^{\la,0}_{\natural}
(\pi^{\la,0}_{\GM}\times \pi^{\la,0}_{\HT})^{*}(\mF_{1}^{\vee}\boxtimes \mF_{1}^{\prime,\vee}).
\]
Applying \( (-)^{\vee} \)
on both sides, we obtain a morphism \[\delta\boxtimes\delta':
h^{\la,0}_{*}
(\pi^{\la,0}_{\GM}\times \pi^{\la,0}_{\HT})^{*}(\mF_{1}\boxtimes \mF_{1}')
\to
h^{\la,0}_{*}
(\pi^{\la,0}_{\GM}\times \pi^{\la,0}_{\HT})^{*}(\mF_{2}\boxtimes \mF_{2}').
\]
\end{construction}
\begin{rmk}\label{rmk:dif-op-L-vs-Cp}
We can consider a variant with (\ref{eq:Cartesian-Ola0-Cp}) replaced by its descent from \( \CC_{p} \) to \( \cfi \) provided by (\ref{eq:Alg-DR-Bdif}). Concretely, we define \( \mX^{\tor,\la,0}_{K^{p},\cfi} \) via the following Cartesian diagram \begin{equation}\label{eq:Cartesian-Ola0-Linfty}
\begin{tikzcd}[column sep=50pt]
\mX^{\tor,\la,0}_{K^{p},\cfi}\arrow[d,"h^{\la,0}_{\cfi}"]
\arrow[r,"\pi^{\la,0}_{\GM,\cfi}\times \pi^{\la,0}_{\HT,\cfi}"]
& \Fl^{\std}_{G,\mu,\cfi}/G^{c}\times \Fl_{G,\mu,\cfi}\arrow[d,"h^{\std}_{\Fl,\cfi}\times h_{\Fl,\cfi}"]\\
\bar{\mX}^{\tor,\la}_{K^{p},\cfi}\arrow[r,"\pi^{\dR}_{\GM,\cfi}\times \pi^{\dR}_{\HT,\cfi}"]
& \Fl^{\std,\pdR/\cfi}_{G,\mu,\cfi}/G^{c}\times \Fl^{\pdR/\cfi}_{G,\mu,\cfi},
\end{tikzcd}
\end{equation}
whose pull-back along \(\mfk{s}: \AnSpec(\CC_{p})\to \AnSpec(\cfi) \)
recovers (\ref{eq:Cartesian-Ola0-Cp}).
Let us denote also by \( \mfk{s} \)
the map from the terms in (\ref{eq:Cartesian-Ola0-Cp}) to the corresponding terms in (\ref{eq:Cartesian-Ola0-Linfty}).

Assume that we are given two vector bundles \( \mF_{1},\mF_{2} \) on \( \Fl^{\std}_{G,\mu,\cfi}/G^{c} \), with a \( G^{c} \)-equivariant differential operator \( \delta \) from \( \mF_{1} \) to \( \mF_{2} \), and 
two vector bundles \( \mF'_{1},\mF_{2}' \) on \( \Fl_{G,\mu,\cfi} \), with a differential operator \( \delta' \) from \( \mF'_{1} \) to \( \mF'_{2} \). Then by the similar construction as in \cref{constr:Diff-op}, we obtain a morphism \[
\delta\boxtimes \delta':(h^{\la,0}_{\cfi})_{*}
(\pi^{\la,0}_{\GM,\cfi}\times \pi^{\la,0}_{\HT,\cfi})^{*}(\mF_{1}\boxtimes \mF_{1}')
\to
(h^{\la,0}_{\cfi})_{*}
(\pi^{\la,0}_{\GM,\cfi}\times \pi^{\la,0}_{\HT,\cfi})^{*}(\mF_{2}\boxtimes \mF_{2}').
\]

We can pull-back \( \delta \) and \( \delta' \)
along \( \mfk{s} \), and we obtain via \cref{constr:Diff-op} \[
\mfk{s}^{*}(\delta)\boxtimes\mfk{s}^{*}(\delta'):
h^{\la,0}_{*}
(\pi^{\la,0}_{\GM}\times \pi^{\la,0}_{\HT})^{*}(\mfk{s}^{*}(\mF_{1})\boxtimes \mfk{s}^{*}(\mF_{1}'))
\to
h^{\la,0}_{*}
(\pi^{\la,0}_{\GM}\times \pi^{\la,0}_{\HT})^{*}(\mfk{s}^{*}(\mF_{2})\boxtimes \mfk{s}^{*}(\mF_{2}')).
\]
Then it is straight-forward to see that we have 
\begin{equation}\label{eq:dif-op-Cp-vs-L}
\mfk{s}_{*}\mfk{s}^{*}(\delta\boxtimes \delta')=\mfk{s}_{*}(\mfk{s}^{*}(\delta)\boxtimes\mfk{s}^{*}(\delta')).
\end{equation}
Note that concretely, \( \mfk{s}_{*}\mfk{s}^{*}M\cong M\otimes_{\cfi}\CC_{p} \).
\end{rmk}

Let us also recall the following construction in \cite[\S \ref{1-subsec:Diff-ope}]{Jiang2026Shla}.
\begin{construction}\label{constr:Diff-op-old}
Given \( (\mF_{1},\mF_{2},\delta) \)
and \( (\mF_{1}',\mF_{2}',\delta') \)
as in \cref{constr:Diff-op}. 
As in the proof of \cite[\cref{1-propBGGFilSVsm}]{Jiang2026Shla},
\( (\mF_{1},\mF_{2},\delta) \) induces a differential operator \( \delta \) from \( (\pi^{\sm}_{\GM})^{*}(\mF_{1}) \)
to \( (\pi^{\sm}_{\GM})^{*}(\mF_{2}) \).

Assume that \( \delta \) (resp. \( \delta' \)) is of degree \( n \) (resp. \( n' \)). 
Then by \cite[\cref{1-lemDiffOpeToLogStack} (2)]{Jiang2026Shla}, they correspond to morphisms \begin{align*}
\delta^{+}&:h^{\sm,+}_{*}(\pi^{\sm}_{\GM})^{*}(\mF_{1})\{-n\}\to h^{\sm,+}_{*}(\pi^{\sm}_{\GM})^{*}(\mF_{2}),\\
\delta^{\prime,+}&:(h^{+}_{\Fl})_{*}(\mF_{1}')\{-n'\}\to (h^{+}_{\Fl})_{*}(\mF_{2}'),
\end{align*} where \( \{i\} \)
refers to twisting by \( i \)-th power of the tautological line bundle on \( [\mA^{1}/\GG_{m}] \).

It is established in \cite[\cref{1-lemKunnFlSm}]{Jiang2026Shla}
that the natural morphisms \begin{align*}
(\pi^{\pdR,++}_{\GM\HT})^{*}(h^{\sm,+}_{*}(\pi^{\sm}_{\GM})^{*}(\mF_{1})\boxtimes (h^{+}_{\Fl})_{*}(\mF_{1}'))
\to h^{\la,0,++}_{*}(\pi^{\la,0}_{\GM}\times \pi^{\la,0}_{\HT})^{*}(\mF_{1}\boxtimes \mF_{1}'),\\
(\pi^{\pdR,++}_{\GM\HT})^{*}(h^{\sm,+}_{*}(\pi^{\sm}_{\GM})^{*}(\mF_{2})\boxtimes (h^{+}_{\Fl})_{*}(\mF_{2}'))
\to h^{\la,0,++}_{*}(\pi^{\la,0}_{\GM}\times \pi^{\la,0}_{\HT})^{*}(\mF_{2}\boxtimes \mF_{2}'),
\end{align*} become isomorphisms when restricted to \( [\hat{\mA}^{1}/\GG_{m}]^{2}\subset [\mA^{1}/\GG_{m}]^{2} \).

Denote \( j:[\hat{\mA}^{1}/\GG_{m}]\hookrightarrow [\mA^{1}/\GG_{m}]\hookleftarrow [\GG_{m}/\GG_{m}]:i_{0} \). Since \( (\pi^{\la,0}_{\GM}\times \pi^{\la,0}_{\HT})^{*}(\mF_{i}\boxtimes \mF_{i}') \) (\( i\in\{1,2\} \)) is endowed with the trivial filtration, it is filtered complete. By \cite[\cref{1-corFilterdatA1modGm}]{Jiang2026Shla}, \[
(\pi^{\la,0}_{\GM}\times \pi^{\la,0}_{\HT})^{*}(\mF_{i}\boxtimes \mF_{i}')\cong (j\times j)_{*}(j\times j)^{*}(\pi^{\la,0}_{\GM}\times \pi^{\la,0}_{\HT})^{*}(\mF_{i}\boxtimes \mF_{i}').
\]
Then \( \delta^{+}\boxtimes \delta^{\prime,+} \)
induces \[
d:h^{\la,0,++}_{*}
(\pi^{\la,0}_{\GM}\times \pi^{\la,0}_{\HT})^{*}(\mF_{1}\boxtimes \mF_{1}')\{-n,-n'\}
\to
h^{\la,0,++}_{*}
(\pi^{\la,0}_{\GM}\times \pi^{\la,0}_{\HT})^{*}(\mF_{2}\boxtimes \mF_{2}').
\]
Pulling back along \( i_{0}\times i_{0} \) and applying \cref{lem:push-forward-Xla0} 
give a morphism \[
\delta\boxtimes \delta':
h^{\la,0}_{*}
(\pi^{\la,0}_{\GM}\times \pi^{\la,0}_{\HT})^{*}(\mF_{1}\boxtimes \mF_{1}')
\to
h^{\la,0}_{*}
(\pi^{\la,0}_{\GM}\times \pi^{\la,0}_{\HT})^{*}(\mF_{2}\boxtimes \mF_{2}').
\] 
\end{construction}
\begin{lem}\label{lem:push-forward-Xla0}
For any vector bundle \( V\in\QCoh(\mX^{\la,0}_{K^{p}}) \), the natural morphism \[
(i_{0}\times i_{0})^{*}h^{\la,0,++}_{*}
(V)\to h^{\la,0}_{*}(V)
\] is an isomorphism.
\end{lem}
\begin{proof}
By \cite[\cref{1-lemBasicFiltered} \& \cref{1-CordfnBiFilAgLa0}]{Jiang2026Shla}, \( h^{\la,0,++} \) is a suave \( ! \)-cover. Thus it suffices to show that \[
(h^{\la,0})^{*}(i_{0}\times i_{0})^{*}h^{\la,0,++}_{*}
(V)\cong 
(i_{0}\times i_{0})^{*}(h^{\la,0,++})^{*}h^{\la,0,++}_{*}
(V)
\to (h^{\la,0})^{*}h^{\la,0}_{*}(V).
\] Using \cite[\cref{1-lemAlgebHTStack} \& \cref{1-CordfnBiFilAgLa0}]{Jiang2026Shla}, we see that \[
\gr^{-i}(h^{\la,0,++})^{*}h^{\la,0,++}_{*}
(V)\cong 0,i\in\ZZ_{\ge1}.
\] Consider the following commutative Cartesian diagram \[
\begin{tikzcd}
{\mX^{\la,0,++}_{K^{p}}\times [\mA^{1}/\GG_{m}]^{2}}\arrow[d,"h^{\la,0,++}"]
\arrow[r,"p_{1}"]
& \mX^{\la,0}_{K^{p}}\arrow[d,"h^{\la,0}"]\\
(\mX^{\tor,\la,0}_{K^{p}})^{\pdR/\CC_{p},++}_{\log}
\arrow[r,"p'"] & \bar{\mX}^{\tor,\la}_{K^{p}},
\end{tikzcd}
\] induced by \( [\mA^{1}/\GG_{m}]^{2}\to \AnSpec(\CC_{p}) \). Recall that \( V \) is viewed as a quasicoherent sheaf on \( {\mX^{\la,0,++}_{K^{p}}\times [\mA^{1}/\GG_{m}]^{2}} \) via \( p_{1}^{*} \).
Then 
\begin{align*}
(i_{0}\times i_{0})^{*}(h^{\la,0,++})^{*}h^{\la,0,++}_{*}
(V)&\cong p_{1,*}(h^{\la,0,++})^{*}h^{\la,0,++}_{*}
(V)
\\&\cong (h^{\la,0})^{*}p'_{*}h^{\la,0,++}_{*}(V)\cong (h^{\la,0})^{*}h^{\la,0}_{*}(V),
\end{align*} where the second isomorphism is given by the smooth base change, and the other two isomorphisms are given by 
\cite[\cref{1-corFilterdatA1modGm}]{Jiang2026Shla}.
\end{proof}
\begin{prop}\label{prop:const2-coincide}
Given \( (\mF_{1},\mF_{2},\delta) \)
and \( (\mF_{1}',\mF_{2}',\delta') \)
as in \cref{constr:Diff-op}. Then \cref{constr:Diff-op}
and \cref{constr:Diff-op-old} coincide.
\end{prop}
\begin{proof}
Let \( n,n' \) be as in \cref{constr:Diff-op-old}. Then \( \delta \)
and \( \delta' \)
also induce \begin{align*}
\delta^{+,\vee}&:h^{\sm,+}_{\natural}(\pi^{\sm}_{\GM})^{*}(\mF_{2}^{\vee})\to h^{\sm,+}_{\natural}(\pi^{\sm}_{\GM})^{*}(\mF_{1}^{\vee})\{n\},\\
\delta^{\prime,+,\vee}&:(h^{+}_{\Fl})_{\natural}(\mF_{2}^{\prime,\vee})\to (h^{+}_{\Fl})_{\natural}(\mF_{1}^{\prime,\vee})\{n'\},
\end{align*} and by the construction in \cite[\cref{1-propBGGFilSVsm}]{Jiang2026Shla}, \( \delta \)
is pulled back from \[
\delta^{+,\vee}:(h^{\std,+}_{\Fl})_{\natural}(\mF_{2}^{\vee})\{-n\}\to (h^{\std,+}_{\Fl})_{\natural}(\mF_{1}^{\vee})
\] along \( \pi^{\sm,\pdR,+}_{\GM}:(\mX^{\tor}_{K^{p}})_{\log}^{\pdR/\CC_{p},\sm,+}\to \Fl^{\std,\pdR/\CC_{p}}_{G,\mu}/G^{c} \), where \( h^{\std,+}_{\Fl} \)
is the natural map from \( \Fl^{\std}_{G,\mu}/G^{c}\times [\mA^{1}/\GG_{m}] \)
to \( \Fl^{\std,\pdR/\CC_{p},+}_{G,\mu}/G^{c} \). 

We can repeat \cref{constr:Diff-op}
replacing \( \delta^{\vee}\boxtimes \delta^{\prime,\vee} \)
with \( \delta^{+,\vee}\boxtimes \delta^{\prime,+,\vee} \) to obtain a morphism over \( (\mX^{\tor,\la,0}_{K^{p}})^{\pdR/\CC_{p},++}_{\log} \) \[
d^{\vee}:h^{\la,0,++}_{\natural}
(\pi^{\la,0}_{\GM}\times \pi^{\la,0}_{\HT})^{*}(\mF_{2}^{\vee}\boxtimes \mF_{2}^{\prime,\vee})
\to
h^{\la,0,++}_{\natural}
(\pi^{\la,0}_{\GM}\times \pi^{\la,0}_{\HT})^{*}(\mF_{1}^{\vee}\boxtimes \mF_{1}^{\prime,\vee})\{n,n'\},
\] whose pull-back along \( i_{0}\times i_{0} \) recovers \cref{constr:Diff-op} by smooth base change.

By applying \cref{lem:pullback-and-dual} to (\ref{eq:la-0-cartesian}), we see that \( d \) in \cref{constr:Diff-op-old} is obtained by applying \( \mbb{D}_{(\mX^{\tor,\la,0}_{K^{p}})^{\pdR/\CC_{p},++}_{\log}} \) to \( d^{\vee} \). Pulling back along \( i_{0}\times i_{0} \) and applying \cref{lem:pullback-and-dual} again, we see that \cref{constr:Diff-op}
and \cref{constr:Diff-op-old} coincide.
\end{proof}

\subsection{The case of modular curves}
\label{subsec:Fontaine-op-mod-curve}
In this subsection, we will explain a new proof of \cite[Theorem 6.2.6]{Pan2209.06II} using \cref{cor:BdR-lambda-E0D-pi-GMHT}, \cref{thm:Vin-GL2-Gr-t}
and \cref{thm:Vin-GL2-Fontaine=Diff}. 
This subsection will serve as a toy example of the argument in \S \ref{subsec:Fontaine-op-HMV} in the Hilbert case, and thus we will be sketchy with the argument.

For this subsection, we will take the Shimura datum to be \( (\GL_{2,\QQ},h) \)
that corresponds to modular curves. 
Recall the following diagram (\ref{eq:Alg-DR-Bdif})
\begin{equation}\label{eq:Vin-Big-Diagram}
\begin{tikzcd}
\Bdifcf^{+}(\mX^{\tor}_{K^{p}})^{\la}\arrow[d,"h^{\Bdif}"]\arrow[r,"\pi_{\GM\HT}^{\Bdif}"] &[-3em] \Vin^{*}_{G^{c},\mu,\cfi}|_{\hat{\mA}^{1}}/G^{c}\arrow[d,"\hat{h}_{0}"]
\\
\bar{\mX}^{\tor,\la}_{K^{p},\cfi}\arrow[r,"\bar{\pi}_{\GM\HT}"]\arrow[rd,"\pi_{\GM}^{\dR}\times \pi_{\HT}^{\dR}"'] & (\Vin^{*}_{G^{c},\mu,\cfi}|_{0})^{\pdR/\cfi}/G^{c}
\arrow[d,"\pi_{\GM}^{\loc,\dR}\times \pi_{\HT}^{\loc,\dR}"] &[-0.5em] \Vin^{*}_{G^{c},\mu,\cfi}|_{0}/G^{c}\arrow[d,"\pi^{\loc}_{\GM}\times \pi^{\loc,\dR}"] \arrow[l,"h_{0}"]
\arrow[lu,"i"]
\\
&
\Fl^{\std,\pdR/\cfi}_{G,\mu,\cfi}/G^{c}\times \Fl_{G,\mu,\cfi}^{\pdR/\cfi}
& \Fl^{\std}_{G,\mu,\cfi}/G^{c}\times \Fl_{G,\mu,\cfi}.\arrow[l,"h_{\Fl}^{\std}\times h_{\Fl}"]
\end{tikzcd}
\end{equation}
\begin{notation}
Given a weight \( (k,w) \), denote \[
\mO^{\la,(k,w)}_{K^{p},\cfi}:=\mbb{D}_{\bar{\mX}^{\tor,\la}_{K^{p},\cfi}}\left((\pi^{\dR}_{\GM})^{*}(h_{\Fl}^{\std})_{\natural}(\omega^{(-k,-w)}_{\Fl^{\std}_{G,\mu}})\otimes 
(\pi^{\dR}_{\HT})^{*}
(h_{\Fl})_{\natural}(\omega^{(-k,w)}_{\Fl_{G,\mu}})\right).
\]
Then by \cref{constr:Diff-op}, for \( k\ge 2 \), there is a differential operator \[
\theta^{k-1}_{\Fl^{\std}_{G,\mu}}\boxtimes
\theta^{k-1}_{\Fl_{G,\mu}}:\mO^{\la,(2-k,w)}_{K^{p},\cfi}\to \mO^{\la,(k,w)}_{K^{p},\cfi}.
\]
If we write \( d^{k-1}:=\theta^{k-1}_{\Fl^{\std}_{G,\mu}}\boxtimes 1 \)
and \( \bar{d}^{k-1}:=1\boxtimes
\theta^{k-1}_{\Fl_{G,\mu}} \), then \[
\theta^{k-1}_{\Fl^{\std}_{G,\mu}}\boxtimes
\theta^{k-1}_{\Fl_{G,\mu}}=d^{k-1}\circ \bar{d}^{k-1}=\bar{d}^{k-1}\circ d^{k-1}.
\]
\end{notation}
\begin{thm}\label{thm:Mod-curve-Fontaine=diff}
Let \( (k,w) \) be a weight with \( k\ge 2 \). 
The colimit system  \[ ``\varinjlim_{i\in\NN}" E_{0}^{\theta}(\BdifKpcf^{+,\la,\lambda^{(k,w)}}\cdot t^{-i}) \] is eventually constant. We denote the colimit as \( E_{0}^{\theta}(\BdifKpcf^{\la,\lambda^{(k,w)}}) \). Then \( E_{0}^{\theta}(\BdifKpcf^{\la,\lambda^{(k,w)}}) \) admits a descending \( t \)-adic filtration \( \Fil^{*}_{t^{\bullet}} \), where the only non-trivial graded pieces are \[
\gr^{\frac{m-w}{2}}_{t^{\bullet}}E_{0}^{\theta}(\BdifKpcf^{\la,\lambda^{(k,w)}})\cong \mO^{\la,(m,w)}_{K^{p},\cfi},\;m\in\{k,2-k\}.
\]

The action of \( \theta \) is zero on the graded pieces, and induces a monodromy operator \( N_{k}:\mO^{\la,(2-k,w)}_{K^{p},\cfi}\to \mO^{\la,(k,w)}_{K^{p},\cfi} \). Moreover, we have \[
N_{k}=c_{k}\cdot \theta^{k-1}_{\Fl^{\std}_{G,\mu}}\boxtimes \theta^{k-1}_{\Fl_{G,\mu}},
\] where \( c_{k}\in\QQ^{\times} \) is as in \cref{thm:Vin-GL2-Fontaine=Diff}.
\end{thm}
\begin{proof}
This follows immediately from \cref{cor:BdR-lambda-E0D-pi-GMHT}, \cref{thm:Vin-GL2-Gr-t}
and \cref{thm:Vin-GL2-Fontaine=Diff}.
\end{proof}


\subsection{The case of Hilbert modular varieties}
\label{subsec:Fontaine-op-HMV}

For this subsection, let \( F \) be a totally real field, and 
\( \Sigma_{\infty}:=\{\tau:F\hookrightarrow \RR\} \).
We will take the Shimura datum to be \( (G:=\Res_{F/\QQ}\GL_{2},h) \)
that corresponds to Hilbert modular varieties. More precisely, \[h:\mbb{S}\to G_{\RR}\cong \prod_{\tau\in \Sigma_{\infty}}\GL_{2,\RR},\;x+y\mi\mapsto (\begin{pmatrix}
	x&y\\-y&x
\end{pmatrix})_{\tau\in \Sigma_{\infty}}.\]
We will conjugate the Hodge cocharacter \( \mu \) to be \( t\mapsto (\mrm{diag}(t,1))_{\tau\in\Sigma_{\infty}} \). 
In this case, \( Z_{c}\cong \Ker(\mrm{Res}_{F/\QQ}\GG_{m}\xrightarrow{\Nm_{F/\QQ}}\GG_{m}) \), and \( G^{c}=G/Z_{c} \).

We take  \( \cf \) to be a finite extension of \( \QQ_{p} \) that contains every embedding of \( F \) into \( \bar{\QQ}_{p} \).
Along the isomorphism \( \iota:\CC\cong\bar{\QQ}_{p} \), \( \Sigma:=\{\tau:F\hookrightarrow \cf\}\cong \Sigma_{\infty}\).
We can still study the action of \( \theta \) on \( E_{0}^{\theta}(\BdifKpcf^{\la,\lambda^{\kw}}) \) as in \cref{thm:Mod-curve-Fontaine=diff}. For this,
we continue working with the diagram (\ref{eq:Vin-Big-Diagram}), where \( \Vin_{G^{c},\mu} \) is studied in \S \ref{subsec:eg-ResF-GL2}. We will make free use of notation in \cref{notation:Naive-dual}, \cref{dfn:weight-HMF} and \cref{notation:Weyl-action}. 
\begin{notation}\label{notation:Ola-kw}
Given a multiweight \( \kw \), denote \[
\mO^{\la,\kw}_{K^{p},\cfi}:=\mbb{D}_{\bar{\mX}^{\tor,\la}_{K^{p},\cfi}}\left((\pi^{\dR}_{\GM})^{*}(h_{\Fl}^{\std})_{\natural}(\omega^{-\kw}_{\Fl^{\std}_{G,\mu}})\otimes 
(\pi^{\dR}_{\HT})^{*}
(h_{\Fl})_{\natural}(\omega^{-(\Bbbk,-w)}_{\Fl_{G,\mu}})\right).
\]
Then by \cref{constr:Diff-op}, for any \( \tau\in\Sigma \), if \( k_{\tau}\ge 2 \), there is a differential operator \[
\theta^{k_{\tau}-1}_{\Fl^{\std}_{G,\mu}}\boxtimes
\theta^{k_{\tau}-1}_{\Fl_{G,\mu}}:\mO^{\la,(s_{\tau}\cdot \Bbbk,w)}_{K^{p},\cfi}\to \mO^{\la,\kw}_{K^{p},\cfi}.
\]
If we write \( d^{k_{\tau}-1}:=\theta^{k_{\tau}-1}_{\Fl^{\std}_{G,\mu}}\boxtimes 1 \)
and \( \bar{d}^{k_{\tau}-1}:=1\boxtimes
\theta^{k_{\tau}-1}_{\Fl_{G,\mu}} \), then \[
\theta^{k_{\tau}-1}_{\Fl^{\std}_{G,\mu}}\boxtimes
\theta^{k_{\tau}-1}_{\Fl_{G,\mu}}=d^{k_{\tau}-1}\circ \bar{d}^{k_{\tau}-1}=\bar{d}^{k_{\tau}-1}\circ d^{k_{\tau}-1}.
\]
\end{notation}
\begin{thm}\label{thm:HMV-Fontaine=diff}
Let \( \kw \) be a multiweight with \( k_{\tau}\ge 2 \) for all \( \tau\in\Sigma \). 

(1)
The colimit system  \[ ``\varinjlim_{i\in\NN}" E_{0}^{\theta}(\BdifKpcf^{+,\la,\lambda^{\kw}}\cdot t^{-i}) \] is eventually constant. We denote the colimit as \( E_{0}^{\theta}(\BdifKpcf^{\la,\lambda^{\kw}}) \). Then \( E_{0}^{\theta}(\BdifKpcf^{\la,\lambda^{\kw}}) \) admits a descending \( \ZZ^{\Sigma} \)-filtration \( \Fil^{*}_{(t_{\tau})_{\Sigma}} \) whose total filtration is the usual \( t \)-adic filtration, and for which the only \emph{non-trivial} graded pieces are 
\begin{equation}\label{eq:HT-grade-HM}
\gr^{\unl{j}_{T}}_{(t_{\tau})_{\Sigma}}E_{0}^{\theta}(\BdifKpcf^{\la,\lambda^{\kw}})\cong \mO^{\la,(s_{T}\cdot \Bbbk,w)}_{K^{p},\cfi},
\end{equation}
 for \( T\subset \Sigma \)
and \( \unl{j}_{T}:=(-\chi^{(s_{T}\cdot\Bbbk,w)}(d\mu_{\tau}))_{\tau} \). 

(2) There exist commuting operators \( \theta_{\tau_{0}}\in \End(E_{0}^{\theta}(\BdifKpcf^{\la,\lambda^{\kw}})) \) for \( \tau_{0}\in\Sigma \)
such that \[
\theta=\sum_{\tau\in\Sigma}\theta_{\tau}.
\]

(3) For \( j\in\ZZ \), denote \[ \Fil_{t_{\tau_{0}}}^{j}E_{0}^{\theta}(\BdifKpcf^{\la,\lambda^{\kw}}):=\varinjlim_{\unl{j}:j_{\tau}=j}\Fil_{(t_{\tau})_{\Sigma}}^{\unl{j}}E_{0}^{\theta}(\BdifKpcf^{\la,\lambda^{\kw}}) \]and denote by \( \gr_{t_{\tau_{0}}}^{j} \) the associated graded piceces, which still have induced \( \ZZ^{\Sigma\backslash \{\tau_{0}\}} \)-filtrations \( \Fil_{(t_{\tau})_{\Sigma\backslash\{\tau_{0}\}}}^{*} \), with \[
\gr_{(t_{\tau})_{\Sigma\backslash\{\tau_{0}\}}}^{(j_{\tau})_{\tau\in\Sigma\backslash\{\tau_{0}\}}}\gr_{t_{\tau_{0}}}^{j_{\tau_{0}}}E_{0}^{\theta}(\BdifKpcf^{\la,\lambda^{\kw}})
\cong \gr_{(t_{\tau})_{\Sigma}}^{(j_{\tau})_{\tau\in\Sigma}}E_{0}^{\theta}(\BdifKpcf^{\la,\lambda^{\kw}}).
\]

For \( \tau_{0}\in\Sigma \), the action of \( \theta_{\tau_{0}} \) preserves the multi-filtration \( \Fil^{*}_{(t_{\tau})_{\Sigma}} \), is zero on \( \gr^{*}_{t_{\tau_{0}}}(E_{0}^{\theta}(\BdifKpcf^{\la,\lambda^{\kw}}) ) \), and its action on \( E_{0}^{\theta}(\BdifKpcf^{\la,\lambda^{\kw}}) \)
induces a monodromy operator \[
\ti{N}_{\tau_{0},k_{\tau_{0}}}:\gr_{t_{\tau_{0}}}^{\frac{2-k_{\tau_{0}}-w}{2}}(E_{0}^{\theta}(\BdifKpcf^{\la,\lambda^{\kw}}))\to \gr_{t_{\tau_{0}}}^{\frac{k_{\tau_{0}}-w}{2}}(E_{0}^{\theta}(\BdifKpcf^{\la,\lambda^{\kw}})),
\] preserving \( \Fil^{(t_{\tau})_{\Sigma\backslash\{\tau_{0}\}}}_{*} \) and thus inducing 
\[
\ti{N}_{\tau_{0},k_{\tau_{0}}}:
\gr^{\unl{j}_{T\cup\{\tau_{0}\}}}_{(t_{\tau})_{\Sigma}}(E_{0}^{\theta}(\BdifKpcf^{\la,\lambda^{\kw}}))
\to 
\gr^{\unl{j}_{T}}_{(t_{\tau})_{\Sigma}}(E_{0}^{\theta}(\BdifKpcf^{\la,\lambda^{\kw}})),
\] for \( \tau_{0}\notin T\subset \Sigma \).

Moreover, up to the isomorphism (\ref{eq:HT-grade-HM}), we have
\begin{equation}\label{eq:HMF-Nk-theta-}
\ti{N}_{\tau_{0},k_{\tau_{0}}}=c_{k}\cdot \theta^{k_{\tau}-1}_{\Fl^{\std}_{G,\mu}}\boxtimes \theta^{k_{\tau}-1}_{\Fl_{G,\mu}},
\end{equation}
where \( c_{k_{\tau}}\in\QQ^{\times} \) is as in \cref{thm:Vin-GL2-Fontaine=Diff}, and \( \theta^{k_{\tau}-1}_{\Fl^{\std}_{G,\mu}}\boxtimes \theta^{k_{\tau}-1}_{\Fl_{G,\mu}} \) is as in \cref{notation:Ola-kw}.
\end{thm}
\begin{proof}
This follows immediately from \cref{cor:BdR-lambda-E0D-pi-GMHT}, and \cref{cor:from-G-to-Gc}.
\end{proof}

\section{Classicality for Hilbert modular forms}\label{sec:Classicality-Non-Parallel}
In this section, we extend the classicality theorem \cite[Theorem \ref{2-IntrothmClassicalityCompleteCoho}]{Jiang2025HMF} to non-parallel weight cases. 
We continue with the setting of \S \ref{subsec:Fontaine-op-HMV}.
Let \( d:=[F:\QQ] \).
	We fix a neat tame level \(K^{p}\subset G(\mA_{f}^{p})\).
	Let \(S\) be a finite set of places of \(\QQ\),
	including \(p\) and \(\infty\)
	and all the ramified places of \(K^{p}\) or \(F \).
	Denote by \(\TT^{S}\) the spherical Hecke algebra over \(\ZZ_{p}\) generated by the spherical Hecke operators of \(G\) at finite places away from \(S\). 
\begin{dfn}\label{dfnAssociatedGaloisRep}
Let \(\rho:\Gal_{F,S}\to \GL_{2}(L)\) be a continuous representation.
 We denote by \(\chi=\chi_{\rho}:\TT^{S}[1/p]\to L\)  the associated character determined by the Eichler-Shimura relation 
 \begin{align}\label{alignEicherShimuraDetermin}
 \det(t-\Fr_{\mfk{l}})=t^{2}-
T_{\mfk{l}}\cdot t+\Nm(\mfk{l})S_{\mfk{l}},
 \end{align}
 that is, \(\chi(S_{\mfk{l}}):=\mrm{Nm}(\mfk{l})^{-1}\det(\rho(\Fr_{\mfk{l}}))\) and \(\chi(T_{\mfk{l}}):=\mrm{Tr}(\Fr_{\mfk{l}})\).
\end{dfn}
\begin{thm}\label{thm:ClassicalityCompleteCoho}
	Let \(\rho:\Gal_{F,S}\to \GL_{2}(L)\) be a continuous representation,
	with the associated character \( \chi \) determined by the Eichler-Shimura relation (\ref{alignEicherShimuraDetermin}).

	We assume that 
	
	(1) \(\ti{H}^{d}(K^{p},L)[\chi]:=\Hom_{\TT^{S}}(\chi,\ti{H}^{d}(K^{p},L))\ne 0\);

	(2) \(\rho\) is absolutely irreducible;
	
	(3) \( \rho \) is de Rham of regular Hodge-Tate weights; that is, there exists a multiweight \( \kw \) with \( k_{\tau}\ge 2,\forall \tau\in\Sigma \), such that \(\forall\tau\in \Sigma\), \(\rho\) is \(\tau\)-de Rham of
	\(\tau\)-Hodge-Tate weights \(\frac{k_{\tau}-w}{2}\) and \(\frac{2-k_{\tau}-w}{2}\).

	Then the character \(\chi\) is classical,
	in the sense that there exists a Hilbert modular eigenform \(f\) of weight \( \kw \),
	such that \(\TT^{S}\) acts on \(f\) via \(\chi\). In other words, 
	\(\rho\cong \rho_{f}\).
\end{thm}
The key ingredient is the following classicality theorem proven in \cite{Jiang2025HMF}:
\begin{dfn}\label{dfn:FontainComplex}
	We fix a multiweight \(\kw\) with \( k_{\tau}\ge 2,\forall\tau\in\Sigma \). We fix a total order \( < \) on \( \Sigma \).
	Using \(\ti{N}_{I_{0},\tau}\), 
	we can define the \emph{Bernstein-Gelfand-Gelfand-Fontaine complex} by 
	\begin{align*}
	&\mrm{Font}(\BdifKpcf^{\la,\lambda^{\kw}}):=\\&\left[
	\mO^{\la,\kw}_{K^{p},\cfi}\to\bigoplus_{\tau\in \Sigma}\mO^{\la,(s_{\tau}\cdot\Bbbk,w)}_{K^{p},\cfi}\to \cdots \bigoplus_{\#(I_{0})=i,I_{0}\subset \Sigma}\mO^{\la,(s_{I_{0}}\cdot\Bbbk,w)}_{K^{p},\cfi}\to \cdots \to \mO^{\la,(s_{\Sigma}\cdot\Bbbk,w)}_{K^{p},\cfi}
\right],
	\end{align*}
where the differential from the \( I_{0} \)-th term to the \( I_{1} \)-th term is  given by \(\epsilon(I_{0},\tau)\ti{N}_{\tau,k_{\tau}}\) (for \( \ti{N}_{\tau,k_{\tau}} \) as in \cref{thm:HMV-Fontaine=diff} and \( \epsilon(I_{0},\tau):=(-1)^{\#(I_{0}\cap \Sigma_{<\tau})} \))
for any \(I_{0}\subsetneq I_{1}=I_{0}\cup\{\tau\} \subset \Sigma\) and zero otherwise. 
Since \( \theta_{\tau} \) in \cref{thm:HMV-Fontaine=diff} commute with each other,
	\( \mrm{Font}(\BdifKpcf^{\la,\lambda^{\kw}}) \) indeed forms a complex.
\end{dfn}
\begin{thm}\label{thm:classicality-Fontaine-complex}
We continue with the setting of \cref{dfn:FontainComplex}. Then \( R\Gamma(\bar{\mX}^{\tor,\la}_{K^{p},\cfi},\mrm{Font}(\BdifKpcf^{\la,\lambda^{\kw}})) \)
is classical of weight \( \kw \) as a Hecke module (\cite[\cref{2-dfnClassicalHecke}]{Jiang2025HMF}). In other words, for any prime ideal \( \mfk{a}\subset \TT^{S} \), if \[
R\Gamma(\bar{\mX}^{\tor,\la}_{K^{p},\cfi},\mrm{Font}(\BdifKpcf^{\la,\lambda^{\kw}}))_{\mfk{a}}\ne 0,
\] then \( \mfk{a} \) is associated to a Hilbert modular eigenform of weight \( \kw \).
\end{thm}
\begin{proof}
The proof is the same as that of  \cite[\cref{2-thmFontainComplexClassical}]{Jiang2025HMF}. Let us give a sketch here. 
Let \[
\nabla_{K^{p},\cfi}^{\kw,(\Bbbk',w')}:=\mbb{D}_{\bar{\mX}^{\tor,\la}_{K^{p},\cfi}}\left((\pi^{\dR}_{\GM})^{*}(h_{\Fl}^{\std})_{\natural}(\omega^{-\kw}_{\Fl^{\std}_{G,\mu}})\otimes 
(\pi^{\dR}_{\HT})^{*}
(h_{\Fl})_{\natural}(\omega^{-(\Bbbk',-w')}_{\Fl_{G,\mu}})\right).
\]
For \( I\subset \Sigma \), if \( k_{\tau}\le 0 \) for \( \tau\in I \), and \( k'_{\tau}\le 0 \)
for \( \tau\in I^{c} \), then we have 
\begin{align*}
\theta^{k_{\tau}-1}_{\Fl^{\std}_{G,\mu}}\boxtimes 1&:\nabla_{K^{p},\cfi}^{\kw,(\Bbbk',w')}\to \nabla_{K^{p},\cfi}^{(s_{\tau}\cdot \Bbbk,w),(\Bbbk',w')},\; \tau\in I;
\\
1\boxtimes\theta^{k_{\tau}-1}_{\Fl_{G,\mu}}&:\nabla_{K^{p},\cfi}^{\kw,(\Bbbk',w')}\to \nabla_{K^{p},\cfi}^{\kw,(s_{\tau}\cdot\Bbbk',w')},\; \tau\in I^{c}.
\end{align*}
We then define the partial de Rham complex \( \dR_{I,I^{c}}(\nabla_{K^{p},\cfi}^{\kw,(\Bbbk',w')}) \) to be the associated ``Koszul complex'' of which the \( m \)-th term is \[
\bigoplus_{I_{0}\subset \Sigma,\#(I_{0})=m}\nabla_{K^{p},\cfi}^{(s_{I\cap I_{0}}\cdot\Bbbk,w),(s_{I^{c}\cap I_{0}}\cdot\Bbbk',w')},
\] 
and the differential from the \( I_{0} \)-th term to the \( (I_{0}\cup\{\tau\}) \)-th term is given by \(\epsilon(I_{0},\tau)\theta^{k_{\tau}-1}_{\Fl^{\std}_{G,\mu}}\boxtimes 1\) (if \( \tau\in I \)) or \(\epsilon(I_{0},\tau)1\boxtimes\theta^{k_{\tau}'-1}_{\Fl_{G,\mu}} \) (if \( \tau\in I^{c} \)).

By \cref{rmk:dif-op-L-vs-Cp} and \cref{prop:const2-coincide}, \[
\dR_{I,I^{c}}(\nabla_{K^{p},\cfi}^{\kw,(\Bbbk',w')})\otimes_{\cfi}\CC_{p}\cong \mfk{s}_{*}\dR_{I,I^{c}}(\nabla_{\SGTKp{\emptyset}}^{\kw,(\Bbbk',w')}),
\] where \( \mfk{s} \) is as in \cref{rmk:dif-op-L-vs-Cp} and \( \dR_{I,I^{c}}(\nabla_{\SGTKp{\emptyset}}^{\kw,(\Bbbk',w')}) \) is as in \cite[\cref{2-dfnDRComplexNotation}]{Jiang2025HMF} (3). Thus by \cite[\cref{2-thmClassicalityCohomoDR}]{Jiang2025HMF}, we see that \[
R\Gamma(\bar{\mX}^{\tor,\la}_{K^{p},\cfi},\dR_{I,I^{c}}(\nabla_{K^{p},\cfi}^{\kw,(\Bbbk',w')}))
\] is classical of weight \( \kw \) as a Hecke module.

As in the proof of \cite[\cref{2-thmFontainComplexClassical}]{Jiang2025HMF},
using (\ref{eq:HMF-Nk-theta-}), the complex \( \mrm{Font}(\BdifKpcf^{\la,\lambda^{\kw}}) \)
is filtered by complexes of the form \[ \dR_{I,I^{c}}(\nabla_{K^{p},\cfi}^{(s_{J}\cdot \Bbbk,w),(s_{J}\cdot\Bbbk,w)}), \]
and thus \( R\Gamma(\bar{\mX}^{\tor,\la}_{K^{p},\cfi},\mrm{Font}(\BdifKpcf^{\la,\lambda^{\kw}})) \)
is classical of weight \( \kw \) as a Hecke module.
\end{proof}
\begin{proof}[Proof of \cref{thm:ClassicalityCompleteCoho}]
The first half of the proof is identical to that of \cite[\cref{2-thmClassicalityCompleteCoho}]{Jiang2025HMF}. We give a sketch here for the convenience of the readers.  
The case where \(\rho\) admits CM can be reduced to global class field theory. By \cite[\cref{2-lemCMorStrongIrr}]{Jiang2025HMF}, we assume without loss of generality that \(\rho\) is strongly irreducible. 

Let \( \p_{\chi}:=\Ker(\chi:\TT^{S}[1/p]\to \cf) \).
By \cite[Theorem 1.1]{JiangPan2025Plectic},
we have
\begin{equation}\label{eq:concentrated-deg-d}
R\Gamma(K^{p},\cf)^{R-\la}_{\p_{\chi}}\cong \ti{H}^{d}(K^{p},\cf)^{\la}_{\p_{\chi}}[-d],
\end{equation}
and 
\begin{equation}\label{eq:plec-action}
\ti{H}^{d}(K^{p},L)^{\la}[\chi]\cong \oInd_{F}^{\QQ}\rho\otimes \Pi^{\la}(\rho),
\end{equation}
where \( \oInd_{F}^{\QQ}\rho \) denotes the tensor induction of \( \rho \) from \( \Gal_{F} \)
to \( \Gal_{\QQ} \), \( \Pi^{\la}(\rho) \) is a locally analytic representation of \( G(\QQ_{p}) \), and \( -[\chi] \) denotes the subspace on which the \( \TT^{S} \)-action factors through \( \chi \). 

[\emph{\textbf{Step 1}: Vanishing of higher \( \Ext^{i}_{Z(\fg^{c})}(\lambda^{\kw},-) \)}]
By (\ref{eq:concentrated-deg-d}), we have
\[ \RHom_{Z(\fg^{c})}(\lambda^{\kw},R\Gamma(K^{p},\cf)^{R-\la}_{\p_{\chi}})\in D^{\ge d}(\cf), \]
with 
\begin{equation}\label{eq:fix-infi-car}
H^{d}(\RHom_{Z(\fg^{c})}(\lambda^{\kw},R\Gamma(K^{p},\cf)^{R-\la}_{\p_{\chi}}))\cong \Hom_{Z(\fg^{c})}(\lambda^{\kw},\ti{H}^{d}(K^{p},\cf)^{\la}_{\p_{\chi}}). 
\end{equation}
On the other hand, by \cite[\cref{2-propInfiniteCharToHorizontal}]{Jiang2025HMF} and the primitive comparison (\cite[Theorem 6.2.6]{Juan2022.09locallyShi}), we have the following Hodge-Tate decomposition 
\begin{equation}\label{eq:HT-decomp}
\RHom_{Z(\fg^{c})}(\lambda^{\kw},R\Gamma(K^{p},\cf)^{R-\la}_{\p_{\chi}}))\otimes_{\QQ_{p}}\CC_{p}
\cong \bigoplus_{I\subset \Sigma}R\Gamma_{\an}(\SGTKpTo{\emptyset},\mO^{\la,(s_{I}\cdot\Bbbk,w)}_{\SGTKp{\emptyset}})_{\p_{\chi}}\otimes\cf,
\end{equation}
and the RHS lies in \( D^{[0,d]}(\cf) \) by Grothendieck vanishing theorem. 
Therefore, we conclude that both sides are concentrated in degree \( d \). 

[\emph{\textbf{Step 2}: Passing from \( \CC_{p} \) to \( \cfi \)}]
By the argument in \cref{rmk:dif-op-L-vs-Cp} and \cite[\cref{2-lemHodgeTateCompareOnTcase}]{Jiang2025HMF},
we also have that
\begin{equation}\label{eq:pass-Cp-Li}
R\Gamma_{\an}(\SGTKpTo{\emptyset},\mO^{\la,(s_{I}\cdot\Bbbk,w)}_{\SGTKp{\emptyset}})\cong 
R\Gamma(\bar{\mX}^{\tor,\la}_{K^{p},\cfi},\mO^{\la,(s_{I}\cdot\Bbbk,w)}_{K^{p},\cfi})\otimes_{\cfi}\CC_{p}(\chi^{\kw}(d\mu)),
\end{equation}
where \( (\chi^{\kw}(d\mu)) \) denotes the Tate twist.
Note that \( \SGTKpTo{\emptyset} \) in \cite{Jiang2025HMF} is denoted as \( \mX_{K^{p},\proket}^{\tor} \) in this paper.
In particular, we also have \[
R\Gamma(\bar{\mX}^{\tor,\la}_{K^{p},\cfi},\mO^{\la,(s_{I}\cdot\Bbbk,w)}_{K^{p},\cfi})_{\p_{\chi}}
\cong H^{d}(\bar{\mX}^{\tor,\la}_{K^{p},\cfi},\mO^{\la,(s_{I}\cdot\Bbbk,w)}_{K^{p},\cfi})_{\p_{\chi}}[-d].
\]

[\emph{\textbf{Step 3}: Partial filtration and partial Fontaine operators}]
By \cref{thm:HMV-Fontaine=diff}, we know that 
\[\ti{M}_{\chi}:=\varinjlim_{i\in\NN} 
R\Gamma(\bar{\mX}^{\tor,\la}_{K^{p},\cfi},E_{0}^{\theta}(\BdifKpcf^{+,\la,\lambda^{\kw}}\cdot t^{-i}))_{\p_{\chi}}[d]
\] is concentrated in degree \( 0 \), and admits a descending \( \ZZ^{\Sigma} \)-filtration \( \Fil^{*}_{(t_{\tau})_{\Sigma}}\ti{M}_{\chi} \) whose total filtration is the usual \( t \)-adic filtration, and for which the only \emph{non-trivial} graded pieces are 
\begin{equation}\label{eq:HT-grade-HM-1}
\gr^{\unl{j}_{T}}_{(t_{\tau})_{\Sigma}}\ti{M}_{\chi}\cong H^{d}(\bar{\mX}^{\tor,\la}_{K^{p},\cfi},\mO^{\la,(s_{T}\cdot\Bbbk,w)}_{K^{p},\cfi})_{\p_{\chi}}.
\end{equation}
 for \( T\subset \Sigma \). Moreover, there exist commuting operators \( \theta_{\tau_{0}} \) acting on \( \ti{M}_{\chi} \) for \( \tau_{0}\in\Sigma \), such that \( \theta=\sum_{\tau\in\Sigma}\theta_{\tau} \), and 
the action of \( \theta_{\tau_{0}} \) preserves the multi-filtration \( \Fil^{*}_{(t_{\tau})_{\Sigma}} \ti{M}_{\chi}\), is zero on \( \gr^{*}_{t_{\tau_{0}}}\ti{M}_{\chi}\) (that is defined as in \cref{thm:HMV-Fontaine=diff}), and induces \[
\ti{N}_{\tau_{0},k_{\tau_{0}}}:
\gr^{\unl{j}_{T\cup\{\tau_{0}\}}}_{(t_{\tau})_{\Sigma}}\cong H^{d}(\bar{\mX}^{\tor,\la}_{K^{p},\cfi},\mO^{\la,(s_{T\cup\{\tau_{0}\}}\cdot\Bbbk,w)}_{K^{p},\cfi})_{\p_{\chi}}
\to 
\gr^{\unl{j}_{T}}_{(t_{\tau})_{\Sigma}}\cong H^{d}(\bar{\mX}^{\tor,\la}_{K^{p},\cfi},\mO^{\la,(s_{T}\cdot\Bbbk,w)}_{K^{p},\cfi})_{\p_{\chi}},
\] for \( \tau_{0}\notin T\subset \Sigma \).

We will prove by contradiction. Let us assume that \( \chi \) is not associated to a Hilbert modular eigenform of weight \( \kw \). Then by \cref{thm:classicality-Fontaine-complex}, the map 
\begin{equation}\label{eq:first-term-Fontaine}
(\ti{N}_{\tau,k_{\tau}})_{\tau\in\Sigma}:\gr^{\unl{j}_{\Sigma}}_{(t_{\tau})_{\Sigma}}\ti{M}_{\chi}
\to \bigoplus_{\tau\in\Sigma}\gr^{\unl{j}_{\Sigma\backslash \{\tau\}}}_{(t_{\tau})_{\Sigma}}\ti{M}_{\chi}
\end{equation}
is injective.

By \cref{thm:cplt-coho-primitive-BdR}, \[
\Darithcf(R\Gamma(K^{p},\QQ_{p})^{R-\la}\hatotimes_{\QQ_{p}}B_{\dR}^{+})^{\wedge}_{t}\cong R\Gamma(\Bdifcf^{+}(\mX^{\tor}_{K^{p}})^{\la}),
\] and thus by \emph{Step 1}, 
\begin{equation}\label{eq:primite-lambda-chi}
\varinjlim_{i\in\NN}E_{0}^{\theta}(\Darithcf(\ti{H}^{d}(K^{p},\cf)^{\la,\lambda^{\kw}}_{\p_{\chi}}\hatotimes_{\QQ_{p}}B_{\dR}^{+})^{\wedge}_{t}\cdot t^{-i})\cong \ti{M}_{\chi}\otimes_{\QQ_{p}}\cf,
\end{equation}
where here \( (-)^{\lambda^{\kw}} \) denotes the subspace where \( Z(\fg^{c}) \) acts through \( \lambda^{\kw} \). 

[\emph{\textbf{Step 4}: Passing to \( -[\chi] \)}]
By \cite[\cref{2-propFixInfCharc}]{Jiang2025HMF}, \( Z(\fg^{c}) \)
acts on \( \ti{H}^{d}(K^{p},\cf)^{\la}[\chi] \) via \( \lambda^{\kw} \). Thus by (\ref{eq:plec-action}), \[
\ti{H}^{d}(K^{p},\cf)^{\la,\lambda^{\kw}}[\chi]\cong \oIndFQ\rho\otimes \Pi^{\la}(\rho)\ne 0.
\]
By (\ref{eq:HT-decomp})
and (\ref{eq:pass-Cp-Li}), we have the following Hodge-Tate decomposition \begin{equation}\label{eq:fix-chi-HT}
\Darithcf(\oIndFQ\rho\otimes\CC_{p})\otimes \Pi^{\la}(\rho)
\cong \bigoplus_{I\subset\Sigma}(\gr^{\unl{j}_{I}}\ti{M}_{\chi}\otimes_{\QQ_{p}}\cf)[\chi]
(\chi^{(s_{I}\cdot\Bbbk,w)}(d\mu)),  
\end{equation}
and combining with (\ref{eq:primite-lambda-chi}), \begin{equation}\label{eq:chi-compo-E0-Darith}
\varinjlim_{i\in\NN}E_{0}^{\theta}(\Darithcf(\oIndFQ\rho\hatotimes_{\QQ_{p}}B_{\dR}^{+})^{\wedge}_{t}\cdot t^{-i})\otimes \Pi^{\la}(\rho)\cong (\ti{M}_{\chi}\otimes_{\QQ_{p}}\cf)[\chi].
\end{equation}
In particular, \[
(\gr^{\unl{j}_{\Sigma}}\ti{M}_{\chi}\otimes_{\QQ_{p}}\cf)[\chi]\ne 0.
\] 
The de Rhamness of \( \rho \) implies that of \( \oIndFQ\rho \), and thus by \cref{lem:Fontaine0de-Rham} the action of \( \theta \) is zero on \( (\ti{M}_{\chi}\otimes_{\QQ_{p}}\cf)[\chi] \).

[\emph{\textbf{Step 5}: Conclusion}]
Through the inclusion \[
(\ti{M}_{\chi}\otimes_{\QQ_{p}}\cf)[\chi] \subset \ti{M}_{\chi}\otimes_{\QQ_{p}}\cf,
\] the \( \ZZ^{\Sigma} \)-filtration on \( \ti{M}_{\chi} \) induces a \( \ZZ^{\Sigma} \)-filtration \( \Fil^{*}_{(t_{\tau})_{\Sigma}} \) on \( (\ti{M}_{\chi}\otimes_{\QQ_{p}}\cf)[\chi]  \), such that 
\begin{equation}\label{eq:induce-gr-inj}
\gr^{*}_{(t_{\tau})_{\Sigma}}((\ti{M}_{\chi}\otimes_{\QQ_{p}}\cf)[\chi] )
\subset (\gr^{*}_{(t_{\tau})_{\Sigma}}\ti{M}_{\chi}\otimes \cf)[\chi]. 
\end{equation}
However, (\ref{eq:fix-chi-HT}) and (\ref{eq:chi-compo-E0-Darith}) force (\ref{eq:induce-gr-inj})
to be an identity by the discussion in \cite[\cref{2-notationArithSenHT}]{Jiang2025HMF}. Now the action of \( \theta=\sum_{\tau}\theta_{\tau} \) on \( (\ti{M}_{\chi}\otimes_{\QQ_{p}} \cf )[\chi]\) preserves the \( \ZZ^{\Sigma} \)-filtration, and is zero. 

For each \( \tau_{0} \), 
let \( \Fil^{*}_{(t_{\tau})_{\Sigma\backslash\{\tau_{0}\}}} \) be the \( \ZZ^{\Sigma\backslash\{\tau_{0}\}} \)-filtration obtained from 
\( \Fil^{*}_{(t_{\tau})_{\Sigma}} \) by forgetting the \( \tau_{0} \)-component. 
We can pass to the quotient \( \gr_{(t_{\tau})_{\Sigma\backslash\{\tau_{0}\}}}^{\unl{j}_{\Sigma}|_{\Sigma\backslash\{\tau_{0}\}}} \), where \( \unl{j}_{\Sigma}|_{\Sigma\backslash\{\tau_{0}}\} \)
is the restriction of \( \unl{j}_{\Sigma} \) to \( \ZZ\backslash \{\tau_{0}\} \). Then we have a short exact sequence \[
0\to (\gr^{\unl{j}_{\Sigma\backslash\{\tau_{0}\}}}_{(t_{\tau})_{\Sigma}}\ti{M}_{\chi}\otimes\cf) [\chi]\to 
\gr_{(t_{\tau})_{\Sigma\backslash\{\tau_{0}\}}}^{\unl{j}_{\Sigma}|_{\Sigma\backslash\{\tau_{0}\}}}\left((\ti{M}_{\chi}\otimes\cf)[\chi]\right)
\to (\gr^{\unl{j}_{\Sigma}}_{(t_{\tau})_{\Sigma}}\ti{M}_{\chi}\otimes\cf)[\chi]
\to 0. 
\] As discussed in \emph{Step 3}, the action of \( \theta_{\tau} \) on \( \gr_{(t_{\tau})_{\Sigma\backslash\{\tau_{0}\}}}^{\unl{j}_{\Sigma}|_{\Sigma\backslash\{\tau_{0}\}}}\left((\ti{M}_{\chi}\otimes\cf)[\chi]\right) \) is zero if \( \tau\ne \tau_{0} \), and the action of \( \theta_{\tau_{0}} \) is induced by \( \ti{N}_{\tau_{0},k_{\tau_{0}}} \). Therefore, \( \theta=0 \) implies that \[ \ti{N}_{\tau_{0},k_{\tau_{0}}}:(\gr^{\unl{j}_{\Sigma}}_{(t_{\tau})_{\Sigma}}\ti{M}_{\chi}\otimes\cf)[\chi]\to (\gr^{\unl{j}_{\Sigma}}_{(t_{\tau})_{\Sigma}}\ti{M}_{\chi}\otimes\cf)[\chi] \]
is zero for any \( \tau_{0}\in\Sigma \). This contradicts the injectivity of (\ref{eq:first-term-Fontaine}).
\end{proof}
As in \cite{Jiang2025HMF}, using \cref{thm:ClassicalityCompleteCoho}, we can prove the following new cases of Langlands-Clozel-Fontaine-Mazur conjecture, removing the parallel weight condition in \cite[\cref{2-thmConjecFontaineMazur}]{Jiang2025HMF}:
\begin{thm}\label{thmConjecFontaineMazur}
Let \(F\) be totally real number field of degree \(d\). Assume that \(p\) is inert in \(F\).
Let \(L/\QQ_{p}\) be a finite extension with residue field \(\FF\), and let
\(\rho:\Gal_{F,S}\to \GL_{2}(L)\) be a continuous representation.
Let \(\bar{\rho}:\Gal_{F}\to \GL_{2}(\FF)\) be the reduction of \(\rho\).

Assume that 

(1) \(\bar{\rho}\) is residual modular, that is, there exists a Hilbert modular eigenform \(f_{0}\)
such that \(\bar{\rho}\cong \bar{\rho}_{f_{0}}\).

(2) \(\bar{\rho}|_{\Gal_{F(\sqrt[p]{1})}}\)
is absolutely irreducible.

(3) \(\bar{\rho}|_{\Gal_{F_{p}}}\) is semisimple generic in the sense of \cite[p.100 (i)]{BreuilHerzigHuMorraSchraen2023GKdim}.

(4) For any \(w\in S\backslash \{p\}\),
the universal framed deformation ring of \(\bar{\rho}|_{\Gal_{F_{w}}}\) over \(W(\FF)\)
is formally smooth over \(W(\FF)\).

(5) Let \(\psi\) denote the Teichmüller lift of \(\chi_{\mrm{cycl}}\det(\bar{\rho})\), and we 
assume that \(\det(\rho)=\psi \chi_{\mrm{cycl}}^{-1}\).

Let \(\chi:\TT^{S}\to \mO_{L}\) be the character determined by \(\rho\) via (\ref{alignEicherShimuraDetermin}).

Then there exists a neat open compact subgroup \(K^{p}\subset G_{\emptyset}(\mA_{f}^{p})\), and a finite set \(S'\supset S\) such that \(K^{p}\) is unramified away from \(S'\), and \[\ti{H}^{d}(K^{p},L)[\chi|_{\TT^{S'}}]\ne 0.
\]

If we assume in addition that 

(6) \(\rho\) is de Rham of regular weights.

Then there exists a Hilbert modular eigenforms \(f\)
such that \(\rho\cong \rho_{f}\).
\end{thm}
\begin{proof}
The proof is identical to that of \cite[\cref{2-thmConjecFontaineMazur}]{Jiang2025HMF}, replacing \cite[\cref{2-thmClassicalityCompleteCoho}]{Jiang2025HMF}
with \cref{thm:ClassicalityCompleteCoho}.
\end{proof}
\printbibliography
\end{document}